\documentclass[11pt,reqno]{amsart}

\usepackage{amsmath,amssymb,amsthm,mathtools}
\usepackage[margin=1.1in]{geometry}
\usepackage{hyperref}
\hypersetup{hidelinks}
\usepackage{enumerate,color}
\usepackage{tikz}
\usetikzlibrary{arrows.meta}
\usepackage{float}
\usepackage{placeins}

\makeatletter
\renewcommand{\@captionfont}{\normalfont\small}
\makeatother

\numberwithin{equation}{section}

\newtheorem{theorem}{Theorem}[section]
\newtheorem{proposition}[theorem]{Proposition}
\newtheorem{lemma}[theorem]{Lemma}
\newtheorem{corollary}[theorem]{Corollary}
\newtheorem*{proposition*}{Proposition}
\newtheorem*{corollary*}{Corollary}
\theoremstyle{definition}

\theoremstyle{remark}
\newtheorem{remark}[theorem]{Remark}

\newcommand{\R}{\mathbb{R}}
\newcommand{\N}{\mathbb{N}}
\newcommand{\sgn}{\operatorname{sgn}}
\newcommand{\ovU}{\overline{U}}

\title[Self-similar cusp formation]{Self-similar cusp formation in a large class of nonlocal transport equations}

\author[Y. Kim, B. Kwon, and W. Shim]{Yunjoo Kim, Bongsuk Kwon, and Wanyong Shim}
\address{Department of Mathematical Sciences, Ulsan National Institute of Science and Technology, Ulsan, 44919, Korea}
\email{gomuli3@unist.ac.kr}

\address{Department of Mathematical Sciences, Ulsan National Institute of Science and Technology, Ulsan, 44919, Korea}
\email{bkwon@unist.ac.kr}

\address{Department of Mathematical Sciences, Korea Advanced Institute of Science and Technology, Daejeon, 34141, Korea}
\email{wyshim25@kaist.ac.kr}

\date{}


\thanks{\textbf{Acknowledgment.}
B.K. was supported by the National Research Foundation of Korea (NRF) grant funded by the Korean government (MSIT) (00560003).}

\begin{document}

\begin{abstract}
We prove finite-time gradient blow-up with sharp \(C^{(5-b)/(b+3)}\) H\"older regularity for a class of nonlocal transport equations extending the \(b\)-family, for all real \(1<b<5\). This class includes the Camassa--Holm and Degasperis--Procesi equations and, after rescaling, the generalized hyperelastic-rod wave equation. In the blow-up regime, the nonlocal equations reduce at leading order to the generalized Hunter--Saxton (gHS) equation with parameter \(\ell=2/(b-1)\), whose self-similar scaling and corresponding profiles depend on \(\ell\). Consequently, standard approaches tailored to a fixed self-similar profile cannot be applied directly. To handle this parameter dependence, we develop a general framework for self-similar singularity formation in evolution equations whose leading-order dynamics are governed by the gHS family. The key ingredient is a family of profile inequalities established for all real \(\ell>1/2\), which we combine with a generalized maximum principle to close a parameter-dependent bootstrap argument. We identify model-specific hypotheses adapted to this framework, under which this argument yields the required stability estimates. This leads to a general theorem on singularity formation, which we apply to the entire class of nonlocal equations, thereby establishing finite-time gradient blow-up and sharp H\"older regularity. The theorem can also be applied to other evolution equations once the corresponding model-specific hypotheses are verified, without repeating the profile analysis or the common stability argument.
\end{abstract}

\maketitle

\tableofcontents

\section{Introduction}\label{sec:intro}

Finite-time gradient blow-up is a fundamental mechanism for singularity formation in a broad class of transport-type evolution equations \cite{ConstantinEscher1998,EscherYin2008}. A classical characteristic argument can often establish the occurrence of such blow-up by showing that a negative slope becomes unbounded in finite time, but by itself it provides little information about the detailed behavior of the solution near the singularity. Buckmaster, Shkoller, and Vicol pioneered a different approach in their work on shock formation for the compressible Euler equations \cite{BSV2D,BSV3D}: they introduced dynamically rescaled coordinates and a logarithmic time variable, and converted finite-time blow-up into a long-time stability problem for a stationary self-similar profile. The blow-up time, location, and amplitude are tracked by modulation parameters, while a nonlinear bootstrap argument controls the perturbation around the profile. This approach thereby provides a quantitative description of how the singularity forms in space and time.

Variants of this strategy have since been used for nonlinear transport-type equations and systems whose leading-order blow-up dynamics are governed by the inviscid Burgers equation. Examples include the Euler--Poisson systems with and without pressure \cite{BKK-ion,BKK-delta}, the Burgers--Hilbert equation \cite{YangBH2021}, the fractal Burgers equation \cite{ChickeringMorenoVasquezPandya2023}, and the generalized Fornberg--Whitham equation \cite{SautWang2025}; see also \cite{OhPasqualotto} for a broader class of dispersive and dissipative perturbations of the Burgers equation. In the stable blow-up regimes considered in these works, the leading-order Burgers dynamics give rise to the prototypical cubic-root cusp with sharp \(C^{1/3}\) H\"older regularity at the blow-up time.

More recent studies have applied related self-similar approaches to nonlocal wave-breaking problems in which the leading-order dynamics differ from those of the inviscid Burgers equation \cite{KKY-CH,KKS-rSV,KKKS-Novikov}. In the Camassa--Holm equation, the transport part has the form of the Burgers equation, but the nonlocal pressure contributes at leading order after differentiation. The resulting blow-up dynamics are governed by the Hunter--Saxton equation, giving rise to a cusp with sharp \(C^{3/5}\) H\"older regularity \cite{KKY-CH}. The same \(C^{3/5}\) cusp formation, governed by Hunter--Saxton self-similar profiles, was subsequently established for the regularized Saint--Venant equations \cite{KKS-rSV}. For the Novikov equation, a Galilean-type change of variables around a nonzero background reveals a Hunter--Saxton-type leading equation with a different coefficient of \(u_x^2\). This difference changes the scaling and far-field behavior of the self-similar profiles, yielding stable cusp formation with sharp \(C^{7/9}\) H\"older regularity \cite{KKKS-Novikov}.

The present paper is motivated by the analyses in \cite{KKY-CH,KKKS-Novikov}, which follow the same sequence: identification of the leading equation, construction and analysis of its self-similar profiles, stability estimates in modulated variables, and recovery of the sharp regularity in the original variables. Although the self-similar profiles and their stability estimates depend on the leading equation, the common structure of the arguments suggests that the profile-dependent parts of the analysis can be treated in a unified way.

We develop such a general framework for self-similar cusp formation in evolution equations whose leading dynamics are governed by the generalized Hunter--Saxton (gHS) family. As a main result, we use this framework to establish finite-time gradient blow-up with sharp H\"older regularity for a large class of nonlocal transport equations including the \(b\)-family and other physically relevant models. The resulting general theorem also applies to other models with leading dynamics in the gHS family: once its model-specific hypotheses are verified, the corresponding singularity formation results follow without repeating the profile analysis or the common stability argument.

\subsection{A class of nonlocal equations and their leading equations}\label{sec:intro-b-family}

We consider the initial value problem for the following class of nonlocal transport equations:
\begin{equation}\label{eq:b-family}
w_t-w_{txx}+\bigl(w+g'(w)\bigr)w_x =bw_xw_{xx}+ww_{xxx},\qquad x\in\R,\quad t>0,
\end{equation}
with
\begin{equation}\label{eq:b-initial-data-intro}
w(x,0)=w_0(x),\qquad x\in\R,
\end{equation}
where \(1<b<5\) and \(g\in C^\infty(\R)\) satisfies \(g(0)=0\). The class \eqref{eq:b-family} covers the following models:
\begin{enumerate}[(i)]
\item The choice \(g(w)=bw^2/2\) gives the \(b\)-family
\begin{equation*} 
m_t+wm_x+bw_xm=0,\qquad m=w-w_{xx}. 
\end{equation*}
The cases \(b=2\) and \(b=3\) correspond to the Camassa--Holm and Degasperis--Procesi equations, respectively \cite{EscherYin2008,LiuYin2006}.

\item The generalized hyperelastic-rod wave equation studied by Coclite, Holden, and Karlsen \cite{CocliteHoldenKarlsen2005} is given by
\begin{equation} \label{ghr}
v_t-v_{txx}+\frac12\partial_x h(v)=\rho(2v_xv_{xx}+vv_{xxx}). 
\end{equation}
Here we take \(\rho\in\R\setminus\{0\}\) and \(h\in C^\infty(\R)\) with \(h(0)=0\). The rescaling \(w=\rho v\) gives \eqref{eq:b-family} with \(b=2\) and \(g(w)=\frac{\rho}{2}h(w/\rho)-\frac12w^2\). This \(b=2\) case includes the normalized rotation--Camassa--Holm equation for shallow water waves with the Coriolis effect \cite{ChenGuiLiu2018}.

\item Equation \eqref{eq:b-family} reduces to the regularized Burgers equation when \(b=2\) and \(g\equiv0\) \cite{GJCP,KKS-rSV} and to the Fornberg--Whitham equation when \(b=3\) and \(g(w)=w\) \cite{FornbergWhitham1978}.
\end{enumerate}

Equation \eqref{eq:b-family} can be rewritten as
\begin{subequations} \label{eq:b-nonlocal}
\begin{align}
& w_t+ww_x=-P_x, \label{u-nonlocal} \\
& P=K*\left(g(w)+\frac{3-b}{2}w_x^2\right),\qquad K(x)=\frac{1}{2}e^{-|x|}. \label{P-nonlocal}
\end{align}
\end{subequations}
Differentiating \eqref{u-nonlocal} with respect to $x$ and using the equation for \(P\), we obtain
\begin{equation} \label{eq:b-slope}
w_{xt}+ww_{xx}+\frac{b-1}{2}w_x^2=F,\qquad F:=g(w)-P.
\end{equation}
Thus the nonlocal pressure $P$ changes the coefficient of \(w_x^2\) from \(1\) to \((b-1)/2\). In the blow-up regime considered here, \(F\) is of lower order than the terms on the left-hand side, although it need not remain bounded. The resulting leading equation therefore takes the canonical form
\begin{equation}\label{eq:ghs}
u_{xt}+uu_{xx}+\frac{1}{\ell}u_x^2=0,\qquad \ell=\frac{2}{b-1}.
\end{equation}
This family is often referred to as the generalized Hunter--Saxton equation; when \(\ell=2\), \eqref{eq:ghs} is the Hunter--Saxton equation \cite{HunterSaxton}. Note that the condition \(1<b<5\) is equivalent to \(\ell>1/2\), for which the corresponding self-similar solutions form bounded cusps; see Remark~\ref{rem:parameter-range}.

\subsection{Main result for the nonlocal transport equations}

We turn to the main result for the nonlocal equations \eqref{eq:b-family}. For a bounded open set \(\Omega\subset\R\) and \(\alpha\in(0,1]\), we write
\begin{equation*}
[f]_{C^\alpha(\Omega)}:=\sup_{x,\widetilde x\in\Omega,\,x\neq\widetilde x}\frac{|f(x)-f(\widetilde x)|}{|x-\widetilde x|^\alpha}
\end{equation*}
for the H\"older seminorm of \(f\) on \(\Omega\). We also write \(A(t)\sim B(t)\) for positive functions \(A\) and \(B\) if there exists a constant \(C>0\), independent of \(t\), such that \(C^{-1}B(t)\leq A(t)\leq CB(t)\). 

The following theorem states our first main result.

\begin{theorem}\label{thm:b-cusp-intro}
Fix \(1<b<5\) and \(g\in C^\infty(\R)\) with \(g(0)=0\), and set
\begin{equation*} \ell:=\frac{2}{b-1},\qquad \gamma_b:=\frac{5-b}{b+3}. \end{equation*}
There exists a non-empty relatively open class of smooth initial data for which the corresponding solution \(w\) of \eqref{eq:b-family}--\eqref{eq:b-initial-data-intro} develops gradient blow-up at a point \((x_*,T_*)\in\R\times(0,\infty)\). More precisely, the solution remains smooth on \([0,T_*)\) and satisfies the following properties:
\begin{enumerate}[(i)]
\item It holds that
\begin{equation*} 
\sup_{0\leq t<T_*}\|w(\cdot,t)\|_{L^\infty(\R)}<\infty, \qquad \lim_{t\uparrow T_*}\|w_x(\cdot,t)\|_{L^\infty(\R)}=\infty. 
\end{equation*}
\item For any bounded open set \(\Omega\subset\R\) and any \(\alpha\in(0,\gamma_b]\),
\begin{equation*} 
\sup_{0\leq t<T_*}[w(\cdot,t)]_{C^\alpha(\Omega)}<\infty. 
\end{equation*}
\item For any bounded open set \(\Omega\subset\R\) and any \(\alpha\in(\gamma_b,1]\),
\begin{equation*} 
\begin{cases} 
\displaystyle\lim_{t\uparrow T_*}[w(\cdot,t)]_{C^\alpha(\Omega)}=\infty,&\text{if }x_*\in\Omega,\\ 
\displaystyle\limsup_{t\uparrow T_*}[w(\cdot,t)]_{C^\alpha(\Omega)}<\infty,&\text{if }x_*\notin\overline{\Omega}. 
\end{cases} 
\end{equation*}
If \(x_*\in\Omega\), the temporal blow-up rate is given by
\begin{equation*} 
[w(\cdot,t)]_{C^\alpha(\Omega)}\sim(T_*-t)^{-((2\ell+1)\alpha-(2\ell-1))/2} 
\end{equation*}
for all \(t\) sufficiently close to \(T_*\).
\end{enumerate}
\end{theorem}

The admissible class of initial data is described in Section~\ref{subsec:b-theorem}, where the result is stated with explicit assumptions in Theorem~\ref{thm:b-cusp}. The exponent \(\gamma_b\) is sharp in the sense that the local \(C^{\gamma_b}\) H\"older seminorm remains uniformly bounded up to \(T_*\), whereas, for any \(\alpha>\gamma_b\), the \(C^\alpha\) H\"older seminorm becomes unbounded in every neighborhood of \(x_*\) as \(t\uparrow T_*\). Through the rescaling in Section~\ref{sec:intro-b-family}, Theorem~\ref{thm:b-cusp-intro} also yields finite-time cusp formation with sharp \(C^{3/5}\) regularity for the generalized hyperelastic-rod wave equation \eqref{ghr}, for any fixed \(\rho\neq0\) and any smooth \(h\) satisfying \(h(0)=0\).

\subsection{A general framework for self-similar singularity formation}\label{sec:intro-profiles-stability}

Theorem~\ref{thm:b-cusp-intro} is obtained through a common analysis of the self-similar profiles associated with \eqref{eq:ghs} and their perturbations. We formulate this analysis as our second main result: a general theorem for singularity formation. For a given evolution equation satisfying model-specific hypotheses, the theorem gives a quantitative description of cusp formation in space and time, including localization of the singularity, sharp H\"older regularity, and temporal blow-up rates. To state it, we first describe the profiles and the corresponding perturbed evolution equations.

\subsubsection{Self-similar blow-up profiles}\label{sec:intro-profiles}

We consider smooth odd profiles whose slope has a nondegenerate negative minimum at the origin. The corresponding self-similar solutions of \eqref{eq:ghs} are given by
\begin{equation} \label{sim_exp}
u(x,t)=(-t)^{(2\ell-1)/2}U(y),\qquad y=\frac{x}{(-t)^{(2\ell+1)/2}},\qquad t<0. 
\end{equation}
Substitution into \eqref{eq:ghs} yields the profile equation
\begin{equation}\label{eq:profile} 
\left(1+\frac{1}{\ell}U'\right)U'+\left(U+\frac{2\ell+1}{2}y\right)U''=0. 
\end{equation}
For every real \(\ell>1/2\), we construct a one-parameter family \(\{U_\beta\}_{\beta>0}\) of smooth odd strictly decreasing solutions of \eqref{eq:profile} on \(\R\); see Proposition~\ref{prop:profile-family}. These profiles satisfy
\begin{equation*} 
U_\beta(0)=0,\qquad U_\beta'(0)=-\ell,\qquad U_\beta''(0)=0,\qquad U_\beta^{(3)}(0)=\frac{\ell}{2}(2\ell)^{2\ell}\beta, 
\end{equation*}
and have the following far-field asymptotics:
\begin{equation}\label{eq:intro-profile-tail} 
\lim_{|y|\to\infty}|y|^{2/(2\ell+1)}|U_\beta'(y)|=\Theta_\ell\beta^{-1/(2\ell+1)} 
\end{equation}
and
\begin{equation*} 
U_\beta(y)=-\frac{2\ell+1}{2\ell-1}\Theta_\ell\beta^{-1/(2\ell+1)}\sgn(y)|y|^{\gamma_\ell}+o\bigl(|y|^{\gamma_\ell}\bigr)\quad\text{as }|y|\to\infty, 
\end{equation*}
where \(\Theta_\ell=\bigl(\ell/(2^{2\ell-2}(2\ell+1)^2)\bigr)^{1/(2\ell+1)}\) and
\begin{equation}\label{eq:cusp-exponent} 
\gamma_\ell=\frac{2\ell-1}{2\ell+1}. 
\end{equation}
The far-field growth of these profiles determines the sharp H\"older exponent \(\gamma_b\) in Theorem~\ref{thm:b-cusp-intro} through the relation \(\ell=2/(b-1)\).

\subsubsection{Perturbative stability and singularity formation}\label{sec:intro-stability}

For \(\ell>1/2\), we establish stability estimates around the profile \(U=U_\beta\) described above, with \(\beta\) chosen as in \eqref{eq:normalized-profile-choice} of Section~\ref{sec:profile-inequalities}.

We establish stability estimates around the solution  \(U=U_\beta\) of the ODE problem \eqref{eq:profile}--\eqref{eq:intro-profile-tail} for \(\ell>1/2\), with \(\beta\) chosen as in \eqref{eq:normalized-profile-choice} of Section~\ref{sec:profile-inequalities}. Consider an evolution equation on \(\R\), or a system of such equations, and let \(w\) denote the scalar unknown under consideration. Let \(\tau,\kappa,\xi\) be modulation functions with \(\tau(t)>t\), and introduce
\begin{equation}\label{eq:intro-modulated-variables}
y(x,t):=\frac{x-\xi(t)}{(\tau(t)-t)^{(2\ell+1)/2}},\qquad s(t):=-\log(\tau(t)-t). 
\end{equation}
We define \(W\) by
\begin{equation}\label{eq:intro-modulated-unknown} 
w(x,t)-\kappa(t)=e^{-(2\ell-1)s/2}W(y,s). 
\end{equation}
The modulation functions are chosen so that
\begin{equation}\label{eq:modulation-constraints-intro} 
W(0,s)=0,\qquad W_y(0,s)=-\ell,\qquad W_{yy}(0,s)=0. 
\end{equation}
Suppose that the equation for \(W_y\) can be written as
\begin{equation}\label{eq:intro-perturbed-slope} 
\partial_sW_y+\left(1+\frac1\ell W_y\right)W_y+\left(\frac{2\ell+1}{2}y+W+\mathcal A\right)W_{yy}=\mathcal R, 
\end{equation}
where \(\mathcal A\) and \(\mathcal R\) contain the model-specific perturbative terms, including those arising from the modulation. They may depend on the solution and, for a system, on its other components. For \eqref{eq:b-nonlocal}, their precise expressions are given in \eqref{eq:b-AR}. Note that setting \(\mathcal A=\mathcal R=0\) and \(W(y,s)=U(y)\) reduces \eqref{eq:intro-perturbed-slope} to the profile equation \eqref{eq:profile}.

The main stability estimates are of the form
\begin{equation} \label{main_stab_est}
|W_y(y,s)-U'(y)|\leq C\min\left\{\frac{y^2}{1+y^2},\frac{1}{1+|y|^{2/(2\ell+1)}}\right\} 
\end{equation}
for all \(y\in\R\) and all self-similar times \(s\geq s_0:=s(t_0)\), with \(C>0\) independent of \(s\). To obtain these estimates, we construct a parameter-dependent bootstrap scheme adapted to \(U\). The key ingredient is a family of profile inequalities that we establish for all real \(\ell>1/2\) and use to close the bootstrap argument. Together with the modulation estimates and the self-similar scaling, these weighted bounds yield sharp H\"older regularity and temporal blow-up rates in the original variables.

We state an informal version of the resulting theorem. The precise statement is given in Theorem~\ref{thm:common-singularity} in Section~\ref{sec:common-singularity}.

\begin{theorem}[Informal statement]\label{thm:common-singularity-intro}
Let \(\ell>1/2\) be real. Consider an evolution equation on \(\R\), or a system of such equations, with smooth initial data, and let \(w\) denote the scalar unknown under consideration. Suppose that the following hold.

\begin{enumerate}[(i)]
\item \emph{(Local existence and initial bounds)}
The equation admits a unique local smooth solution and modulation functions for which \eqref{eq:intro-modulated-variables}--\eqref{eq:intro-perturbed-slope} hold. The initial data $W(\cdot,s_0)$ satisfy the bootstrap bounds with sufficient margin.

\item \emph{(Bounds on the model-specific perturbative terms)}
On every bootstrap interval, the model-specific perturbative terms satisfy the required bounds with uniform constants.

\item \emph{(Modulation estimates)}
On every bootstrap interval, the modulation functions satisfy
\begin{equation*} c\leq1-\dot\tau(t)\leq C,\qquad |\kappa(t)|\leq C,\qquad \int_{t_0}^{t}|\dot\xi(t')|\,dt'\leq C \end{equation*}
with uniform constants \(c>0\) and \(C>0\), where \(t_0\) is the initial time.

\item \emph{(Continuation)}
The solution and the modulation functions can be continued beyond any finite \(S\), provided the bootstrap bounds hold on \([s_0,S)\).
\end{enumerate}
Then \(W\) is defined for all self-similar times \(s\geq s_0:=s(t_0)\), and the bounds in \eqref{eq:bootstrap-improved-bounds} hold on \(\R\times[s_0,\infty)\). In the original variables, the solution \(w\) develops finite-time gradient blow-up and exhibits sharp \(C^{\gamma_\ell}\) H\"older regularity up to the blow-up time, where \(\gamma_\ell\) is given by \eqref{eq:cusp-exponent}.
\end{theorem}

The required bounds on the model-specific perturbative terms are specified in Section~\ref{sec:remainder-bounds}. A convenient sufficient condition for these bounds is the following: on any finite bootstrap interval $[s_0,s_1]$, there exists a nonnegative continuous function $e$ such that
\begin{equation}\label{eq:sufficient-remainder-bounds}
\left\|\frac{\mathcal A(y,s)}{1+|y|}\right\|_{L_y^\infty}+\left\|(1+|y|^{2/(2\ell+1)})\mathcal R(y,s)\right\|_{L_y^\infty}+\sum_{j=2}^4\|\mathcal R_j(y,s)\|_{L_y^\infty}\leq e(s)
\end{equation}
for all $s\in[s_0,s_1]$, and
\begin{equation}\label{eq:sufficient-remainder-smallness}
\int_{s_0}^{s_1}e(s)\,ds\leq C_*,
\end{equation}
where $C_*>0$ is sufficiently small, depending only on the reference profile $U$ and bootstrap parameters, and is independent of $s_0$ and $s_1$. Here $\mathcal R_j$, $j=2,3,4$, are the remainders in the differentiated equations, defined in \eqref{eq:differentiated-remainders}; see Remark~\ref{rem:sufficient-remainder-bounds}.

For the nonlocal equations~\eqref{eq:b-nonlocal}, we verify the sufficient condition~\eqref{eq:sufficient-remainder-bounds}--\eqref{eq:sufficient-remainder-smallness}, together with the remaining hypotheses of Theorem~\ref{thm:common-singularity}. Applying this general theorem then yields Theorem~\ref{thm:b-cusp-intro}.

For a system, local existence and continuation concern all components, while the regularity conclusions concern \(w\). The construction of initial data satisfying the required margins and the verification of the model-specific estimates are carried out for the given equation. Once these conditions are verified, the profile inequalities, the weighted stability estimates, and the derivation of sharp regularity need not be repeated.

The leading dynamics of the regularized Saint--Venant and Novikov equations are also governed by the gHS family, with \(\ell=2\) and \(\ell=4\), respectively \cite{KKS-rSV,KKKS-Novikov}. The present framework gives a common formulation of the self-similar stability analyses carried out in these works. Further possible applications include the generalized Camassa--Holm--Novikov family and the asymptotic variational wave equations of Bressan, Zhang, and Zheng \cite{BressanZhangZheng2007}. Their formal leading-order reductions to \eqref{eq:ghs} are given in Appendix~\ref{app:formal-reductions}.

\subsection{Related works}\label{sec:related-works}

\noindent\textit{Self-similar approaches to singularity formation.}
Buckmaster, Shkoller, and Vicol developed a self-similar approach to shock formation for the compressible Euler equations \cite{BSV2D,BSV3D}, using self-similar Burgers profiles and obtaining \(C^{1/3}\) cusp regularity. Related analyses based on the stable Burgers profiles yield the same \(C^{1/3}\) cusp regularity for several perturbations of the Burgers equation and for Euler--Poisson systems \cite{BKK-delta, BKK-ion, ChickeringMorenoVasquezPandya2023, OhPasqualotto, SautWang2025, WangWang2026,  YangBH2021}.

The leading equation in the self-similar analysis of nonlocal wave breaking need not be the Burgers equation. For the Camassa--Holm equation and Hamiltonian regularizations of the Saint--Venant and Burgers equations, the leading-order dynamics are governed by the Hunter--Saxton equation, and the resulting cusps have sharp \(C^{3/5}\) regularity \cite{KKS-rSV,KKY-CH}. For the Novikov equation, the leading dynamics are governed by another member of the generalized Hunter--Saxton family, whose profiles yield sharp \(C^{7/9}\) regularity \cite{KKKS-Novikov}.

Elgindi, Ghoul, and Masmoudi \cite{ElgindiGhoulMasmoudi2021} established stable self-similar blow-up for a family of one-dimensional nonlocal vorticity models generalizing the Constantin--Lax--Majda equation. Their analysis combines modulation and weighted energy estimates to control perturbations of previously constructed profiles. They also obtained uniform H\"older bounds for the associated velocity up to blow-up, with exponents tending to zero as the negative advection parameter tends to zero for smooth profiles, and tending to one as the regularity parameter tends to zero for nonsmooth profiles with fixed advection.

\medskip
\noindent\textit{Wave breaking for nonlocal equations.}
The class \eqref{eq:b-nonlocal} covers models for shallow water waves and nonlinear waves in hyperelastic rods, for which wave-breaking criteria have been established. Constantin and Escher established wave-breaking criteria for the Camassa--Holm equation and related shallow-water models, and studied the blow-up rate and blow-up set for the periodic Camassa--Holm equation \cite{ConstantinEscher1998, ConstantinEscher2000}. For the \(b\)-family, Escher and Yin and Gui, Liu, and Tian obtained blow-up scenarios and sufficient conditions on the initial data for finite-time blow-up \cite{EscherYin2008,GuiLiuTian2008}. For the Degasperis--Procesi equation, Liu and Yin \cite{LiuYin2006} showed that the first loss of regularity occurs through wave breaking and gave sufficient conditions for its occurrence. Coclite, Gargano, and Sciacca \cite{CocliteGarganoSciacca2012} also studied singularity formation in the \(b\)-family numerically; their computations suggest a cubic-root-type singularity for the Degasperis--Procesi equation.

For generalized hyperelastic-rod equations, Yin \cite{Yin2004GeneralizedCH} obtained blow-up scenarios and finite-time blow-up results for smooth solutions. Further sufficient conditions for wave breaking for the classical hyperelastic-rod equation, corresponding to \(h(v)=3v^2\), were established by Guo and Zhou \cite{GuoZhou2009}. Brandolese and Cortez \cite{BrandoleseCortez2014} developed local-in-space blow-up criteria for a related class of generalized rod equations, under structural assumptions on the nonlinearities. A complementary global weak-solution theory for generalized hyperelastic-rod equations was developed by Coclite, Holden, and Karlsen \cite{CocliteHoldenKarlsen2005} for initial data in \(H^1(\R)\).

\medskip
\noindent\textit{Generic regularity of conservative solutions.}
A different approach studies the regularity of global conservative solutions via semilinear formulations in characteristic coordinates. Suitable transformed variables remain regular through wave breaking, while gradient blow-up corresponds to degeneracy of the map back to physical coordinates. Bressan and Chen \cite{BressanChen2017} established generic piecewise smoothness of conservative solutions to a nonlinear variational wave equation. Bressan, Huang, and Yu \cite{BressanHuangYu2015} described the corresponding structurally stable singularities using local space-time expansions. Generic regularity has also been studied for equations of Hunter--Saxton type and equations of Camassa--Holm type \cite{CaiChenShenTan2017,LiZhang2017}.

Using this approach, Armstrong, Chen, Huang, and Shen \cite{ArmstrongChenHuangShen2026} studied two families of equations depending on a parameter, one of Hunter--Saxton type and the other of Camassa--Holm type. For generic smooth initial data, they obtained local expansions in space and time near two types of structurally stable singularities: points on singular characteristic curves and isolated points at which such curves are created or terminate. The H\"older exponent at points of the second type, called Type~II in that work, is \((2-\lambda)/(2+\lambda)\). Their stated results concern \(\lambda=1/n\), with integers \(n\geq2\), in the Hunter--Saxton case, and a prescribed family of nonlocal equations with \(\lambda=1/(2a)\), with integers \(a\geq1\), in the Camassa--Holm case.

\medskip
\noindent\textit{The present work.}
We use the self-similar approach to study cusp formation for the class of nonlocal transport equations \eqref{eq:b-nonlocal}, which includes a variety of physically relevant models, as described in Section~\ref{sec:intro-b-family}. For a class of smooth initial data, we prove finite-time gradient blow-up with sharp \(C^{(5-b)/(b+3)}\) H\"older regularity, together with quantitative estimates describing the formation of the singularity in space and time. As \(b\) varies over \((1,5)\), the sharp cusp exponent \((5-b)/(b+3)\) ranges over the entire interval \((0,1)\).

For the Camassa--Holm, regularized Burgers, and Fornberg--Whitham equations, Theorem~\ref{thm:b-cusp-intro} recovers the sharp cusp regularity established in \cite{KKY-CH,KKS-rSV,SautWang2025}. For the Degasperis--Procesi equation, it gives finite-time cusp formation with sharp \(C^{1/3}\) regularity, complementing the earlier wave-breaking criteria and numerical observations \cite{LiuYin2006,CocliteGarganoSciacca2012}. For the generalized hyperelastic-rod family, it establishes sharp \(C^{3/5}\) cusp regularity, adding to the earlier wave-breaking results.

As discussed above, Armstrong, Chen, Huang, and Shen \cite{ArmstrongChenHuangShen2026} described the singular set and local space-time behavior of generic conservative solutions for families of Hunter--Saxton and Camassa--Holm type, including beyond the first breaking time. Under the correspondence \(\lambda=1/\ell\), their Type~II H\"older exponent agrees with the exponent \(\gamma_\ell=(2\ell-1)/(2\ell+1)\) associated with the gHS profiles; their results concern discrete parameter values corresponding to \(\ell\geq2\). In our self-similar framework, we treat all real \(\ell>1/2\) through profile inequalities and a common bootstrap argument for perturbations of the gHS profiles. Under the model-specific hypotheses, the framework provides explicit sufficient conditions on smooth initial data for the formation of a first cusp at a single point. It controls the rescaled evolution from these data up to cusp formation, with stability estimates yielding sharp regularity and temporal blow-up rates.

\subsection{Main ideas of the proof}\label{sec:intro-strategy}

The main task is to establish the stability estimates \eqref{main_stab_est} for all self-similar times $s\geq s_0$ by a bootstrap argument. These estimates control the perturbation $W_y-U'$ near $y=0$ and give a bound for $W_y$ with the same far-field decay rate as $U'$. Together with the local derivative and modulation bounds, they yield the sharp H\"older regularity and temporal blow-up rates in the original variables.

\begin{figure}[htbp]
\centering
\begin{tikzpicture}[
  x=1cm,y=1cm,
  box/.style={draw,rounded corners=1pt,align=center,inner sep=5pt,
    font=\small,text width=6.1cm,minimum height=1.65cm},
  wide/.style={box,text width=10.8cm,minimum height=1.15cm},
  match/.style={{Stealth[length=2mm]}-{Stealth[length=2mm]},line width=.5pt},
  flow/.style={-{Stealth[length=2.2mm]},line width=.5pt}
]
\node[box,line width=.8pt] (profiles) at (-3.8,0)
  {\textbf{Profile inequalities}\\[2pt]
   Damping estimates for\\each reference profile};
\node[box] (comparison) at (3.8,0)
  {\textbf{Generalized maximum principle}\\[2pt]
   Allow comparison bounds\\to vary in time};
\node[box] (bootstrap) at (-3.8,-2.9)
  {\textbf{Bootstrap scheme}\\[2pt]
   Choose stability bounds\\adapted to the profile};
\node[box] (remainders) at (3.8,-2.9)
  {\textbf{Model-specific bounds}\\[2pt]
   Control errors within\\the available margins};
\node[wide] (conclusion) at (0,-6.1)
  {\textbf{Common singularity-formation theorem}\\[2pt]
   Stability estimates and sharp regularity};
\draw[match] (profiles.east) -- (comparison.west);
\draw[match] (profiles.south) -- (bootstrap.north);
\draw[match] (comparison.south) -- (remainders.north);
\draw[match] (bootstrap.east) -- (remainders.west);
\coordinate (join) at (0,-4.4);
\draw[line width=.5pt] (bootstrap.south) |- (join);
\draw[line width=.5pt] (remainders.south) |- (join);
\draw[flow] (join) -- node[right,font=\footnotesize,align=left]
  {Strict improvement\\and continuation} (conclusion.north);
\end{tikzpicture}
\caption{The ingredients of the common bootstrap argument. All four ingredients depend on $\ell$ and must be adjusted together to remain compatible as $\ell$ varies over $(1/2,\infty)$. The double arrows indicate this interdependence.}
\label{fig:common-proof-structure}
\end{figure}
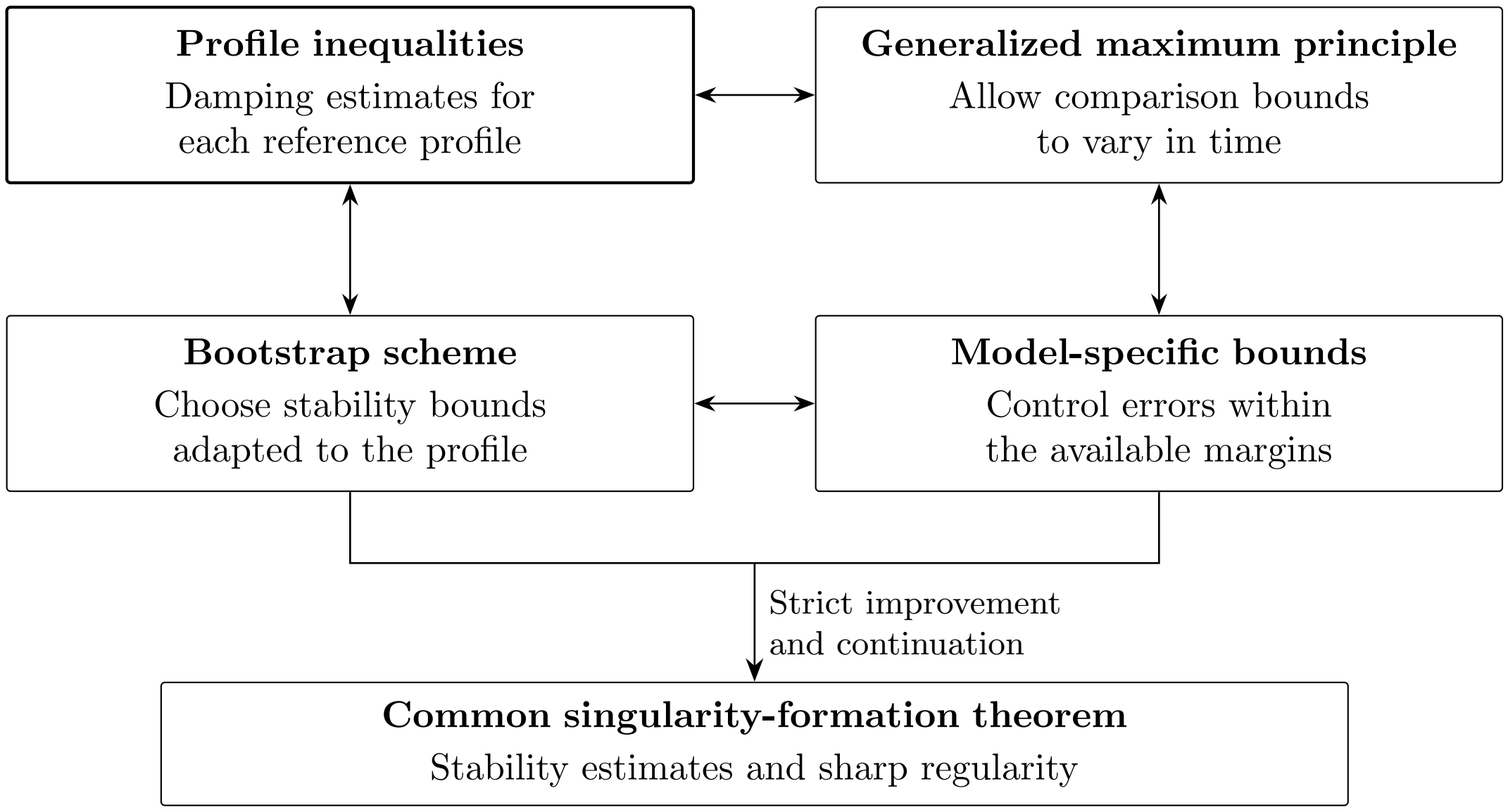

The main difficulty here is that the leading equation itself varies with $\ell$ through the gHS family \eqref{eq:ghs}. The self-similar scaling, the far-field behavior of the reference profile, and the coefficients in the perturbation equations all depend on $\ell$. Consequently, the bootstrap scheme and the admissible bounds on model-specific errors must be tailored to the estimates available for each profile. Our aim is to develop a common framework for all real $\ell>1/2$. The main ingredients and their interplay are illustrated in Figure~\ref{fig:common-proof-structure}.

We first construct the self-similar profiles of the gHS family and derive the evolution equation for the general perturbation $W_y-U'$. We observe that this equation contains the term
\begin{equation*}
-U''(y)(W-U)(y,s) =-U''(y)\int_0^y\bigl(W_y(z,s)-U'(z)\bigr)\,dz,
\end{equation*}
where we used $W(0,s)=U(0)=0$. In the weighted perturbation equations used to prove \eqref{main_stab_est}, we control this term using the corresponding supremum norm of the perturbation. To close these estimates by a maximum principle argument, we need profile inequalities showing that the damping dominates this integral contribution, with sufficient margin for the nonlinear and model-specific terms.

To obtain these inequalities, we use the following parametric representation for $y\geq0$:
\begin{equation*}
y=\frac{2}{(2\ell)^\ell\sqrt\beta}\int_0^\rho(1+z^2)^\ell\,dz,\qquad U_\beta'(y)=-\frac{\ell}{1+\rho^2},\qquad \rho\geq0.
\end{equation*}
This representation makes explicit how $\ell$ and $\beta$ determine the relation between the profile and the spatial variable $y$. We exploit this relation to control the profile in different spatial regions and choose $\beta$ depending on $\ell$ so that the inequalities required for the bootstrap argument hold for every $\ell>1/2$. This dependence is essential near $\ell=1/2$, since no fixed choice of $\beta$ satisfies the required conditions for all $\ell>1/2$; see Remark~\ref{rmk_beta}. The resulting inequalities are established in Propositions~\ref{prop:damping-inequalities}--\ref{prop:far-field}.

To use these inequalities in the presence of model-specific errors, we generalize the maximum principle to bound the supremum norm of the weighted perturbation by a time-dependent comparison function. The strict margins in the profile inequalities allow us to absorb some of the error terms, while the time-dependent comparison function controls the accumulated effect of the remaining errors. This leads to time-integral conditions on the model-specific terms. For each $\ell$, we choose the bootstrap bounds and the allowed errors so that this function remains strictly below the assumed bootstrap bound. Thus the profile inequalities determine the available margins, the model-specific assumptions control the errors within these margins, and the generalized maximum principle gives the bootstrap improvement.

With suitable initial margins, these estimates improve the assumed bootstrap bounds on every finite bootstrap interval. The local existence, modulation, and continuation hypotheses then extend the stability estimates to all self-similar times and yield Theorem~\ref{thm:common-singularity}. The constants may depend on $\ell$, but the structure of the argument is common to the whole family. For a particular equation, the remaining task is to construct admissible initial data and verify the model-specific hypotheses.

To facilitate this verification, we also record the sufficient condition~\eqref{eq:sufficient-remainder-bounds}--\eqref{eq:sufficient-remainder-smallness}, which bounds the relevant remainders by a common function $e(s)$ with a sufficiently small time integral. This condition is useful for the nonlocal equations~\eqref{eq:b-nonlocal} because estimates on the same nonlocal forcing also control the modulation and differentiated remainder terms. The slope bounds assumed in the bootstrap, together with an amplitude estimate, allow us to control the nonlocal pressure directly. Although the pressure can grow as the singular time is approached, its growth is integrable in physical time. After rescaling, we obtain a common remainder bound $e(s)$, independent of the upper endpoint of the bootstrap interval, such that
\begin{equation*}
\int_{s_0}^{\infty}e(s)\,ds\longrightarrow0
\qquad\text{as }s_0\to\infty.
\end{equation*}
Thus the accumulated effect of the remainders is sufficiently small when the initial self-similar time is sufficiently large. These estimates verify the perturbative hypotheses. Together with the remaining existence and modulation estimates and the construction of admissible initial data, they allow us to apply the common theorem and prove Theorem~\ref{thm:b-cusp-intro}.

\medskip
\emph{Plan of the paper.} Section~\ref{sec:profiles} constructs and classifies the self-similar profiles and states the inequalities used in the stability analysis. Section~\ref{sec:perturbative} states and proves the common bootstrap estimates and the sharp regularity result, and combines them into Theorem~\ref{thm:common-singularity}. Section~\ref{sec:coercivity} provides the deferred proofs of the profile inequalities. Section~\ref{sec:b-family} specifies the initial conditions and applies Theorem~\ref{thm:common-singularity} to prove Theorem~\ref{thm:b-cusp}. Appendix~\ref{subsec:b-data} constructs a non-empty relatively open class of admissible initial data. Appendix~\ref{app:b-local} establishes local existence and a blow-up
criterion for the nonlocal equations~\eqref{eq:b-nonlocal}. Appendix~\ref{app:transport-estimates} gives the far-field estimate and the maximum principle used in the bootstrap argument. Appendix~\ref{app:formal-reductions} gives formal reductions of the generalized Camassa--Holm--Novikov family and the asymptotic variational wave equations to \eqref{eq:ghs}.

\section{Construction of the self-similar profiles}\label{sec:profiles}

For each real $\ell>1/2$, we construct a one-parameter family of smooth odd solutions of the profile equation \eqref{eq:profile} whose slopes attain a nondegenerate minimum at the origin. We begin by deriving the profile equation and determining the similarity exponents in \eqref{sim_exp} from this nondegeneracy condition. After constructing and classifying the profiles, we record the inequalities needed for the stability analysis.

\subsection{Construction and classification}\label{sec:profile-construction}

For \eqref{eq:ghs}, we seek a self-similar solution centered at the origin, with blow-up time \(t=0\), of the form
\begin{equation*} 
u(x,t)=(-t)^aU(y),\qquad y=\frac{x}{(-t)^{a+1}},\qquad t<0. 
\end{equation*}
This choice of the amplitude and spatial scales gives \(u_x(x,t)=(-t)^{-1}U'(y)\). Substitution into \eqref{eq:ghs} yields
\begin{equation} \label{eq:general-profile-a} 
\left(1+\frac{1}{\ell}U'\right)U'+\left(U+(a+1)y\right)U''=0. 
\end{equation}
We consider smooth odd profiles whose negative slope has a nondegenerate minimum at the origin. The corresponding conditions are
\begin{equation*} 
U(0)=0,\qquad U'(0)<0,\qquad U''(0)=0,\qquad U^{(3)}(0)>0. 
\end{equation*}
Here nondegeneracy means that the slope \(U'\) increases quadratically away from its minimum, so that \(U^{(3)}(0)\) is positive. This corresponds to the conditions \(u_0''(0)=0\) and \(u_0^{(3)}(0)>0\) on the initial slope, where \(u_0:=u(\cdot,t_0)\) for a fixed initial time \(t_0<0\).

Evaluating \eqref{eq:general-profile-a} at the origin gives \(U'(0)=-\ell\). Writing \(U(y)=-\ell y+cy^3+O(y^5)\), where \(c>0\), and substituting into \eqref{eq:general-profile-a}, we obtain
\begin{equation*} 
a=\ell-\frac{1}{2},\qquad a+1=\ell+\frac{1}{2}. 
\end{equation*}
With this choice of \(a\), \eqref{eq:general-profile-a} becomes \eqref{eq:profile}.

The following proposition establishes the existence and key properties of these profiles.

\begin{proposition}\label{prop:profile-family}
Let \(\ell>1/2\). The profile equation \eqref{eq:profile} admits a one-parameter family \(\{U_\beta\}_{\beta>0}\) of smooth odd monotonically decreasing solutions on \(\mathbb R\) with the following properties for each \(\beta>0\):

\begin{enumerate}[(i)]
\item \emph{(Behavior near the origin)}
The profile \(U_\beta\) satisfies
\begin{equation}\label{eq:origin-derivatives}
U_\beta(0)=0,\qquad U_\beta'(0)=-\ell,\qquad U_\beta''(0)=0,\qquad U_\beta^{(3)}(0)=\frac{\ell}{2}(2\ell)^{2\ell}\beta,
\end{equation}
and
\begin{equation}\label{eq:origin-expansion}
U_\beta(y)=-\ell y+\frac{\ell(2\ell)^{2\ell}}{12}\beta y^3+O(y^5),\qquad y\to0.
\end{equation}

\item \emph{(Scaling relation)}
The profiles satisfy
\begin{equation}\label{eq:beta-scaling}
U_\beta(y)=\beta^{-1/2}U_1(\beta^{1/2}y),\qquad y\in\mathbb R.
\end{equation}

\item \emph{(Far-field asymptotics)}
Let
\begin{equation}\label{eq:theta-ell}
\Theta_\ell:=\left(\frac{\ell}{2^{2\ell-2}(2\ell+1)^2}\right)^{1/(2\ell+1)}.
\end{equation}
Then, as \(|y|\to\infty\),
\begin{align}
|y|^{-\gamma_\ell}|U_\beta(y)| &\longrightarrow \frac{2\ell+1}{2\ell-1}\Theta_\ell\beta^{-1/(2\ell+1)}, \label{eq:far-U} \\
|y|^{2/(2\ell+1)}|U_\beta'(y)| &\longrightarrow \Theta_\ell\beta^{-1/(2\ell+1)}, \label{eq:far-Uprime} \\
|y|^{(2\ell+3)/(2\ell+1)}|U_\beta''(y)| &\longrightarrow \frac{2}{2\ell+1}\Theta_\ell\beta^{-1/(2\ell+1)}. \label{eq:far-Usecond}
\end{align}

\item \emph{(Parametric representation)}
There exists a unique smooth strictly increasing bijection \(\rho_\beta:[0,\infty)\to[0,\infty)\) satisfying
\begin{equation}\label{eq:y-rho}
y=\frac{2}{(2\ell)^\ell\sqrt{\beta}}\int_0^{\rho_\beta(y)}(1+s^2)^\ell\,ds,\qquad y\geq0.
\end{equation}
For \(y\geq0\), writing \(\rho=\rho_\beta(y)\), we have
\begin{align}
U_\beta(y)&=\frac{1}{(2\ell)^\ell\sqrt{\beta}}\left(\rho(1+\rho^2)^\ell-(2\ell+1)\int_0^\rho(1+s^2)^\ell\,ds\right), \label{eq:U-rho} \\
U_\beta'(y)&=-\frac{\ell}{1+\rho^2}, \qquad U_\beta''(y)=\ell(2\ell)^\ell\sqrt{\beta}\,\frac{\rho}{(1+\rho^2)^{\ell+2}}, \label{eq:Uprime-rho}
\end{align}
and
\begin{equation}\label{eq:Uthird-rho}
U_\beta^{(3)}(y)=\frac{\ell}{2}(2\ell)^{2\ell}\beta\,\frac{1-(2\ell+3)\rho^2}{(1+\rho^2)^{2\ell+3}}.
\end{equation}
For \(y<0\), the identities \eqref{eq:U-rho}--\eqref{eq:Uthird-rho} hold with \(\rho=-\rho_\beta(-y)\).
\end{enumerate}
\end{proposition}

\begin{proof}
We first construct the desired one-parameter family of solutions. Since the equation \eqref{eq:profile} is invariant under the odd reflection, it suffices to work on the half-line \(\{y\in\mathbb R:y\geq0\}\). We set \(\ovU(y):=U(y)+(2\ell+1)y/2\). Then \eqref{eq:profile} can be reformulated as the autonomous ODE:
\begin{equation}\label{eq:W-equation}
\left(\ovU'-\frac12\right)\left(\ovU'-\frac{2\ell+1}{2}\right)+\ell \ovU\ovU''=0.
\end{equation}
We look for solutions \(\ovU\) to \eqref{eq:W-equation} satisfying \(\ovU'=V(\ovU)\) for some function \(V\). Using the relation \(\ovU''=V\frac{dV}{d\ovU}\), we obtain
\begin{equation*}
\left(V-\frac12\right)\left(V-\frac{2\ell+1}{2}\right)+\ell \ovU V\frac{dV}{d\ovU}=0.
\end{equation*}
By separation of variables, this gives
\begin{equation}\label{eq:first-integral}
\beta \ovU^2=\frac{2V-1}{(2\ell+1-2V)^{2\ell+1}},
\end{equation}
where \(\beta>0\) is a constant of integration. The right-hand side of \eqref{eq:first-integral} is smooth for \(1/2\leq V<(2\ell+1)/2\), and
\begin{equation*}
\frac{d}{dV}\left(\frac{2V-1}{(2\ell+1-2V)^{2\ell+1}}\right)=\frac{8\ell V}{(2\ell+1-2V)^{2\ell+2}}>0.
\end{equation*}
Moreover, it vanishes at \(V=1/2\) and tends to infinity as \(V\uparrow(2\ell+1)/2\). Hence, for each \(\beta>0\), the relation \eqref{eq:first-integral} defines a smooth even function \(V=V_\beta(\ovU)\) on \(\mathbb R\) such that
\begin{equation}\label{eq:V-bounds}
\frac12\leq V(\ovU)<\frac{2\ell+1}{2},\qquad \lim_{|\ovU|\to\infty}V(\ovU)=\frac{2\ell+1}{2},\qquad V(0)=\frac12.
\end{equation}
Differentiating \eqref{eq:first-integral} with respect to \(\ovU\), we obtain
\begin{equation}\label{eq:V-derivative}
\frac{dV}{d\ovU}=\frac{\beta \ovU(2\ell+1-2V)^{2\ell+2}}{4\ell V}.
\end{equation}
Using \eqref{eq:first-integral}, we also have
\begin{equation*}
\left|\frac{dV}{d\ovU}\right|=\frac{\sqrt{\beta}}{4\ell V}(2V-1)^{1/2}(2\ell+1-2V)^{\ell+3/2}.
\end{equation*}
Thus \(V\) and \(dV/d\ovU\) are uniformly bounded on \(\mathbb R\). Hence, by standard ODE theory, for each \(\beta>0\), the initial value problem
\begin{equation}\label{eq:W-ivp}
\ovU'=V(\ovU),\qquad \ovU(0)=0
\end{equation}
admits a unique smooth solution \(\ovU_\beta\) on \([0,\infty)\). Since \(V\) is bounded, this solution extends smoothly to the whole real line. The evenness of \(V\) and uniqueness for \eqref{eq:W-ivp} imply that \(\ovU_\beta\) is odd. Consequently, \(U_\beta(y):=\ovU_\beta(y)-(2\ell+1)y/2\) is a smooth odd solution of \eqref{eq:profile} on \(\mathbb R\). Moreover, since \(U_\beta'=V(\ovU_\beta)-(2\ell+1)/2\), \eqref{eq:V-bounds} gives
\begin{equation}\label{eq:profile-slope-bounds}
-\ell\leq U_\beta'(y)<0,\qquad y\in\mathbb R.
\end{equation}
Thus \(U_\beta\) is monotonically decreasing.

Next, we prove \eqref{eq:origin-expansion} for the function \(U_\beta\) constructed above. We first verify \eqref{eq:origin-derivatives}. The oddness of \(\ovU_\beta\), together with \eqref{eq:W-ivp} and \(V(0)=1/2\), gives
\begin{equation*}
U_\beta(0)=\ovU_\beta(0)=0,\qquad U_\beta'(0)=\ovU_\beta'(0)-\frac{2\ell+1}{2}=-\ell,\qquad U_\beta''(0)=\ovU_\beta''(0)=0.
\end{equation*}
Differentiating \eqref{eq:V-derivative} with respect to \(\ovU\) gives \(V''(0)=\beta(2\ell)^{2\ell+1}\). Hence, differentiating \eqref{eq:W-ivp} twice in \(y\) and using \(\ovU_\beta''(0)=0\), we obtain
\begin{equation*}
U_\beta^{(3)}(0)=\ovU_\beta^{(3)}(0)=\ovU_\beta'(0)^2V''(0)=\frac{\ell}{2}(2\ell)^{2\ell}\beta.
\end{equation*}
This proves \eqref{eq:origin-derivatives}. Then, using Taylor's theorem together with \eqref{eq:origin-derivatives} and the oddness of \(U_\beta\), we obtain \eqref{eq:origin-expansion}.

To prove \eqref{eq:beta-scaling}, we note from \eqref{eq:first-integral} that
\begin{equation*}
V_\beta(\ovU)=V_1(\sqrt{\beta}\ovU).
\end{equation*}
Therefore, \(\beta^{-1/2}\ovU_1(\sqrt{\beta}y)\) solves the initial value problem \(\ovU'=V_\beta(\ovU)\), \(\ovU(0)=0\). By uniqueness, \(\ovU_\beta(y)=\beta^{-1/2}\ovU_1(\sqrt{\beta}y)\), which gives \eqref{eq:beta-scaling} after subtracting \((2\ell+1)y/2\).

We now prove the parametric representation \eqref{eq:y-rho}--\eqref{eq:Uthird-rho}. By the oddness of \(U_\beta\), it suffices to consider \(y\geq0\). From \eqref{eq:V-bounds} and \eqref{eq:W-ivp}, we have \(\ovU_\beta(y)\geq y/2\). Define \(\rho\geq0\) by
\begin{equation*}
\rho^2=\frac{\ell+U_\beta'}{-U_\beta'}.
\end{equation*}
This is well defined by \eqref{eq:profile-slope-bounds}. Since \(U_\beta'=V-(2\ell+1)/2\), we have
\begin{equation}\label{eq:V-rho}
V=\frac{1+(2\ell+1)\rho^2}{2(1+\rho^2)},\qquad 2\ell+1-2V=\frac{2\ell}{1+\rho^2},\qquad 2V-1=\frac{2\ell\rho^2}{1+\rho^2}.
\end{equation}
The first identity in \eqref{eq:V-rho} gives the formula for \(U_\beta'\) in \eqref{eq:Uprime-rho}. Substituting the other two identities into \eqref{eq:first-integral} and using \(\ovU_\beta\geq0\), we obtain
\begin{equation}\label{eq:W-rho}
\ovU_\beta=\frac{\rho(1+\rho^2)^\ell}{(2\ell)^\ell\sqrt{\beta}}.
\end{equation}
The right-hand side of \eqref{eq:W-rho} has strictly positive derivative for \(\rho\geq0\). Since \(\ovU_\beta'(y)>0\), the relation \eqref{eq:W-rho} defines \(y\) as a smooth function of \(\rho\). Using \eqref{eq:V-rho} and \eqref{eq:W-rho}, we obtain
\begin{equation}\label{eq:dy-drho}
\frac{dy}{d\rho}=\frac{1}{V}\frac{d\ovU_\beta}{d\rho}=\frac{2(1+\rho^2)^\ell}{(2\ell)^\ell\sqrt{\beta}}.
\end{equation}
Since \(\rho=0\) when \(y=0\), integration gives \eqref{eq:y-rho}. Its right-hand side, viewed as a function of \(\rho\geq0\), is smooth, has strictly positive derivative, vanishes at \(\rho=0\), and tends to infinity as \(\rho\to\infty\). It therefore has a unique smooth strictly increasing inverse \(\rho_\beta:[0,\infty)\to[0,\infty)\). In particular, \(\rho=\rho_\beta(y)\). Combining \eqref{eq:y-rho} and \eqref{eq:W-rho} with \(U_\beta=\ovU_\beta-(2\ell+1)y/2\), we obtain \eqref{eq:U-rho}.

Differentiating the formula for \(U_\beta'\) with respect to \(y\) and using \eqref{eq:dy-drho}, we have
\begin{equation} \label{Ubeta''}
U_\beta''(y)=\left(\frac{dy}{d\rho}\right)^{-1}\frac{d}{d\rho}\left(-\frac{\ell}{1+\rho^2}\right)=\ell(2\ell)^\ell\sqrt{\beta}\,\frac{\rho}{(1+\rho^2)^{\ell+2}}.
\end{equation}
This proves the second identity in \eqref{eq:Uprime-rho}. Since
\begin{equation*}
\frac{d}{d\rho}\left(\frac{\rho}{(1+\rho^2)^{\ell+2}}\right)=\frac{1-(2\ell+3)\rho^2}{(1+\rho^2)^{\ell+3}},
\end{equation*}
we obtain \eqref{eq:Uthird-rho} by differentiating \eqref{Ubeta''} with respect to \(y\) and using \eqref{eq:dy-drho}. The identities for \(y<0\) follow from the oddness of \(U_\beta\), with \(\rho=-\rho_\beta(-y)\).

It remains to prove \eqref{eq:far-U}--\eqref{eq:far-Usecond}. By the oddness of \(U_\beta\), it suffices to consider \(y\to\infty\), or equivalently \(\rho\to\infty\). We have
\begin{equation*}
\lim_{\rho\to\infty}\rho^{-(2\ell+1)}\int_0^\rho(1+s^2)^\ell\,ds=\frac{1}{2\ell+1}.
\end{equation*}
Hence \eqref{eq:y-rho} gives
\begin{equation}\label{eq:y-rho-infinity}
\lim_{y\to\infty}\frac{y}{\rho_\beta(y)^{2\ell+1}}=\frac{2}{(2\ell+1)(2\ell)^\ell\sqrt{\beta}}.
\end{equation}
Using the first identity in \eqref{eq:Uprime-rho}, together with \eqref{eq:y-rho-infinity} and the definition \eqref{eq:theta-ell}, we obtain
\begin{equation*}
\lim_{y\to\infty}y^{2/(2\ell+1)}U_\beta'(y)=-\ell\left(\frac{2}{(2\ell+1)(2\ell)^\ell\sqrt{\beta}}\right)^{2/(2\ell+1)}=-\Theta_\ell\beta^{-1/(2\ell+1)}.
\end{equation*}
This proves \eqref{eq:far-Uprime} as \(y\to\infty\). Similarly, using the second identity in \eqref{eq:Uprime-rho}, we have
\begin{equation*}
\lim_{y\to\infty}y^{(2\ell+3)/(2\ell+1)}U_\beta''(y)=\frac{2}{2\ell+1}\Theta_\ell\beta^{-1/(2\ell+1)},
\end{equation*}
which proves \eqref{eq:far-Usecond} as \(y\to\infty\). Finally, the formula for \(U_\beta'\), together with \eqref{eq:dy-drho} and \(U_\beta(0)=0\), gives
\begin{equation*}
U_\beta(y)=-\frac{2\ell}{(2\ell)^\ell\sqrt{\beta}}\int_0^\rho(1+s^2)^{\ell-1}\,ds.
\end{equation*}
Since \(\ell>1/2\), we have
\begin{equation*}
\lim_{\rho\to\infty}\rho^{-(2\ell-1)}\int_0^\rho(1+s^2)^{\ell-1}\,ds=\frac{1}{2\ell-1}.
\end{equation*}
Combining these identities with \eqref{eq:y-rho-infinity} and \(\gamma_\ell=(2\ell-1)/(2\ell+1)\), we obtain
\begin{equation*}
\lim_{y\to\infty}y^{-\gamma_\ell}U_\beta(y)=-\frac{2\ell+1}{2\ell-1}\Theta_\ell\beta^{-1/(2\ell+1)}.
\end{equation*}
This proves \eqref{eq:far-U} as \(y\to\infty\). The case \(y\to-\infty\) follows from the oddness of \(U_\beta\). This completes the proof.
\end{proof}

\begin{remark}[The parameter range]\label{rem:parameter-range}
The parametrization in Proposition~\ref{prop:profile-family} is well defined for all real \(\ell>0\). The restriction \(\ell>1/2\) comes from the type of singularity considered in this paper. We study bounded cusp formation, in which the solution remains bounded near the singular point while its gradient blows up. For the corresponding self-similar solutions, this occurs precisely when \(\ell>1/2\).

For \(\ell>1/2\), the identity
\begin{equation*} (-t)^a\left|\frac{x}{(-t)^{a+1}}\right|^{\gamma_\ell}=|x|^{\gamma_\ell},\qquad \gamma_\ell=\frac{a}{a+1}, \end{equation*}
shows why growth of order \(|y|^{\gamma_\ell}\) gives a nontrivial limiting cusp for an exact self-similar solution. Here large \(|y|\) does not necessarily mean large \(|x|\). Indeed, a fixed small nonzero distance from the center corresponds to \(|y|\to\infty\) as \(t\uparrow0\).

For \(0<\ell\leq1/2\), the gradient still blows up, but the solution itself becomes unbounded in every fixed neighborhood of the origin as \(t\uparrow0\). Indeed, at \(\ell=1/2\), the logarithmic growth of the profile at infinity leads to divergence of the solution at each fixed \(x\neq0\). When \(0<\ell<1/2\), the profile has finite nonzero limits at infinity, but the amplitude factor \((-t)^{\ell-1/2}\) diverges. Thus, although the profiles can still be constructed in these cases, the associated solutions do not remain bounded up to the gradient blow-up time.
\end{remark}

\begin{remark}[The Burgers profile] \label{rem:burgers-profile}
When \(\ell=1\), the formulas \eqref{eq:y-rho} and \eqref{eq:U-rho} reduce to
\begin{equation*}
y=\frac{1}{\sqrt{\beta}}\left(\rho+\frac{\rho^3}{3}\right),\qquad U_\beta(y)=-\frac{\rho}{\sqrt{\beta}}.
\end{equation*}
Equivalently, \(U_\beta\) is determined by
\begin{equation*} 
y=-U_\beta(y)-\frac{\beta}{3}U_\beta(y)^3.
\end{equation*}
This cubic relation defines the self-similar profile for the formation of pre-shocks in the Burgers equation. In particular, it satisfies
\begin{equation*}
U_\beta(y)=-\left(\frac{3}{\beta}\right)^{1/3}\sgn(y)|y|^{1/3}+o\left(|y|^{1/3}\right) \qquad \text{as } |y|\to\infty.
\end{equation*}
\end{remark}

The following corollary shows that every smooth odd solution of \eqref{eq:profile} whose slope has a nondegenerate negative minimum at the origin belongs to the family constructed in Proposition~\ref{prop:profile-family}. In particular, such a profile is uniquely determined by its third derivative at the origin.

\begin{corollary} \label{cor:classification}
Let \(U\) be a smooth odd solution of \eqref{eq:profile} such that \(U'(0)<0\) and \(U^{(3)}(0)>0\). Then there exists a unique \(\beta>0\) such that \(U=U_\beta\). More precisely, the unique $\beta$ is given by
\begin{equation} \label{eq:beta-from-jet}
\beta=\frac{2U^{(3)}(0)}{\ell(2\ell)^{2\ell}}.
\end{equation}
\end{corollary}

\begin{proof}
Evaluating \eqref{eq:profile} at \(y=0\), we obtain \(U'(0)=-\ell\). Set \(\ovU(y):=U(y)+(2\ell+1)y/2\) and \(a_3:=U^{(3)}(0)>0\). By Taylor's theorem and oddness,
\begin{equation*}
\ovU(y)=\frac y2+\frac{a_3}{6}y^3+O(y^5),\qquad \ovU'(y)=\frac12+\frac{a_3}{2}y^2+O(y^4).
\end{equation*}
Hence, for sufficiently small \(y>0\), we have \(\ovU(y)>0\) and \(1/2<\ovU'(y)<(2\ell+1)/2\). On this interval, \(\ovU\) is invertible. Applying the separation-of-variables argument from the proof of Proposition~\ref{prop:profile-family} with \(V(\ovU(y))=\ovU'(y)\), we obtain \eqref{eq:first-integral} for some \(\beta>0\). Taking the limit as \(y\downarrow0\), we find
\begin{equation*}
\beta=\lim_{y\downarrow0}\frac{2\ovU'(y)-1}{\ovU(y)^2(2\ell+1-2\ovU'(y))^{2\ell+1}}=\frac{4a_3}{(2\ell)^{2\ell+1}},
\end{equation*}
which is \eqref{eq:beta-from-jet}. For this value of \(\beta\), \eqref{eq:first-integral} gives \(\ovU'=V_\beta(\ovU)\) for sufficiently small \(y\geq0\), where \(V_\beta\) is the smooth function constructed in the proof of Proposition~\ref{prop:profile-family}. Since \(\ovU(0)=\ovU_\beta(0)=0\), uniqueness for \eqref{eq:W-ivp} gives \(\ovU=\ovU_\beta\) on this interval.

Since \(\ovU_\beta(y)>0\) for any \(y>0\), the equation \eqref{eq:W-equation} is nonsingular along \(\ovU_\beta\). The usual uniqueness and continuation argument extends the equality to all \(y>0\), and oddness gives \(U=U_\beta\) on \(\mathbb R\). Finally, \eqref{eq:beta-from-jet} gives uniqueness of the parameter.
\end{proof}

\begin{remark}\label{rem:first-members}
For \(\ell=2\) and \(\ell=4\), the construction in Proposition~\ref{prop:profile-family} reduces to those in \cite{KKY-CH} and \cite{KKKS-Novikov}, respectively, with the same normalization of \(\beta\). More precisely, when \(\ell=2\), we have
\begin{equation*}
U_\beta^{(3)}(0)=256\beta,\qquad \Theta_2=50^{-1/5},
\end{equation*}
while for \(\ell=4\),
\begin{equation*}
U_\beta^{(3)}(0)=2^{25}\beta,\qquad \Theta_4=1296^{-1/9}.
\end{equation*}
\end{remark}

\subsection{Useful inequalities}\label{sec:profile-inequalities}

In this subsection, we record the profile inequalities used in the stability estimates of Section~\ref{sec:perturbative}, with the proofs deferred to Section~\ref{sec:coercivity}.  We take $U=U_\beta$ from the family constructed in Proposition~\ref{prop:profile-family}, with $\beta$ fixed in terms of $\ell$ by
\begin{equation}\label{eq:normalized-profile-choice}
\beta:=\begin{cases}
\displaystyle \max\left\{576,\left[2\Theta_\ell\left(\frac4\ell+\frac6{2\ell-1}\right)\right]^{2\ell+1}\right\}, &1/2<\ell<2, \\
1, &\ell\geq2.
\end{cases}
\end{equation}
Here, $\Theta_\ell$ is defined in \eqref{eq:theta-ell}. For the motivation behind this choice of \(\beta\), see Remark~\ref{rmk_beta}.

We first record bounds for \(U\) and its derivatives.

\begin{lemma} \label{lem:derivative-bounds}
Let \(\ell>1/2\), and let \(U=U_\beta\) be the profile constructed in Proposition~\ref{prop:profile-family}, with \(\beta\) chosen as in \eqref{eq:normalized-profile-choice}. Then the following holds:
\begin{subequations}
\begin{align}
&\|U'\|_{L^\infty}=\ell, \label{eq:Uprime-norm} \\
&\|U''\|_{L^\infty}= \ell(2\ell)^{\ell}\sqrt{\beta}\,\frac{(2\ell+3)^{\ell+3/2}}{(2\ell+4)^{\ell+2}}, \label{eq:Usecond-norm} \\
& \|U^{(3)}\|_{L^\infty}=\frac{\ell}{2}(2\ell)^{2\ell}\beta. \label{eq:Uthird-norm}
\end{align}
\end{subequations}
Moreover, there exists \(C_{\ell}>0\) such that
\begin{equation}\label{eq:Ufourth-bound}
\|U^{(4)}\|_{L^\infty}\leq C_\ell.
\end{equation}
\end{lemma}

Let us define
\begin{subequations}
\begin{align}
&D_{1,\ell}(y):=\begin{cases}
	1+\frac{2}{\ell}U'(y)+\frac{2}{1+y^2}\left(\frac{2\ell+1}{2}+\frac{U(y)}y\right), & y\neq 0,
	\\
	0, & y=0,
\end{cases}\label{eq:D1} \\
&D_{2,\ell}(y):=\begin{cases}
	\frac{2\ell+3}{2}+\left(1+\frac{2}{\ell}\right)U'(y)+\frac{1}{1+y^2}\left(\frac{2\ell+1}{2}+\frac{U(y)}y\right), & y\neq 0,
	\\
	0, & y=0.
\end{cases}\label{eq:D2}
\end{align}
\end{subequations}
By \eqref{eq:origin-derivatives} in Proposition~\ref{prop:profile-family}, both functions are continuous at \(y=0\).

\begin{proposition}\label{prop:damping-inequalities}
The following estimates hold for all $y \in \mathbb{R}$:
\begin{enumerate}[(i)]
\item For \(\ell>2\),
\begin{subequations}
\begin{align}
&\ell+U'(y)\geq(\ell-1)\frac{y^2}{1+y^2},\label{eq:first-damping}\\
&D_{1,\ell}(y)\geq\frac{\ell-2}{\ell}\frac{y^2}{1+y^2},\label{eq:D1-lower}\\
&D_{2,\ell}(y)\geq\left(\ell+\frac12-\frac{2}{\ell}\right)\frac{y^2}{1+y^2}.\label{eq:D2-lower}
\end{align}
\end{subequations}
\item For \(\ell=2\),
\begin{equation} \label{eq:CH-damping}
2+U'(y)\geq\frac65\frac{y^2}{1+y^2}, \qquad D_{1,2}(y)\geq\frac15\frac{y^2}{1+y^2}, \qquad D_{2,2}(y)\geq\frac{19}{10}\frac{y^2}{1+y^2}. 
\end{equation}
\item For \(1/2<\ell<2\),
\begin{subequations}\label{eq:small-ell-damping}
\begin{align}
&\ell+U'(y)\geq\frac{3\ell}{4}\frac{y^2}{1+y^2},\label{eq:small-ell-first-damping}\\
&D_{1,\ell}(y)\geq\frac12\frac{y^2}{1+y^2},\label{eq:small-ell-D1-lower}\\
&D_{2,\ell}(y)\geq\left(\frac{3\ell}{4}+1\right)\frac{y^2}{1+y^2}. \label{eq:small-ell-D2-lower}
\end{align}
\end{subequations}
\end{enumerate}
\end{proposition}

For the weighted perturbation estimate adapted to the behavior near $y=0$, we also define
\begin{equation}\label{eq:local-kernel}
K_{\ell}^{\mathrm{loc}}(y):=|U''(y)|\frac{1+y^2}{y^2}\int_0^{|y|}\frac{z^2}{1+z^2}\,dz, \qquad y \neq 0,
\end{equation}
with $K_\ell^{\mathrm{loc}}(0):=0$. Then, we have the following inequality.
\begin{proposition}\label{prop:local-integral}
For any \(\ell>1/2\), there exist constants \(\lambda_{\ell}^{\mathrm{loc}}>1\) and \(\eta_{\ell}>0\) such that
\begin{equation}\label{eq:local-integral-inequality}
\lambda_{\ell}^{\mathrm{loc}}K_{\ell}^{\mathrm{loc}}(y) \leq D_{1,\ell}(y)-\eta_{\ell}\frac{y^2}{1+y^2},\qquad y\in\R.
\end{equation}
\end{proposition}

Finally, for the far-field comparison, we define
\begin{equation}\label{eq:Jr}
J_\ell(y):=\int_0^y\frac{1}{1+|z|^{2/(2\ell+1)}} \, dz.
\end{equation}
For \(B>0\), set
\begin{equation}\label{eq:far-damping}
D_{\ell,B}^{\mathrm{far}}(y):=\begin{cases}
1-\frac{B}{\ell(1+|y|^{2/(2\ell+1)})}+\frac{2}{\ell}U'(y)-\frac{\frac{2}{2\ell+1}|y|^{2/(2\ell+1)}}{1+|y|^{2/(2\ell+1)}}\left(\frac{2\ell+1}{2}+\frac{U(y)}y+\frac B yJ_\ell(y)\right), & y\neq 0 \\
-1-\frac{B}{\ell}, & y=0
\end{cases}
\end{equation}
and
\begin{equation}\label{eq:far-kernel}
K_{\ell}^{\mathrm{far}}(y):=(1+|y|^{2/(2\ell+1)})|U''(y)|\int_0^{|y|}\frac{dz}{1+|z|^{2/(2\ell+1)}}.
\end{equation}
We choose the constants \(B_{\ell}\) and \(\Lambda_{\ell}\) as follows:
\begin{equation}\label{eq:B-choice}
\begin{split}
B_{\ell} &:= \begin{cases}
2\Theta_\ell\beta^{-1/(2\ell+1)},&1/2<\ell<2,\\
\frac{6}{13},&\ell=2,\\
\Theta_\ell+\frac13\left(1-\Theta_\ell\left(\frac3{\ell}+\frac2{2\ell-1}\right)\right),&2<\ell<4,\\
1,&\ell\geq4,
\end{cases} \\
\Lambda_{\ell} &:= \begin{cases}
2,&1/2<\ell<2,\\
\frac{10001}{10000},&\ell=2,\\
1+\frac13\left(1-\Theta_\ell\left(\frac3{\ell}+\frac2{2\ell-1}\right)\right),&2<\ell<4,\\
\frac{63}{62},&\ell\geq4.
\end{cases}
\end{split}
\end{equation}

\begin{proposition}\label{prop:far-field}
For any real $\ell>1/2$, the constant $\Lambda_\ell$ in \eqref{eq:B-choice} satisfies $\Lambda_\ell>1$, and there exists a constant $m_0(\ell)>0$ such that
\begin{equation}\label{eq:far-field-inequality}
\Lambda_{\ell}K_{\ell}^{\mathrm{far}}(y)\leq D_{\ell,B_{\ell}}^{\mathrm{far}}(y),\qquad |y|\geq m_0(\ell).
\end{equation}
\end{proposition}

The proofs of Lemma~\ref{lem:derivative-bounds} and Propositions~\ref{prop:damping-inequalities}--\ref{prop:local-integral} are given in Section~\ref{sec:localized-inequality-proofs}, while the proof of Proposition~\ref{prop:far-field} is given in Section~\ref{sec:far-field-inequality-proof}, where a choice of \(m_0(\ell)\) is specified in \eqref{eq:m0-definition}.

\section{Perturbative stability and sharp H\texorpdfstring{\"o}{o}lder regularity at blow-up} \label{sec:perturbative}

In this section, we establish stability estimates around the self-similar profile \(U=U_\beta\) for every \(\ell>1/2\) by a bootstrap argument under suitable bounds on the model-specific perturbative terms. We then show that these estimates, together with appropriate modulation assumptions, yield gradient blow-up with sharp H\"older regularity and temporal blow-up rates. Combining these results with assumptions on the initial data, local existence, and continuation, we obtain a general conditional theorem for singularity formation; see Theorem~\ref{thm:common-singularity}.

\subsection{Self-similar variables and the perturbation equations}\label{sec:similarity-variables}

Fix $\ell>1/2$ and let $U:=U_\beta$, with $\beta$ chosen as in \eqref{eq:normalized-profile-choice}. Let \(w=w(x,t)\) be a smooth scalar function, and let \(\tau,\kappa,\xi\) be modulation functions with \(\tau(t)>t\). We introduce the self-similar variables
\begin{equation}\label{eq:modulated-variables}
y(x,t):=\frac{x-\xi(t)}{(\tau(t)-t)^{(2\ell+1)/2}},\qquad s(t):=-\log(\tau(t)-t),
\end{equation}
and define \(W\) by
\begin{equation}\label{eq:modulated-unknown}
w(x,t)-\kappa(t)=e^{-(2\ell-1)s/2}W(y,s).
\end{equation}
The modulation functions are chosen so that \(W\) satisfies \eqref{eq:modulation-constraints-intro}. Suppose that the equation for \(W_y\) has the form
\begin{equation}\label{eq:rescaled-slope-toolkit}
\partial_sW_y+\left(1+\frac1\ell W_y\right)W_y+\left(\frac{2\ell+1}{2}y+W\right)W_{yy}=\mathcal E,
\end{equation}
where $\mathcal E$ collects the model-specific perturbative terms, including those arising from modulation. We seek the bounds
\begin{equation}\label{eq:regularity-profile-bounds}
|W_y(y,s(t))|\leq\frac{C}{1+|y|^{2/(2\ell+1)}},\qquad \sup_{|y|\leq1}|W_{yy}(y,s(t))|\leq C,
\end{equation}
with a constant $C>0$ independent of time. The corresponding bounds for \(U\) follow from \eqref{eq:far-Uprime} and its smoothness. Together with the normalization and the modulation bounds in Proposition~\ref{prop:sharp-regularity}, these estimates give the sharp regularity in the original variables. We prove \eqref{eq:regularity-profile-bounds} by estimating \(W_y-U'\) near the origin and in the far field, and by controlling the derivatives up to order four.

The transport coefficient in the rescaled equation may differ from \((2\ell+1)y/2+W\). We denote this difference by \(\mathcal A\) and retain it on the left-hand side. Accordingly, we write
\begin{equation}\label{eq:remainder-transport-split}
\mathcal E=\mathcal R-\mathcal A W_{yy},\qquad \mathcal V(y,s):=\frac{2\ell+1}{2}y+W(y,s)+\mathcal A(y,s).
\end{equation}
Then \eqref{eq:rescaled-slope-toolkit} becomes
\begin{equation}\label{eq:perturbed-slope-equation}
\partial_sW_y+\mathcal V W_{yy}+\left(1+\frac1\ell W_y\right)W_y=\mathcal R.
\end{equation}

Let
\begin{equation}\label{eq:f-definition-toolkit}
\widetilde W(y,s):=W(y,s)-U(y),\qquad \widetilde W(y,s)=\int_0^y\widetilde W_y(z,s)\,dz,
\end{equation}
where the second identity follows from \(W(0,s)=U(0)=0\). Subtracting \eqref{eq:profile} from \eqref{eq:rescaled-slope-toolkit}, we obtain
\begin{equation*}
\partial_s\widetilde W_y+\left(\frac{2\ell+1}{2}y+U\right)\widetilde W_{yy}+\left(1+\frac2\ell U'\right)\widetilde W_y+U''\widetilde W=\mathcal E-\frac1\ell\widetilde W_y^2-\widetilde W\widetilde W_{yy}.
\end{equation*}
Using \eqref{eq:remainder-transport-split}, we rewrite this as
\begin{equation}\label{eq:f-full-transport}
\partial_s\widetilde W_y+\mathcal V\widetilde W_{yy}+\left(1+\frac{\widetilde W_y+2U'}\ell\right)\widetilde W_y=-U''\widetilde W+\mathcal R-\mathcal A U''.
\end{equation}

For the derivative estimates, we apply $\partial_y^j$ to \eqref{eq:perturbed-slope-equation} for $j=1,2,3$ and obtain
\begin{align}
& \partial_s W_{yy}+\mathcal V\partial_y^3W+\left(\frac{2\ell+3}{2}+\left(1+\frac2\ell\right)W_y\right)W_{yy}=\mathcal R_2, \label{eq:second-derivative-common} \\
& \partial_s\partial_y^3W+\mathcal V\partial_y^4W+\left(2\ell+2+\left(2+\frac2\ell\right)W_y\right)\partial_y^3W=-\left(1+\frac2\ell\right)W_{yy}^2+\mathcal R_3, \label{eq:third-derivative-common} \\
& \partial_s\partial_y^4W+\mathcal V\partial_y^5W+\left(\frac{6\ell+5}{2}+\left(3+\frac2\ell\right)W_y\right)\partial_y^4W=-\left(4+\frac6\ell\right)W_{yy}\partial_y^3W+\mathcal R_4, \label{eq:fourth-derivative-common}
\end{align}
where
\begin{subequations} \label{eq:differentiated-remainders}
\begin{align}
&\mathcal R_2(y,s):=\mathcal R_y(y,s)-\mathcal A_y(y,s)W_{yy}(y,s),\label{eq:remainder-second}\\ &\mathcal R_3(y,s):=\mathcal R_{yy}(y,s)-2\mathcal A_y(y,s)\partial_y^3W(y,s)-\mathcal A_{yy}(y,s)W_{yy}(y,s),\label{eq:remainder-third}\\ &\mathcal R_4(y,s):=\partial_y^3\mathcal R(y,s)-3\mathcal A_y(y,s)\partial_y^4W(y,s)\notag\\ &\qquad\qquad\quad-3\mathcal A_{yy}(y,s)\partial_y^3W(y,s)-\partial_y^3\mathcal A(y,s)W_{yy}(y,s).\label{eq:remainder-fourth}
\end{align}
\end{subequations}
Using \eqref{eq:modulation-constraints-intro}, we evaluate \eqref{eq:perturbed-slope-equation} and \eqref{eq:second-derivative-common} at \(y=0\) to obtain
\begin{equation}\label{eq:remainder-compatibility}
\mathcal R(0,s)=0,\qquad \mathcal R_y(0,s)=\mathcal A(0,s)\partial_y^3W(0,s).
\end{equation}
In particular, neither \(\mathcal A(0,s)\) nor \(\mathcal R_y(0,s)\) is required to vanish separately.

\subsection{Bootstrap estimates and sharp regularity}\label{sec:common-bootstrap}

We formulate the bootstrap assumptions and specify the required bounds on the model-specific perturbative terms. The resulting stability estimates are then used to derive sharp regularity at blow-up and to establish a conditional theorem for singularity formation.

\subsubsection{Bootstrap assumptions} \label{sec:bootstrap-assumptions}

Choose \(m\geq\max\{1,m_0(\ell)\}\), where \(m_0(\ell)\) is the threshold in Proposition~\ref{prop:far-field}. We use a parameter
\begin{equation*}
0<\vartheta<\frac12
\end{equation*}
to specify the allowed increase in the bootstrap quantities.

Let \(N=N_{\ell}(m,\vartheta)>0\) be sufficiently small and \(M>0\) sufficiently large, to be specified below. On a finite interval \([s_0,s_1]\), we assume that
\begin{subequations}\label{eq:common-bootstrap-assumptions}
\begin{align}
&|W_y(y,s)-U'(y)|\leq\frac{Ny^2}{1+y^2},\label{eq:boot-local}\\
&|W_y(y,s)-U'(y)|\leq\frac{B_{\ell}}{1+|y|^{2/(2\ell+1)}},\label{eq:boot-far}\\
&|W_{yy}(y,s)|\leq\frac{M^{1/8}|y|}{(1+y^2)^{1/2}},\label{eq:boot-second}\\
&|\partial_y^3W(0,s)-U^{(3)}(0)|\leq1,\label{eq:boot-origin}\\
&\|\partial_y^3W(\cdot,s)\|_{L^\infty}\leq M^{3/4},\label{eq:boot-third}\\
&\|\partial_y^4W(\cdot,s)\|_{L^\infty}\leq M,\label{eq:boot-fourth}
\end{align}
\end{subequations}
for all \((y,s)\in\R\times[s_0,s_1]\), where $B_{\ell}$ is given by \eqref{eq:B-choice}. The two bounds \eqref{eq:boot-local} and \eqref{eq:boot-far} on \(W_y-U'\) control the perturbation near the origin and its decay at infinity, respectively. The derivative bounds \eqref{eq:boot-second}--\eqref{eq:boot-fourth} are used to estimate the solution near the origin and to close the bootstrap, as in \cite{KKY-CH,KKS-rSV,KKKS-Novikov}.

We specify the smallness of \(N\). By Proposition~\ref{prop:damping-inequalities}, we have
\begin{equation} \label{eq:bootstrap-profile-damping}
\ell+U'(y)\geq c_{1,\ell}\frac{y^2}{1+y^2},\qquad D_{2,\ell}(y)\geq c_{2,\ell}\frac{y^2}{1+y^2},\qquad y\in\R,
\end{equation}
where
\begin{equation*}
c_{1,\ell}:=\begin{cases}
3\ell/4,&1/2<\ell<2,\\
6/5,&\ell=2,\\\ell-1,&\ell>2,
\end{cases}\qquad c_{2,\ell}:= \begin{cases}
3\ell/4+1,&1/2<\ell<2,\\19/10,&\ell=2,\\
\ell+1/2-2/\ell,&\ell>2.
\end{cases}
\end{equation*}
Let $\eta_{\ell}$ and $\lambda_{\ell}^{\mathrm{loc}}$ be the constants in Proposition~\ref{prop:local-integral}. We choose \(N\) so that
\begin{subequations}\label{eq:common-smallness}
\begin{align}
&N\leq\min\left\{1,\frac{c_{1,\ell}}2\right\},\qquad \left(\frac1\ell+\frac23\right)N\leq\frac{\eta_{\ell}}{4}, \qquad \left(\frac43+\frac2\ell\right)N\leq\frac{c_{2,\ell}}4 \label{eq:N-local}\\
&N(1+m^{2/(2\ell+1)})\leq(1-2\vartheta)B_{\ell}.\label{eq:N-second-far}
\end{align}
\end{subequations}
The first and third conditions in \eqref{eq:N-local} use the lower bounds for \(\ell+U'\) and \(D_{2,\ell}\) in \eqref{eq:bootstrap-profile-damping}, while the second uses the comparison between \(D_{1,\ell}\) and \(\lambda_\ell^{\mathrm{loc}}K_\ell^{\mathrm{loc}}\) in Proposition~\ref{prop:local-integral}. These conditions ensure that the damping provided by the profile inequalities is sufficient for the maximum principle estimates. The condition \eqref{eq:N-second-far} ensures that the local bootstrap bound is strictly below the far-field bound \eqref{eq:boot-far} on \([-m,m]\).

The constant \(M\) is chosen sufficiently large in the proof of Proposition~\ref{prop:common-bootstrap}, depending only on the parameters fixed there and on the profile, but not on \(s_1\).

\subsubsection{Bounds on the model-specific perturbative terms}\label{sec:remainder-bounds}

We state the assumptions on the model-specific perturbative terms and provide a sufficient condition in Remark~\ref{rem:sufficient-remainder-bounds}, which will be used for the nonlocal equations \eqref{eq:b-family} in Section~\ref{sec:b-family}. All pointwise bounds below are assumed to hold for \(s\in[s_0,s_1]\) on the indicated spatial regions.

At the origin, we require
\begin{equation}\label{eq:bootstrap-origin-forcing}
\int_{s_0}^{s_1}\left|\mathcal R_3(0,s)-\mathcal A(0,s)\partial_y^4W(0,s)\right|\,ds\leq\frac{\vartheta N}{8}.
\end{equation}

For the estimates of \(\widetilde W_y\), we define
\begin{equation}\label{eq:bootstrap-core-radii}
h:=\frac{\vartheta N}{8(M+\|U^{(4)}\|_{L^\infty})}.
\end{equation}
The following bound \eqref{eq:common-local-comparison} is imposed only for \(|y|\geq h\), since \(\widetilde W_y\) is estimated on \([-h,h]\) using Taylor's theorem in Lemma~\ref{lem:bootstrap-local}. We suppose that there exist nonnegative continuous functions \(a_0,e_0,a_\infty,e_\infty\) on \([s_0,s_1]\) such that
\begin{equation}\label{eq:common-local-comparison}
\begin{split}
&\left|\frac{1+y^2}{y^2}\bigl(\mathcal R(y,s)-\mathcal A(y,s)U''(y)\bigr)-\frac{2\mathcal A(y,s)}{y^3}\widetilde W_y(y,s)\right|\\
&\quad\leq\frac{\eta_{\ell}}{4}|\widetilde W_y(y,s)|+(\lambda_{\ell}^{\mathrm{loc}}-1)K_{\ell}^{\mathrm{loc}}(y)\left\|\frac{1+y^2}{y^2}\widetilde W_y(y,s)\right\|_{L^\infty_y}\\
&\qquad +a_0(s)\frac{1+y^2}{y^2}|\widetilde W_y(y,s)|+e_0(s), \qquad |y|\geq h,
\end{split}
\end{equation}
and
\begin{equation}\label{eq:common-far-comparison}
\begin{split}
&\left|(1+|y|^{2/(2\ell+1)})\bigl(\mathcal R(y,s)-\mathcal A(y,s)U''(y)\bigr)+\frac{2|y|^{2/(2\ell+1)}\mathcal A(y,s)}{(2\ell+1)y}\widetilde W_y(y,s)\right|\\ &\quad\leq(\Lambda_{\ell}-1)K_{\ell}^{\mathrm{far}}(y)\left\|(1+|y|^{2/(2\ell+1)})\widetilde W_y(y,s)\right\|_{L^\infty_y} \\
&\qquad +a_\infty(s)(1+|y|^{2/(2\ell+1)})|\widetilde W_y(y,s)|+e_\infty(s), \qquad |y|\geq m,
\end{split}
\end{equation}
where \(K_{\ell}^{\mathrm{loc}}\) and \(K_{\ell}^{\mathrm{far}}\) are defined in \eqref{eq:local-kernel} and \eqref{eq:far-kernel}, respectively.

For the estimates of higher-order derivatives, we suppose that there exist nonnegative continuous functions \(a_j,e_j\) on \([s_0,s_1]\), \(j=2,3,4\), and constants \(\Gamma_3,\Gamma_4\geq0\) such that
\begin{equation}\label{eq:bootstrap-second-forcing}
\begin{split}
&\frac{\sqrt{1+y^2}}{|y|}\left|\mathcal R_2(y,s)-\frac{\mathcal A(y,s)}{y(1+y^2)}W_{yy}(y,s)\right|\\ &\quad\leq\left(\frac{c_{2,\ell}y^2}{4(1+y^2)}+a_2(s)\right)\frac{\sqrt{1+y^2}}{|y|}|W_{yy}(y,s)|+e_2(s),\qquad |y|\geq M^{-1},
\end{split}
\end{equation}
and
\begin{align}
|\mathcal R_3(y,s)|
&\leq \begin{aligned}[t]
&\left(\frac{(1+1/\ell)c_{1,\ell}y^2}{4(1+y^2)}+a_3(s)\right)|\partial_y^3W(y,s)|\\
&\quad+\Gamma_3|W_{yy}(y,s)|^2+e_3(s),\qquad |y|\geq\frac{\vartheta}{16M^{1/4}},
\end{aligned} \label{eq:bootstrap-third-forcing}\\
|\mathcal R_4(y,s)|
&\leq \begin{aligned}[t]
&\left(\frac18+a_4(s)\right)|\partial_y^4W(y,s)|\\
&\quad +\Gamma_4|W_{yy}(y,s)\partial_y^3W(y,s)|+e_4(s), \qquad y\in\R.
\end{aligned}\label{eq:bootstrap-fourth-forcing}
\end{align}

We impose the following time-integral bounds on \(a_j\) and \(e_j\):
\begin{equation}\label{eq:bootstrap-coefficient-integrals}
\int_{s_0}^{s_1}a_j(s)\,ds\leq\log\left(1+\frac{\vartheta}{8}\right),\qquad j=0,2,3,4,\infty,
\end{equation}
and
\begin{subequations}\label{eq:bootstrap-damped-forcing}
\begin{align}
&\sup_{s_0\leq s\leq s_1}\int_{s_0}^{s}e^{-\mu_0(s-\sigma)}e_0(\sigma)\,d\sigma\leq\frac{\vartheta N}{8},\label{eq:forcing-local-time}\\
&\sup_{s_0\leq s\leq s_1}\int_{s_0}^{s}e^{-\mu_2(s-\sigma)}e_2(\sigma)\,d\sigma\leq\frac{\vartheta M^{1/8}}8,\label{eq:forcing-second-time}\\
&\sup_{s_0\leq s\leq s_1}\int_{s_0}^{s}e^{-\mu_3(s-\sigma)}e_3(\sigma)\,d\sigma\leq\frac{\vartheta M^{3/4}}8,\label{eq:forcing-third-time}\\
&\sup_{s_0\leq s\leq s_1}\int_{s_0}^{s}e^{-(s-\sigma)/4}e_4(\sigma)\,d\sigma\leq\frac{\vartheta M}{8}, \label{eq:forcing-fourth-time}
\end{align}
\end{subequations}
where
\begin{equation}\label{eq:bootstrap-damping-rates}
\mu_0:=\frac{\eta_{\ell}h^2}{2(1+h^2)},\qquad \mu_2:=\frac{c_{2,\ell}}{2(1+M^2)},\qquad \mu_3:=\frac{(1+1/\ell)c_{1,\ell}\vartheta^2}{2(\vartheta^2+256M^{1/2})}.
\end{equation}
For \(e_\infty\), we require
\begin{equation}\label{eq:common-forcing-smallness}
\int_{s_0}^{s_1}e_\infty(s)\,ds\leq\frac{\vartheta B_{\ell}}{8}.
\end{equation}

\begin{remark}[A sufficient condition for the remainder bounds]\label{rem:sufficient-remainder-bounds}
Fix the profile \(U\) and the bootstrap parameters as in Section~\ref{sec:bootstrap-assumptions}. There exists \(C_*>0\), independent of \(s_0\) and \(s_1\), such that \eqref{eq:sufficient-remainder-bounds}--\eqref{eq:sufficient-remainder-smallness}, together with \eqref{eq:common-bootstrap-assumptions}, imply \eqref{eq:bootstrap-origin-forcing} and \eqref{eq:common-local-comparison}--\eqref{eq:common-forcing-smallness}.

By \eqref{eq:far-Usecond} and the smoothness of \(U\), we have \(\|(1+|y|^{2/(2\ell+1)})(1+|y|)U''\|_{L^\infty}<\infty\). Hence \eqref{eq:sufficient-remainder-bounds} gives
\begin{equation}\label{eq:sufficient-A-bounds}
|\mathcal A(0,s)|+\left\|(1+|y|^{2/(2\ell+1)})\mathcal A(y,s) U''(y)\right\|_{L_y^\infty}+\sup_{|y|\geq d}\left|\frac{\mathcal A(y,s)}{y}\right|\leq C_d e(s)
\end{equation}
for each \(d>0\), where \(C_d\) depends on \(d\) and the profile \(U\), but is independent of \(s_0\) and \(s_1\).

We first verify the pointwise bounds \eqref{eq:common-local-comparison}--\eqref{eq:bootstrap-fourth-forcing}. Using \eqref{eq:boot-local}, \eqref{eq:sufficient-remainder-bounds}, and \eqref{eq:sufficient-A-bounds} with \(d=h\), we obtain
\begin{equation*}
\begin{split}
&\sup_{|y|\geq h}\left|\frac{1+y^2}{y^2}\bigl(\mathcal R-\mathcal A U''\bigr)-\frac{2\mathcal A}{y^3}\widetilde W_y\right|\\
&\quad\leq(1+h^{-2})\bigl(\|\mathcal R\|_{L^\infty}+\|\mathcal A U''\|_{L^\infty}\bigr)+2N\sup_{|y|\geq h}\frac{|\mathcal A|}{|y|(1+y^2)} \leq Ce(s).
\end{split}
\end{equation*}
For \eqref{eq:common-far-comparison}, we use \eqref{eq:boot-far} and \eqref{eq:sufficient-A-bounds} with \(d=m\) to obtain
\begin{equation*}
\begin{split}
&\sup_{|y|\geq m}\left|(1+|y|^{2/(2\ell+1)})\bigl(\mathcal R-\mathcal A U''\bigr)+\frac{2|y|^{2/(2\ell+1)}\mathcal A}{(2\ell+1)y}\widetilde W_y\right|\\
&\quad\leq\|(1+|y|^{2/(2\ell+1)})\mathcal R\|_{L^\infty}+\|(1+|y|^{2/(2\ell+1)})\mathcal A U''\|_{L^\infty}+\frac{2B_\ell}{2\ell+1}\sup_{|y|\geq m}\left|\frac{\mathcal A}{y}\right| \leq Ce(s).
\end{split}
\end{equation*}
Here we used \(|y|^{2/(2\ell+1)}|\widetilde W_y|\leq B_\ell\), which follows from \eqref{eq:boot-far}. For \eqref{eq:bootstrap-second-forcing}, \eqref{eq:boot-second} gives
\begin{equation*}
\begin{split}
&\sup_{|y|\geq M^{-1}}\frac{\sqrt{1+y^2}}{|y|}\left|\mathcal R_2-\frac{\mathcal A}{y(1+y^2)}W_{yy}\right|\\
&\quad\leq\sqrt{1+M^2}\|\mathcal R_2\|_{L^\infty}+M^{1/8}\sup_{|y|\geq M^{-1}}\frac{|\mathcal A|}{|y|(1+y^2)} \leq Ce(s),
\end{split}
\end{equation*}
where the last inequality follows from \eqref{eq:sufficient-remainder-bounds} and \eqref{eq:sufficient-A-bounds} with \(d=M^{-1}\). The bounds \(|\mathcal R_3|\leq e(s)\) and \(|\mathcal R_4|\leq e(s)\) follow directly from \eqref{eq:sufficient-remainder-bounds}. Moreover, \eqref{eq:sufficient-remainder-bounds} and \eqref{eq:boot-fourth} give
\begin{equation*}
|\mathcal R_3(0,s)-\mathcal A(0,s)\partial_y^4W(0,s)|\leq(1+M)e(s).
\end{equation*}
Taking \(C\geq1+M\) sufficiently large to dominate the constants above, we may choose
\begin{equation*}
a_j=0,\qquad e_j=Ce,\qquad j=0,2,3,4,\infty,\qquad \Gamma_3=\Gamma_4=0.
\end{equation*}
Since \(\lambda_\ell^{\mathrm{loc}}>1\) and \(\Lambda_\ell>1\), the remaining terms on the right-hand sides of \eqref{eq:common-local-comparison}--\eqref{eq:bootstrap-fourth-forcing} are nonnegative. Thus the pointwise bounds hold with these choices.

It remains to verify \eqref{eq:bootstrap-origin-forcing}, \eqref{eq:bootstrap-coefficient-integrals}, \eqref{eq:bootstrap-damped-forcing}, and \eqref{eq:common-forcing-smallness}. The coefficient bounds \eqref{eq:bootstrap-coefficient-integrals} hold immediately by \(a_j=0\). For \eqref{eq:bootstrap-damped-forcing}, we use the estimate
\begin{equation}\label{eq:forcing-time-tests}
\sup_{s_0\leq s\leq s_1}\int_{s_0}^{s}e^{-\mu(s-\sigma)}e(\sigma)\,d\sigma\leq\min\left\{\int_{s_0}^{s_1}e(\sigma)\,d\sigma,\frac{\|e\|_{L^\infty([s_0,s_1])}}{\mu}\right\}
\end{equation}
for any \(\mu>0\). By \eqref{eq:sufficient-remainder-smallness}, each integral in \eqref{eq:bootstrap-damped-forcing} is therefore bounded by \(CC_*\). The bounds \eqref{eq:bootstrap-origin-forcing} and \eqref{eq:common-forcing-smallness} follow directly from \eqref{eq:sufficient-remainder-smallness}, since \(1+M\leq C\) and \(e_\infty=Ce\). Thus all the time-integral bounds hold if \(C_*\) is chosen so that
\begin{equation*}
CC_*\leq\frac{\vartheta}{8}\min\{N,M^{1/8},M^{3/4},M,B_\ell\}.
\end{equation*}
This proves the assertion.
\end{remark}

\subsubsection{Bootstrap estimates}

We now state the estimates that improve the bootstrap bounds \eqref{eq:common-bootstrap-assumptions} under suitable assumptions on the initial data; see Remark~\ref{rem:bootstrap-data}.

\begin{proposition} \label{prop:common-bootstrap}
 Let \(\ell>1/2\) be real, and choose \(m,\vartheta,N\) as in Section~\ref{sec:bootstrap-assumptions}. Fix \(\Gamma_3,\Gamma_4\geq0\). There exists \(M_0\geq1\), depending only on these constants and the corresponding profile $U$, such that the following statement holds for any \(M\geq M_0\).

Let \(W,\mathcal V,\mathcal R\) be smooth functions on \(\R\times[s_0,s_1]\) satisfying \eqref{eq:perturbed-slope-equation} and \eqref{eq:modulation-constraints-intro}, with \(\mathcal A\) determined by \eqref{eq:remainder-transport-split}. Assume that
\begin{equation}\label{eq:bootstrap-transport-growth}
|\mathcal V(y,s)|\leq C(1+|y|),\qquad (y,s)\in\R\times[s_0,s_1]
\end{equation}
for some constant \(C>0\). Assume also that
\begin{equation*}
\frac{1+y^2}{y^2}(W_y-U'),\quad (1+|y|^{2/(2\ell+1)})(W_y-U'),\quad \frac{\sqrt{1+y^2}}{y}W_{yy}, \quad \partial_y^3W,\quad \partial_y^4W
\end{equation*}
belong to $C([s_0,s_1];L^\infty(\R))$, with the quotients at the origin defined by continuity. If the bootstrap assumptions \eqref{eq:common-bootstrap-assumptions} and the remainder bounds \eqref{eq:bootstrap-origin-forcing} and \eqref{eq:common-local-comparison}--\eqref{eq:common-forcing-smallness} hold, then
\begin{subequations}\label{eq:common-bootstrap-conclusion}
\begin{align}
&|W_y(y,s)-U'(y)|\leq\frac{y^2}{1+y^2}\left(\left\|\frac{1+z^2}{z^2}(W_y-U')(z,s_0)\right\|_{L_z^\infty}+\vartheta N\right),\label{eq:improved-local}\\
&\bigl(1+|y|^{2/(2\ell+1)}\bigr)|W_y(y,s)-U'(y)|\notag\\
&\quad\leq\max\left\{\left\|\bigl(1+|z|^{2/(2\ell+1)}\bigr)(W_y-U')(z,s_0)\right\|_{L_z^\infty},N\bigl(1+m^{2/(2\ell+1)}\bigr)\right\}+\vartheta B_{\ell},\label{eq:improved-far}\\
&|W_{yy}(y,s)|\leq\frac{|y|}{\sqrt{1+y^2}}\left(\left\|\frac{\sqrt{1+z^2}}{z}W_{yy}(z,s_0)\right\|_{L_z^\infty}+\vartheta M^{1/8}\right),\label{eq:improved-second}\\
&|\partial_y^3W(0,s)-U^{(3)}(0)|\leq|\partial_y^3W(0,s_0)-U^{(3)}(0)|+\vartheta N,\label{eq:improved-origin}\\
&\|\partial_y^3W(\cdot,s)\|_{L^\infty}\leq\|\partial_y^3W(\cdot,s_0)\|_{L^\infty}+\vartheta M^{3/4},\label{eq:improved-third}\\
&\|\partial_y^4W(\cdot,s)\|_{L^\infty}\leq\|\partial_y^4W(\cdot,s_0)\|_{L^\infty}+\vartheta M \label{eq:improved-fourth}
\end{align}
\end{subequations}
for all \((y,s)\in\R\times[s_0,s_1]\).
\end{proposition}

\begin{proof}
The proof is deferred to Section~\ref{sec:proof-common-bootstrap}, where it is decomposed into a sequence of lemmas. The estimates \eqref{eq:improved-local}, \eqref{eq:improved-far}, \eqref{eq:improved-second}, \eqref{eq:improved-origin}, \eqref{eq:improved-third}, and \eqref{eq:improved-fourth} are proved in Lemmas~\ref{lem:bootstrap-local}, \ref{lem:bootstrap-far}, \ref{lem:bootstrap-second}, \ref{lem:bootstrap-origin}, \ref{lem:bootstrap-third}, and~\ref{lem:bootstrap-fourth}, respectively.
\end{proof}

\begin{remark}[Bootstrap improvement]\label{rem:bootstrap-data}
Proposition~\ref{prop:common-bootstrap} gives a strict improvement of \eqref{eq:common-bootstrap-assumptions} under suitable initial bounds. In addition to its hypotheses, suppose that
\begin{subequations}\label{eq:bootstrap-initial-bounds}
\begin{align}
&|W_y(y,s_0)-U'(y)|\leq(1-2\vartheta)\frac{Ny^2}{1+y^2},\label{eq:bootstrap-initial-local}\\
&|W_y(y,s_0)-U'(y)|\leq(1-2\vartheta)\frac{B_\ell}{1+|y|^{2/(2\ell+1)}},\label{eq:bootstrap-initial-far}\\
&|W_{yy}(y,s_0)|\leq(1-2\vartheta)\frac{M^{1/8}|y|}{\sqrt{1+y^2}},\label{eq:bootstrap-initial-second}\\
&|\partial_y^3W(0,s_0)-U^{(3)}(0)|\leq1-2\vartheta,\label{eq:bootstrap-initial-origin}\\
&\|\partial_y^3W(\cdot,s_0)\|_{L^\infty}\leq(1-2\vartheta)M^{3/4},\label{eq:bootstrap-initial-third}\\
&\|\partial_y^4W(\cdot,s_0)\|_{L^\infty}\leq(1-2\vartheta)M.\label{eq:bootstrap-initial-fourth}
\end{align}
\end{subequations}
for all \(y\in\R\). Then Proposition~\ref{prop:common-bootstrap} gives
\begin{equation}\label{eq:bootstrap-improved-bounds}
\begin{aligned}
&|W_y(y,s)-U'(y)|\leq(1-\vartheta)\frac{Ny^2}{1+y^2}, \\
&|W_y(y,s)-U'(y)|\leq(1-\vartheta)\frac{B_{\ell}}{1+|y|^{2/(2\ell+1)}}, \\
&|W_{yy}(y,s)|\leq(1-\vartheta)\frac{M^{1/8}|y|}{\sqrt{1+y^2}}, \\
&|\partial_y^3W(0,s)-U^{(3)}(0)|\leq1-\vartheta, \\
&\|\partial_y^3W(\cdot,s)\|_{L^\infty}\leq(1-\vartheta)M^{3/4}, \\
&\|\partial_y^4W(\cdot,s)\|_{L^\infty}\leq(1-\vartheta)M.
\end{aligned}
\end{equation}
for all \((y,s)\in\R\times[s_0,s_1]\). Thus all six bootstrap bounds are improved with a margin independent of \(s_1\).

The initial bounds \eqref{eq:bootstrap-initial-bounds} are one sufficient choice. The six initial margins may be chosen separately, provided that the bounds in \eqref{eq:common-bootstrap-conclusion} are strictly below the corresponding bootstrap thresholds.

For a given evolution equation, the initial conditions corresponding to \eqref{eq:bootstrap-initial-bounds}, together with the matching conditions at the origin and the conditions needed for the model-specific estimates, identify a class of smooth initial data for which the bootstrap argument yields stability estimates around \(U\). For \eqref{eq:b-family}, we construct a non-empty relatively open class of such initial data in Appendix~\ref{subsec:b-data}. Related constructions for other models are given in Appendix~A of \cite{KKS-rSV} and Section~5.2 of \cite{KKKS-Novikov}.
\end{remark}

\subsubsection{From stability estimates to sharp regularity at blow-up}

Proposition~\ref{prop:common-bootstrap} provides stability estimates for \(W\) around \(U\) on a bootstrap interval. In particular, \eqref{eq:improved-far} and \eqref{eq:improved-second}, together with \eqref{eq:far-Uprime}, imply \eqref{eq:regularity-profile-bounds}. Assuming \eqref{eq:regularity-profile-bounds}, we derive sharp regularity and blow-up rates in the original variables under the normalization \eqref{eq:modulation-constraints-intro} and the modulation assumptions in Proposition~\ref{prop:sharp-regularity}.

\begin{proposition}\label{prop:sharp-regularity}
Let \(\ell>1/2\), and let \(W\) be a smooth function on \(\R\times[s_0,\infty)\) satisfying \eqref{eq:modulation-constraints-intro} and \eqref{eq:regularity-profile-bounds}. Let \(t_0<T_*<\infty\), and let \(\tau,\kappa,\xi\) be modulation functions on \([t_0,T_*)\), with \(\tau,\xi\) continuous and \(\kappa\) bounded. Assume that \(\tau(t)>t\) and \(s(t):=-\log(\tau(t)-t)\geq s_0\) for all $t \in [t_0,T_*)$, and that
\begin{equation}\label{eq:regularity-modulation}
\tau(t)-t\sim T_*-t,\qquad \xi(t)\longrightarrow x_*\quad\text{as }t\uparrow T_*
\end{equation}
for some \(x_*\in\R\). Define \(w\) on \(\R\times[t_0,T_*)\) by \eqref{eq:modulated-variables} and \eqref{eq:modulated-unknown}, and let \(\gamma_\ell\) be given by \eqref{eq:cusp-exponent}. Then the following holds.

\begin{enumerate}[(i)]
\item For any bounded open set \(\Omega\subset\R\),
\begin{equation}\label{eq:regularity-local-boundedness}
\sup_{t_0\leq t<T_*} \|w(\cdot,t)\|_{L^\infty(\Omega)} < \infty,
\end{equation}
whereas
\begin{equation}\label{eq:regularity-gradient-blowup}
\lim_{t\uparrow T_*}\|w_x(\cdot,t)\|_{L^\infty(\R)}=\infty.
\end{equation}
\item For any bounded open set \(\Omega\subset\R\) and any \(\alpha\in(0,\gamma_\ell]\),
\begin{equation}\label{eq:regularity-critical}
\sup_{t_0\leq t<T_*}[w(\cdot,t)]_{C^\alpha(\Omega)}<\infty.
\end{equation}

\item For any bounded open set \(\Omega\subset\R\) and any \(\alpha\in(\gamma_\ell,1]\),
\begin{equation}\label{eq:regularity-holder-dichotomy}
\begin{cases}
\displaystyle \lim_{t\uparrow T_*}[w(\cdot,t)]_{C^\alpha(\Omega)}=\infty,&\text{if }x_*\in\Omega,\\
\displaystyle \limsup_{t\uparrow T_*}[w(\cdot,t)]_{C^\alpha(\Omega)}<\infty,&\text{if }x_*\notin\overline{\Omega}.
\end{cases}
\end{equation}
If \(x_*\in\Omega\), the temporal blow-up rate is given by
\begin{equation}\label{eq:regularity-rate}
[w(\cdot,t)]_{C^\alpha(\Omega)}\sim(T_*-t)^{-((2\ell+1)\alpha-(2\ell-1))/2}
\end{equation}
for all \(t\) sufficiently close to \(T_*\).
\end{enumerate}
\end{proposition}

Assertions \emph{(ii)} and \emph{(iii)} show that the H\"older exponent \(\gamma_\ell\) is sharp in the sense that, on any bounded open set containing \(x_*\), the \(C^{\gamma_\ell}\) seminorm of \(w(\cdot,t)\) remains uniformly bounded up to \(T_*\), whereas the \(C^\alpha\) seminorm blows up as \(t\uparrow T_*\) for any \(\alpha\in(\gamma_\ell,1]\).

\begin{proof}
We prove the assertions \emph{(i)}--\emph{(iii)} separately. In what follows, \(s=s(t)\).

\emph{(i)} Differentiating \eqref{eq:modulated-unknown} with respect to \(x\) and using \eqref{eq:regularity-profile-bounds}, we obtain
\begin{equation*}
|w_x(x,t)|=\frac{|W_y(y(x,t),s)|}{\tau(t)-t}\leq\frac{C}{\tau(t)-t+|x-\xi(t)|^{2/(2\ell+1)}}
\end{equation*}
for all \((x,t)\in\R\times[t_0,T_*)\). In particular, by \eqref{eq:modulation-constraints-intro}, we have
\begin{equation*}
w_x(\xi(t),t)=\frac{W_y(0,s)}{\tau(t)-t}=-\frac{\ell}{\tau(t)-t}.
\end{equation*}
Since \(\tau(t)>t\), \eqref{eq:regularity-modulation} gives \(\tau(t)-t\to0^+\) as \(t\uparrow T_*\). Hence
\begin{equation*}
\|w_x(\cdot,t)\|_{L^\infty(\R)}\geq|w_x(\xi(t),t)|=\frac{\ell}{\tau(t)-t}\longrightarrow\infty\quad\text{as }t\uparrow T_*.
\end{equation*}
This proves \eqref{eq:regularity-gradient-blowup}.

To prove \eqref{eq:regularity-local-boundedness}, we first note from \eqref{eq:regularity-profile-bounds} that, for any \(y,\widetilde y\in\R\),
\begin{equation}\label{eq:regularity-rescaled-increments}
|W(y,s)-W(\widetilde y,s)|=\left|\int_{\widetilde y}^y W_y(\widehat y,s)\,d\widehat y\right|\leq C\left|\int_{\widetilde y}^y\frac{d\widehat y}{1+|\widehat y|^{2/(2\ell+1)}}\right|\leq C|y-\widetilde y|^{\gamma_\ell}.
\end{equation}
By the change of variables \eqref{eq:modulated-variables}--\eqref{eq:modulated-unknown} together with \eqref{eq:cusp-exponent}, for any \(x,\widetilde x\in\R\), we have
\begin{equation}\label{eq:regularity-increments}
|w(x,t)-w(\widetilde x,t)|=e^{-(2\ell-1)s/2}|W(y,s)-W(\widetilde y,s)|\leq C|x-\widetilde x|^{\gamma_\ell}
\end{equation}
for all \(t\in[t_0,T_*)\), where \(y=y(x,t)\) and \(\widetilde y=y(\widetilde x,t)\). Since \(\xi\) is continuous on \([t_0,T_*)\) and converges to \(x_*\) as \(t\uparrow T_*\), it is bounded on \([t_0,T_*)\). By \eqref{eq:modulated-unknown} and \eqref{eq:modulation-constraints-intro}, \(w(\xi(t),t)=\kappa(t)\). Thus, for any bounded open set \(\Omega\subset\R\), \eqref{eq:regularity-increments} gives
\begin{equation*}
|w(x,t)|\leq|\kappa(t)|+|w(x,t)-w(\xi(t),t)|\leq|\kappa(t)|+C|x-\xi(t)|^{\gamma_\ell}\leq C
\end{equation*}
for all \(x\in\Omega\) and \(t\in[t_0,T_*)\). This proves \eqref{eq:regularity-local-boundedness}.

\medskip
\emph{(ii)} Let \(\Omega\subset\R\) be a bounded open set, and let \(\alpha\in(0,\gamma_\ell]\). For any \(x,\widetilde x\in\Omega\) with \(x\ne\widetilde x\), \eqref{eq:regularity-increments} gives
\begin{equation*}
\frac{|w(x,t)-w(\widetilde x,t)|}{|x-\widetilde x|^\alpha}\leq C|x-\widetilde x|^{\gamma_\ell-\alpha}\leq C,\qquad t\in[t_0,T_*).
\end{equation*}
Taking the supremum over \(x,\widetilde x\in\Omega\) with \(x\ne\widetilde x\) and then over \(t\in[t_0,T_*)\), we obtain \eqref{eq:regularity-critical}.

\medskip
\emph{(iii)} Let \(\Omega\subset\R\) be a bounded open set, and let \(\alpha\in(\gamma_\ell,1]\). We note from \eqref{eq:modulated-variables} and \eqref{eq:modulated-unknown} that
\begin{equation}\label{eq:regularity-holder-scaling}
\frac{|w(x,t)-w(\widetilde x,t)|}{|x-\widetilde x|^\alpha}=e^{((2\ell+1)\alpha-(2\ell-1))s/2}\frac{|W(y,s)-W(\widetilde y,s)|}{|y-\widetilde y|^\alpha}
\end{equation}
for any \(x,\widetilde x\in\R\) with \(x\ne\widetilde x\), where \(y=y(x,t)\) and \(\widetilde y=y(\widetilde x,t)\).

Suppose first that \(x_*\in\Omega\). By \eqref{eq:modulation-constraints-intro} and the bound for \(W_{yy}\) in \eqref{eq:regularity-profile-bounds}, we have
\begin{equation*}
|W_y(y,s)+\ell|=\left|\int_0^y W_{yy}(z,s)\,dz\right|\leq C_0|y|,\qquad |y|\leq1.
\end{equation*}
Choose \(\eta\in(0,1]\) sufficiently small that \(C_0\eta\leq\ell/2\). Then
\begin{equation*}
W_y(y,s)\leq-\frac{\ell}{2},\qquad 0\leq y\leq\eta.
\end{equation*}
Fix \(y_0\in(0,\eta)\). Integrating over \([0,y_0]\), we obtain
\begin{equation}\label{eq:regularity-profile-increment-lower}
|W(y_0,s)-W(0,s)|=-\int_0^{y_0}W_y(y,s)\,dy\geq\frac{\ell y_0}{2}.
\end{equation}
We define the corresponding points in the original variables by
\begin{equation*}
x_0(t):=\xi(t)+e^{-(2\ell+1)s/2}y_0,\qquad \widetilde x_0(t):=\xi(t).
\end{equation*}
Then, by \eqref{eq:regularity-modulation}, \(x_0(t),\widetilde x_0(t)\to x_*\) as \(t\uparrow T_*\). Since \(x_*\in\Omega\) and \(\Omega\) is open, \(x_0(t),\widetilde x_0(t)\in\Omega\) for all \(t\) sufficiently close to \(T_*\). Using \eqref{eq:regularity-holder-scaling} and \eqref{eq:regularity-profile-increment-lower}, we obtain
\begin{equation}\label{eq:regularity-holder-lower}
\begin{aligned}
[w(\cdot,t)]_{C^\alpha(\Omega)}&\geq\frac{|w(x_0(t),t)-w(\widetilde x_0(t),t)|}{|x_0(t)-\widetilde x_0(t)|^\alpha} \\
&=e^{((2\ell+1)\alpha-(2\ell-1))s/2}\frac{|W(y_0,s)-W(0,s)|}{y_0^\alpha}\\
&\geq c_0e^{((2\ell+1)\alpha-(2\ell-1))s/2},
\end{aligned}
\end{equation}
where \(c_0:=\ell y_0^{1-\alpha}/2>0\). Since \(\alpha>\gamma_\ell\) and \(s\to\infty\) as \(t\uparrow T_*\), we have
\begin{equation*}
[w(\cdot,t)]_{C^\alpha(\Omega)}\longrightarrow\infty\quad\text{as }t\uparrow T_*.
\end{equation*}
This proves the first assertion in \eqref{eq:regularity-holder-dichotomy}.

On the other hand, suppose that \(x_*\notin\overline\Omega\). We prove that $\limsup_{t\uparrow T_*}[w(\cdot,t)]_{C^\alpha(\Omega)}<\infty$. Set
\begin{equation*}
d_*:=\inf_{x\in\Omega}|x-x_*|>0.
\end{equation*}
By the convergence \(\xi(t)\to x_*\), we may choose \(t_1\in(t_0,T_*)\) such that
\begin{equation*}
|\xi(t)-x_*|\leq\frac{d_*}{2},\qquad t\in[t_1,T_*).
\end{equation*}
For this choice of \(t_1\), we have
\begin{equation*}
|x-\xi(t)|\geq|x-x_*|-|\xi(t)-x_*|\geq\frac{d_*}{2},\qquad x\in\Omega,\quad t\in[t_1,T_*),
\end{equation*}
and hence
\begin{equation*}
|y(x,t)|=e^{(2\ell+1)s/2}|x-\xi(t)|\geq\frac{d_*}{2}e^{(2\ell+1)s/2},\qquad x\in\Omega,\quad t\in[t_1,T_*).
\end{equation*}
Fix \(t\in[t_1,T_*)\) and \(x,\widetilde x\in\Omega\) with \(x\ne\widetilde x\). We estimate the right-hand side of \eqref{eq:regularity-holder-scaling} according to whether \(y\) and \(\widetilde y\) have the same or opposite signs.

If \(y\) and \(\widetilde y\) have the same sign, then
\begin{equation*}
|\widehat y|\geq\frac{d_*}{2}e^{(2\ell+1)s/2}
\end{equation*}
for every \(\widehat y\) between \(y\) and \(\widetilde y\). Since \(\gamma_\ell-\alpha<0\), \eqref{eq:regularity-profile-bounds} gives
\begin{equation*}
|W_y(\widehat y,s)| \leq C|\widehat y|^{-2/(2\ell+1)} =C|\widehat y|^{\alpha-1}|\widehat y|^{\gamma_\ell-\alpha} \leq Ce^{-((2\ell+1)\alpha-(2\ell-1))s/2}|\widehat y|^{\alpha-1}.
\end{equation*}
Combining this estimate with \eqref{eq:regularity-holder-scaling}, we obtain
\begin{equation*}
\begin{aligned}
\frac{|w(x,t)-w(\widetilde x,t)|}{|x-\widetilde x|^\alpha} &=e^{((2\ell+1)\alpha-(2\ell-1))s/2}\frac{|W(y,s)-W(\widetilde y,s)|}{|y-\widetilde y|^\alpha}\\
&\leq\frac{C}{|y-\widetilde y|^\alpha}\left|\int_{\widetilde y}^y|\widehat y|^{\alpha-1}\,d\widehat y\right| \leq C\frac{\bigl||y|^\alpha-|\widetilde y|^\alpha\bigr|}{|y-\widetilde y|^\alpha}\leq C,
\end{aligned}
\end{equation*}
where we used
\begin{equation*}
\bigl||y|^\alpha-|\widetilde y|^\alpha\bigr|\leq|y-\widetilde y|^\alpha,\qquad 0<\alpha\leq1.
\end{equation*}

If \(y\) and \(\widetilde y\) have opposite signs, then \(x\) and \(\widetilde x\) lie on opposite sides of \(\xi(t)\), and hence
\begin{equation*}
|x-\widetilde x|=|x-\xi(t)|+|\widetilde x-\xi(t)|\geq d_*.
\end{equation*}
Therefore,
\begin{equation*}
|y-\widetilde y|=e^{(2\ell+1)s/2}|x-\widetilde x|\geq d_*e^{(2\ell+1)s/2}.
\end{equation*}
Using \eqref{eq:regularity-rescaled-increments} and \eqref{eq:regularity-holder-scaling}, we obtain
\begin{equation*}
\begin{aligned}
\frac{|w(x,t)-w(\widetilde x,t)|}{|x-\widetilde x|^\alpha} &\leq Ce^{((2\ell+1)\alpha-(2\ell-1))s/2}|y-\widetilde y|^{\gamma_\ell-\alpha}\\
&\leq Ce^{((2\ell+1)\alpha-(2\ell-1))s/2}\left( d_*e^{(2\ell+1)s/2}\right)^{\gamma_\ell-\alpha} \leq C.
\end{aligned}
\end{equation*}

Combining the estimates in the two cases and taking the supremum over \(x,\widetilde x\in\Omega\) with \(x\ne\widetilde x\), we obtain
\begin{equation*}
[w(\cdot,t)]_{C^\alpha(\Omega)}\leq C,\qquad t\in[t_1,T_*).
\end{equation*}
Hence
\begin{equation*}
\limsup_{t\uparrow T_*}[w(\cdot,t)]_{C^\alpha(\Omega)}<\infty.
\end{equation*}
This proves the second assertion in \eqref{eq:regularity-holder-dichotomy}.

Finally, we prove \eqref{eq:regularity-rate} under the assumption \(x_*\in\Omega\). By \eqref{eq:regularity-profile-bounds} and \(\alpha>\gamma_\ell\),
\begin{equation*}
|W_y(y,s)| \leq\frac{C}{1+|y|^{2/(2\ell+1)}} \leq C(1+y^2)^{-1/(2\ell+1)} \leq C(1+y^2)^{(\alpha-1)/2}.
\end{equation*}
Hence, for any \(y,\widetilde y\in\R\),
\begin{equation*}
|W(y,s)-W(\widetilde y,s)|\leq\left|\int_{\widetilde y}^y|W_y(\widehat y,s)|\,d\widehat y\right|\leq C\left|\int_{\widetilde y}^y(1+\widehat y^2)^{(\alpha-1)/2}\,d\widehat y\right|\leq C|y-\widetilde y|^\alpha.
\end{equation*}
Using \eqref{eq:regularity-holder-scaling}, for any \(x,\widetilde x\in\Omega\) with \(x\ne\widetilde x\), we have
\begin{equation*}
\frac{|w(x,t)-w(\widetilde x,t)|}{|x-\widetilde x|^\alpha}=e^{((2\ell+1)\alpha-(2\ell-1))s/2}\frac{|W(y,s)-W(\widetilde y,s)|}{|y-\widetilde y|^\alpha}\leq Ce^{((2\ell+1)\alpha-(2\ell-1))s/2}.
\end{equation*}
Taking the supremum over these points and combining the resulting upper bound with \eqref{eq:regularity-holder-lower}, we obtain
\begin{equation*}
[w(\cdot,t)]_{C^\alpha(\Omega)}\sim e^{((2\ell+1)\alpha-(2\ell-1))s/2}
\end{equation*}
for all \(t\) sufficiently close to \(T_*\). Since \(e^{-s}=\tau(t)-t\), \eqref{eq:regularity-modulation} gives
\begin{equation*}
e^{((2\ell+1)\alpha-(2\ell-1))s/2}=(\tau(t)-t)^{-((2\ell+1)\alpha-(2\ell-1))/2}\sim(T_*-t)^{-((2\ell+1)\alpha-(2\ell-1))/2}.
\end{equation*}
This proves \eqref{eq:regularity-rate} and completes the proof.
\end{proof}

\subsubsection{A general theorem for singularity formation} \label{sec:common-singularity}

We now combine Propositions~\ref{prop:common-bootstrap} and~\ref{prop:sharp-regularity} with the assumptions on local existence, modulation, and continuation for a given evolution equation. This provides the precise version of Theorem~\ref{thm:common-singularity-intro}.

\begin{theorem}[Singularity formation and sharp regularity]\label{thm:common-singularity}
Let \(\ell>1/2\), and let \(U=U_\beta\) be the profile constructed in Proposition~\ref{prop:profile-family}, with \(\beta\) chosen as in \eqref{eq:normalized-profile-choice}. Let \(\Theta_\ell\) and \(B_\ell\) be given by \eqref{eq:theta-ell} and \eqref{eq:B-choice}, respectively. Fix \(m,\vartheta,N,\Gamma_3,\Gamma_4\) as in Proposition~\ref{prop:common-bootstrap}, with
\begin{equation}\label{eq:common-singularity-theta}
0<\vartheta<\frac12\left(1-\frac{\Theta_\ell\beta^{-1/(2\ell+1)}}{B_\ell}\right).
\end{equation}
Let \(M_0\) be as in Proposition~\ref{prop:common-bootstrap}, and fix \(M\geq M_0\). Consider an initial-value problem for a scalar evolution equation on \(\R\), with unknown \(w\) and smooth initial data prescribed at \(t=t_0\). Suppose that the following assumptions hold.

\begin{enumerate}[(i)]
\item \emph{(Local existence and self-similar formulation)}
The initial-value problem admits a local solution \(w\) and \(C^1\) modulation functions \(\tau,\kappa,\xi\) such that \(\tau(t)>t\) and
\begin{equation*}
s(t):=-\log(\tau(t)-t)
\end{equation*}
is strictly increasing. Set \(s_0:=s(t_0)\), and define \(W\) by \eqref{eq:modulated-variables}--\eqref{eq:modulated-unknown}. There exist functions \(\mathcal A\) and \(\mathcal R\) such that, with
\begin{equation*}
\mathcal V(y,s):=\frac{2\ell+1}{2}y+W(y,s)+\mathcal A(y,s),
\end{equation*}
the function \(W\) satisfies \eqref{eq:perturbed-slope-equation} and \eqref{eq:modulation-constraints-intro}. On every finite interval \([s_0,s_1]\) contained in the self-similar time interval of existence, the functions \(W,\mathcal V,\mathcal R\) are smooth on \(\R\times[s_0,s_1]\) and satisfy
\begin{equation}\label{eq:common-singularity-velocity}
|\mathcal V(y,s)|\leq C (1+|y|),\qquad (y,s)\in\R\times[s_0,s_1],
\end{equation}
for some \(C>0\). Moreover,
\begin{equation*}
\frac{1+y^2}{y^2}(W_y-U'),\quad (1+|y|^{2/(2\ell+1)})(W_y-U'),\quad \frac{\sqrt{1+y^2}}{y}W_{yy}, \quad \partial_y^3W,\quad \partial_y^4W
\end{equation*}
belong to \(C([s_0,s_1];L^\infty(\R))\), with the quotients at the origin defined by continuity.

\item \emph{(Initial bounds)}
At \(s=s_0\), the initial data $W(y,s_0)$ satisfy \eqref{eq:bootstrap-initial-bounds}.

\item \emph{(Bounds on the model-specific perturbative terms)}
Whenever the bootstrap assumptions \eqref{eq:common-bootstrap-assumptions} hold on a finite interval \([s_0,s_1]\), the bounds \eqref{eq:bootstrap-origin-forcing} and \eqref{eq:common-local-comparison}--\eqref{eq:common-forcing-smallness} hold on $[s_0,s_1]$.

\item \emph{(Modulation estimates)}
There exist constants $c>0$ and $C>0$ such that, for every finite interval \([s_0,s_1]\) on which \eqref{eq:common-bootstrap-assumptions} holds,
\begin{equation}\label{eq:common-singularity-modulation}
c\leq1-\dot\tau(t)\leq C,\qquad
|\kappa(t)|\leq C,\qquad
\int_{t_0}^{t}|\dot\xi(t')|\,dt'\leq C,
\qquad t_0\leq t\leq t(s_1),
\end{equation}
where \(t=t(s)\) is the inverse of \(s=s(t)\).

\item \emph{(Continuation)}
The local solution has the following continuation property. Suppose that \(S<\infty\), that the solution and its modulation functions are defined for \(s_0\leq s<S\), and that the self-similar formulation and regularity properties in \emph{(i)} hold on every compact subinterval of \([s_0,S)\). If the bootstrap assumptions \eqref{eq:common-bootstrap-assumptions} hold for \(s_0\leq s<S\), then there exists \(\varepsilon>0\) such that the solution and its modulation functions extend to \(s_0\leq s<S+\varepsilon\), with the self-similar formulation and the regularity properties in \emph{(i)} preserved.
\end{enumerate}
Then \(W\) is defined for all self-similar times \(s\geq s_0\), and the bounds in \eqref{eq:bootstrap-improved-bounds} hold on \(\R\times[s_0,\infty)\). In the original variables, there exist \(T_*\in(t_0,\infty)\) and \(x_*\in\R\) such that \(w\) extends as a solution to \(\R\times[t_0,T_*)\) and
\begin{equation}\label{eq:common-singularity-time-center}
\tau(t)-t\sim T_*-t,\qquad \xi(t)\longrightarrow x_*\quad\text{as }t\uparrow T_*.
\end{equation}
Moreover, \(w\) satisfies conclusions \emph{(i)}--\emph{(iii)} of Proposition~\ref{prop:sharp-regularity}, with these values of \(T_*\) and \(x_*\).
\end{theorem}

The proof of Theorem~\ref{thm:common-singularity} is given in Section~\ref{sec:proof-common-singularity}. The same argument applies to a system of evolution equations if a scalar component \(w\) admits the self-similar formulation and satisfies the initial, perturbative, and modulation assumptions above, while local existence and continuation are established for the full system. The perturbative terms may depend on the remaining components, and the estimates needed to control and continue those components must be verified for the particular system.

\begin{remark}[Choice and compatibility of the parameters]\label{rem:common-singularity-parameters}
In Theorem~\ref{thm:common-singularity}, once \(\ell>1/2\) is fixed, the profile \(U=U_\beta\), with \(\beta\) given by \eqref{eq:normalized-profile-choice}, and the asymptotic coefficient \(\Theta_\ell\beta^{-1/(2\ell+1)}\) of \(U'\) are determined. We also fix the constants in Propositions~\ref{prop:damping-inequalities} and~\ref{prop:local-integral}, the constants \(B_\ell,\Lambda_\ell\) in \eqref{eq:B-choice}, and a threshold \(m_0(\ell)\) as in Proposition~\ref{prop:far-field}. We may then choose \(m\geq\max\{1,m_0(\ell)\}\) and \(\vartheta\) satisfying \eqref{eq:common-singularity-theta}. These choices are independent.  Since \(\Theta_\ell\beta^{-1/(2\ell+1)}<B_\ell\), the admissible interval for \(\vartheta\) is nonempty, and \eqref{eq:common-singularity-theta} gives \((1-2\vartheta)B_\ell>\Theta_\ell\beta^{-1/(2\ell+1)}\). In particular, the initial far-field bound in \eqref{eq:bootstrap-initial-bounds} is compatible with initial slopes whose decay is faster than that of \(U'\).

After \(m\) and \(\vartheta\) have been fixed, \(N>0\) is chosen to satisfy \eqref{eq:common-smallness}. These conditions are simultaneous upper bounds with positive right-hand sides, and hence such a choice is always possible. The parameter \(N\) may depend on \(m\) and \(\vartheta\); in particular, increasing \(m\) may require \(N\) to be smaller.

We next fix \(\Gamma_3,\Gamma_4\geq0\) and choose \(M\geq M_0\), where \(M_0\) is given by Proposition~\ref{prop:common-bootstrap} and depends only on the parameters fixed above and the profile \(U\), independently of \(s_1\). Thus an admissible order of selection is
\begin{equation*}
\ell\ \longrightarrow\ \text{profile and comparison constants}\ \longrightarrow\ (m,\vartheta)\ \longrightarrow\ N\ \longrightarrow\ (\Gamma_3,\Gamma_4)\ \longrightarrow\ M.
\end{equation*}
This establishes the compatibility of the numerical restrictions in the common argument. The construction of a nonempty class of initial data and the verification of assumptions \emph{(i)}--\emph{(v)} remain model-dependent. All constants required uniformly on bootstrap intervals, including those in \eqref{eq:common-singularity-modulation}, are independent of \(s_1\), whereas the constant \(C\) in \eqref{eq:common-singularity-velocity} may depend on the finite interval.
\end{remark}

\subsection{Proof of Proposition~\ref{prop:common-bootstrap}} \label{sec:proof-common-bootstrap}

We prove the estimates \eqref{eq:improved-local}--\eqref{eq:improved-fourth} in Proposition~\ref{prop:common-bootstrap}. We first estimate the third derivative at the origin and the fourth derivative on the whole line, and then prove the weighted estimates for \(W_y-U'\) and \(W_{yy}\). The global third-derivative estimate and the far-field estimate complete the proof. Throughout the proof, we assume the hypotheses of Proposition~\ref{prop:common-bootstrap} and choose \(M_0\) sufficiently large as needed, independently of \(s_1\).

We first prove \eqref{eq:improved-origin} and \eqref{eq:improved-fourth}.

\begin{lemma}\label{lem:bootstrap-origin}
Under the assumptions of Proposition~\ref{prop:common-bootstrap}, the estimate \eqref{eq:improved-origin} holds; that is,
\begin{equation*}
|\partial_y^3W(0,s)-U^{(3)}(0)|\leq|\partial_y^3W(0,s_0)-U^{(3)}(0)|+\vartheta N
\end{equation*}
for all \(s\in[s_0,s_1]\).
\end{lemma}

\begin{proof}
Evaluating \eqref{eq:third-derivative-common} at \(y=0\) and using \eqref{eq:modulation-constraints-intro}, we obtain
\begin{equation}\label{eq:third-derivative-origin}
\partial_s\partial_y^3W(0,s)=\mathcal R_3(0,s)-\mathcal A(0,s)\partial_y^4W(0,s).
\end{equation}
Integrating \eqref{eq:third-derivative-origin} and applying \eqref{eq:bootstrap-origin-forcing}, we have
\begin{equation}\label{eq:origin-increment-bound}
\left|\partial_y^3W(0,s)-\partial_y^3W(0,s_0)\right|\leq\int_{s_0}^{s}\left|\mathcal R_3(0,s')-\mathcal A(0,s')\partial_y^4W(0,s')\right|\,ds'\leq\frac{\vartheta N}{8}
\end{equation}
for all $s\in[s_0,s_1]$. By the triangle inequality and \eqref{eq:origin-increment-bound},
\begin{equation*}
|\partial_y^3W(0,s)-U^{(3)}(0)|
\leq|\partial_y^3W(0,s_0)-U^{(3)}(0)|+\frac{\vartheta N}{8}.
\end{equation*}
Since \(\vartheta N/8\leq\vartheta N\), this proves \eqref{eq:improved-origin}.
\end{proof}

\begin{lemma}\label{lem:bootstrap-fourth}
Under the assumptions of Proposition~\ref{prop:common-bootstrap}, the estimate \eqref{eq:improved-fourth} holds; that is,
\begin{equation*}
\|\partial_y^4W(\cdot,s)\|_{L^\infty}\leq\|\partial_y^4W(\cdot,s_0)\|_{L^\infty}+\vartheta M
\end{equation*}
for all \(s\in[s_0,s_1]\).
\end{lemma}

\begin{proof}
We rewrite \eqref{eq:fourth-derivative-common} in the form
\begin{equation}\label{eq:fourth-transport-proof}
\partial_s\partial_y^4W+D_4^W\partial_y^4W+\mathcal V\partial_y^5W=F_4^W,
\end{equation}
where \(\mathcal V\) is defined in \eqref{eq:remainder-transport-split}, and
\begin{equation*}
\begin{aligned}
D_4^W(y,s)&:=\frac{6\ell+5}{2}+\left(3+\frac2\ell\right)W_y(y,s),\\
F_4^W(y,s)&:=-\left(4+\frac6\ell\right)W_{yy}(y,s)\partial_y^3W(y,s)+\mathcal R_4(y,s).
\end{aligned}
\end{equation*}
By \eqref{eq:bootstrap-profile-damping}, \eqref{eq:boot-local}, and \eqref{eq:N-local},
\begin{equation}\label{eq:bootstrap-slope-lower}
\ell+W_y(y,s)\geq\ell+U'(y)-|W_y(y,s)-U'(y)|\geq(c_{1,\ell}-N)\frac{y^2}{1+y^2}\geq\frac{c_{1,\ell}}2\frac{y^2}{1+y^2}\geq0.
\end{equation}
It follows that
\begin{equation}\label{eq:fourth-damping-lower}
D_4^W(y,s)=\frac12+\left(3+\frac2\ell\right)(\ell+W_y(y,s))\geq\frac12.
\end{equation}
Next, using \eqref{eq:bootstrap-fourth-forcing}, \eqref{eq:boot-second}, and \eqref{eq:boot-third}, we estimate \(F_4^W\) as follows:
\begin{equation}\label{eq:fourth-forcing-bound}
\begin{split}
|F_4^W(y,s)|&\leq\left(4+\frac6\ell+\Gamma_4\right)|W_{yy}(y,s)\partial_y^3W(y,s)|+\left(\frac18+a_4(s)\right)|\partial_y^4W(y,s)|+e_4(s)\\
&\leq\left(4+\frac6\ell+\Gamma_4\right)M^{7/8}+\left(\frac18+a_4(s)\right)|\partial_y^4W(y,s)|+e_4(s).
\end{split}
\end{equation}
Combining \eqref{eq:fourth-damping-lower} and \eqref{eq:fourth-forcing-bound} with \eqref{eq:fourth-transport-proof}, we have
\begin{equation}\label{eq:fourth-differential-bound}
(\partial_s+\mathcal V\partial_y)|\partial_y^4W|+\left(\frac38-a_4(s)\right)|\partial_y^4W|\leq\left(4+\frac6\ell+\Gamma_4\right)M^{7/8}+e_4(s).
\end{equation}

Let \(\psi(y,s)\) denote the characteristic flow defined by
\begin{equation*}
\partial_s\psi(y,s)=\mathcal V(\psi(y,s),s),\qquad \psi(y,s_0)=y.
\end{equation*}
By the smoothness of \(\mathcal V\) and \eqref{eq:bootstrap-transport-growth}, the flow is defined on \([s_0,s_1]\), and \(y\mapsto\psi(y,s)\) is onto for each \(s\). Integrating \eqref{eq:fourth-differential-bound} along \(\psi\) and taking the supremum in \(y\), we obtain
\begin{equation*}
\begin{split}
\|\partial_y^4W(\cdot,s)\|_{L^\infty} &\leq e^{\int_{s_0}^sa_4(s')\,ds'} \bigg( \|\partial_y^4W(\cdot,s_0)\|_{L^\infty}e^{-3(s-s_0)/8} \\
& \qquad \qquad \qquad \quad +\int_{s_0}^{s}e^{-3(s-s')/8}\left[\left(4+\frac6\ell+\Gamma_4\right)M^{7/8}+e_4(s')\right] \, ds' \bigg)\\
&\leq e^{\int_{s_0}^sa_4(s')\,ds'} \bigg(\|\partial_y^4W(\cdot,s_0)\|_{L^\infty}+\frac83\left(4+\frac6\ell+\Gamma_4\right)M^{7/8} \\
& \qquad \qquad \qquad \quad +\int_{s_0}^{s}e^{-(s-s')/4}e_4(s')\,ds' \bigg).
\end{split}
\end{equation*}
We choose \(M\) sufficiently large so that
\begin{equation*}
64\left(4+\frac6\ell+\Gamma_4\right)M^{7/8}\leq\vartheta M.
\end{equation*}
Then \eqref{eq:bootstrap-coefficient-integrals} and \eqref{eq:forcing-fourth-time} give
\begin{equation*}
\|\partial_y^4W(\cdot,s)\|_{L^\infty}\leq\left(1+\frac{\vartheta}{8}\right)\left(\|\partial_y^4W(\cdot,s_0)\|_{L^\infty}+\frac{\vartheta M}{6}\right)\leq\|\partial_y^4W(\cdot,s_0)\|_{L^\infty}+\vartheta M,
\end{equation*}
where the last inequality follows from \eqref{eq:boot-fourth} at \(s=s_0\) and \(0<\vartheta<1/2\). This completes the proof.
\end{proof}

We next prove the estimates \eqref{eq:improved-local}, \eqref{eq:improved-second}, and \eqref{eq:improved-third}, using Propositions~\ref{prop:damping-inequalities} and~\ref{prop:local-integral} from Section~\ref{sec:profile-inequalities}.

\begin{lemma}\label{lem:bootstrap-local}
Under the assumptions of Proposition~\ref{prop:common-bootstrap}, the estimate \eqref{eq:improved-local} holds; that is,
\begin{equation*}
|W_y(y,s)-U'(y)|\leq\frac{y^2}{1+y^2}\left(\left\|\frac{1+z^2}{z^2}(W_y-U')(z,s_0)\right\|_{L_z^\infty}+\vartheta N\right)
\end{equation*}
for all \((y,s)\in\R\times[s_0,s_1]\).
\end{lemma}

\begin{proof}
Let \(\widetilde W:=W-U\), and define
\begin{equation*}
Z(y,s):=\frac{1+y^2}{y^2}\widetilde W_y(y,s),
\end{equation*}
with the quotient at the origin defined by continuity. Then, by \eqref{eq:f-full-transport},
\begin{equation}\label{eq:weighted-full-transport}
\partial_sZ+D^ZZ+\mathcal V\partial_yZ = F^Z-U''(y)\frac{1+y^2}{y^2}\int_0^y\frac{y'^2}{1+y'^2}Z(y',s)\,dy',
\end{equation}
where \(\mathcal V\) is given in \eqref{eq:remainder-transport-split}, and
\begin{subequations}
\begin{align}
&D^Z(y,s):=1+\frac{\widetilde W_y(y,s)+2U'(y)}\ell+\frac{2}{1+y^2}\left(\frac{2\ell+1}{2}+\frac{W(y,s)}y\right),\label{eq:local-damping-general}\\
& F^Z(y,s):=\frac{1+y^2}{y^2}\bigl(\mathcal R(y,s)-\mathcal A(y,s)U''(y)\bigr)-\frac{2\mathcal A(y,s)}{y(1+y^2)}Z(y,s).\label{eq:local-forcing-general}
\end{align}
\end{subequations}

First, we estimate \(Z\) near \(y=0\). By \eqref{eq:modulation-constraints-intro}, \(\widetilde W_y(0,s)=\widetilde W_{yy}(0,s)=0\). Hence Taylor's theorem gives
\begin{equation*}
\widetilde W_y(y,s)=\frac{y^2}{2}\bigl(\partial_y^3W(0,s)-U^{(3)}(0)\bigr)+\frac{y^3}{6}\bigl(\partial_y^4W(\bar{y},s)-U^{(4)}(\bar{y})\bigr)
\end{equation*}
for some \(\bar y\) between \(0\) and \(y\). Multiplying by \((1+y^2)/y^2\), and using \eqref{eq:boot-fourth} and \eqref{eq:origin-increment-bound}, we obtain
\begin{equation}\label{eq:local-core-control}
|Z(y,s)|\leq(1+y^2)\left(\frac12|\partial_y^3W(0,s_0)-U^{(3)}(0)|+\frac{\vartheta N}{16}+\frac{M+\|U^{(4)}\|_{L^\infty}}6|y|\right).
\end{equation}
Also, the continuous extension of \(Z\) at the origin satisfies
\begin{equation}\label{eq:origin-initial-weight}
|Z(0,s_0)|=\frac12|\partial_y^3W(0,s_0)-U^{(3)}(0)|\leq\|Z(\cdot,s_0)\|_{L^\infty}\leq N,
\end{equation}
where the last inequality follows from \eqref{eq:boot-local} at \(s=s_0\). With \(h\) defined by \eqref{eq:bootstrap-core-radii}, we have \(h\leq\vartheta/8\). Hence \eqref{eq:local-core-control} and \eqref{eq:origin-initial-weight} imply
\begin{equation}\label{eq:core-local}
|Z(y,s)|\leq(1+h^2)\left(\|Z(\cdot,s_0)\|_{L^\infty}+\frac{\vartheta N}{12}\right)\leq\|Z(\cdot,s_0)\|_{L^\infty}+\frac{\vartheta N}{8},\qquad |y|\leq h.
\end{equation}
By the definition of \(Z\), this proves \eqref{eq:improved-local} for \(|y|\leq h\).

We next consider the region \(|y|\geq h\). To apply Lemma~\ref{lem:bootstrap-comparison}, we derive a lower bound for the damping coefficient \(D^Z\), estimate the nonlocal integral and forcing terms in \eqref{eq:weighted-full-transport}, and verify the required bound at spatial infinity.

For the damping coefficient, using \(W=U+\widetilde W\) in \eqref{eq:local-damping-general}, we write
\begin{equation*}
D^Z(y,s)=D_{1,\ell}(y)+\frac{\widetilde W_y(y,s)}\ell+\frac{2\widetilde W(y,s)}{y(1+y^2)}.
\end{equation*}
Since \(\widetilde W(0,s)=0\), the bootstrap assumption \eqref{eq:boot-local} gives
\begin{equation}\label{eq:bootstrap-primitive-bound}
|\widetilde W(y,s)|\leq\int_0^{|y|}\frac{Ny'^2}{1+y'^2}\,dy'\leq\frac N3|y|^3.
\end{equation}
Thus, by \eqref{eq:boot-local} and \eqref{eq:N-local},
\begin{equation*}
\begin{split}
\left|\frac{\widetilde W_y(y,s)}\ell+\frac{2\widetilde W(y,s)}{y(1+y^2)}\right|&\leq\frac{Ny^2}{\ell(1+y^2)}+\frac{2N}{|y|(1+y^2)}\int_0^{|y|}\frac{y'^2}{1+y'^2}\,dy'\\ &\leq\left(\frac1\ell+\frac23\right)\frac{Ny^2}{1+y^2}\leq\frac{\eta_{\ell}y^2}{4(1+y^2)}.
\end{split}
\end{equation*}
Consequently,
\begin{equation}\label{eq:local-damping-lower}
D^Z(y,s)\geq D_{1,\ell}(y)-\frac{\eta_{\ell}y^2}{4(1+y^2)}.
\end{equation}

We now compare the integral term in \eqref{eq:weighted-full-transport} with the damping. By \eqref{eq:local-kernel},
\begin{equation}\label{eq:local-integral-bound}
\left|U''(y)\frac{1+y^2}{y^2}\int_0^y\frac{y'^2}{1+y'^2}Z(y',s)\,dy'\right|\leq K_{\ell}^{\mathrm{loc}}(y)\|Z(\cdot,s)\|_{L^\infty}.
\end{equation}
Moreover, \eqref{eq:local-integral-inequality} gives
\begin{equation}\label{eq:local-kernel-bound}
\lambda_{\ell}^{\mathrm{loc}}K_{\ell}^{\mathrm{loc}}(y)\leq D_{1,\ell}(y)-\frac{\eta_{\ell}y^2}{1+y^2}.
\end{equation}

For the forcing term, using the definition of \(F^Z\) in \eqref{eq:local-forcing-general} and the identity \(\widetilde W_y=y^2(1+y^2)^{-1}Z\), we obtain from \eqref{eq:common-local-comparison} that
\begin{equation}\label{eq:local-forcing-bound}
|F^Z(y,s)|\leq\left(\frac{\eta_{\ell}y^2}{4(1+y^2)}+a_0(s)\right)|Z(y,s)|+(\lambda_{\ell}^{\mathrm{loc}}-1)K_{\ell}^{\mathrm{loc}}(y)\|Z(\cdot,s)\|_{L^\infty}+e_0(s)
\end{equation}
for $|y|\geq h$. Combining \eqref{eq:local-damping-lower}, \eqref{eq:local-integral-bound}, and \eqref{eq:local-forcing-bound} with \eqref{eq:weighted-full-transport}, we have
\begin{equation}\label{eq:local-comparison-inequality}
\begin{split}
&(\partial_s+\mathcal V\partial_y)|Z|+\left(D_{1,\ell}(y)-\frac{\eta_{\ell}y^2}{2(1+y^2)}-a_0(s)\right)|Z|\\ &\quad\leq\lambda_{\ell}^{\mathrm{loc}}K_{\ell}^{\mathrm{loc}}(y)\|Z(\cdot,s)\|_{L^\infty}+e_0(s),\qquad |y|\geq h.
\end{split}
\end{equation}
By \eqref{eq:local-kernel-bound} and the definition of \(\mu_0\) in \eqref{eq:bootstrap-damping-rates},
\begin{equation}\label{eq:local-comparison-gap}
D_{1,\ell}(y)-\frac{\eta_{\ell}y^2}{2(1+y^2)}-\lambda_{\ell}^{\mathrm{loc}}K_{\ell}^{\mathrm{loc}}(y)\geq\frac{\eta_{\ell}y^2}{2(1+y^2)}\geq\mu_0,\qquad |y|\geq h.
\end{equation}
Finally, the assumption \eqref{eq:boot-far} gives
\begin{equation}\label{eq:local-tail-control}
|Z(y,s)|\leq\frac{B_{\ell}(1+y^2)}{y^2(1+|y|^{2/(2\ell+1)})}\longrightarrow0\quad\text{as }|y|\to\infty,
\end{equation}
uniformly for \(s\in[s_0,s_1]\).

We apply Lemma~\ref{lem:bootstrap-comparison} to \eqref{eq:local-comparison-inequality} with
\begin{equation*}
f:=Z,\qquad \Omega:=[-h,h],\qquad K(y,s):=\lambda_\ell^{\mathrm{loc}}K_\ell^{\mathrm{loc}}(y),
\end{equation*}
and
\begin{equation*}
D(y,s):=D_{1,\ell}(y)-\frac{\eta_\ell y^2}{2(1+y^2)}-a_0(s),\qquad H(y,s):=e_0(s).
\end{equation*}
Let \(\Phi\) be given by \eqref{eq:comparison-function} with
\begin{equation*}
L:=\|Z(\cdot,s_0)\|_{L^\infty}+\frac{\vartheta N}{8},\qquad \mu:=\mu_0,\qquad a:=a_0,\qquad e:=e_0.
\end{equation*}
Since \(\Phi(s_0)=L\), the initial condition \eqref{eq:maximum-principle-initial} holds. Since \(\Phi(s)\geq L\), \eqref{eq:core-local} gives \eqref{eq:maximum-principle-core}. Moreover, \eqref{eq:local-tail-control} gives \eqref{eq:maximum-principle-tail}. By \eqref{eq:local-comparison-gap},
\begin{equation*}
D(y,s)-K(y,s)\geq\mu_0-a_0(s),\qquad |y|\geq h.
\end{equation*}
Together with \eqref{eq:comparison-function-derivative}, this gives
\begin{equation*}
\Phi'(s)+(D(y,s)-K(y,s))\Phi(s)\geq e_0(s)=H(y,s),\qquad (y,s)\in\Omega^c\times[s_0,s_1],
\end{equation*}
which is \eqref{eq:comparison-barrier}. Lemma~\ref{lem:bootstrap-comparison} therefore gives
\begin{equation*}
\|Z(\cdot,s)\|_{L^\infty}\leq\Phi(s).
\end{equation*}
Using \eqref{eq:bootstrap-coefficient-integrals} and \eqref{eq:forcing-local-time}, we obtain
\begin{equation*}
\|Z(\cdot,s)\|_{L^\infty}\leq\Phi(s)\leq\left(1+\frac{\vartheta}{8}\right)\left(\|Z(\cdot,s_0)\|_{L^\infty}+\frac{\vartheta N}{4}\right)\leq\|Z(\cdot,s_0)\|_{L^\infty}+\vartheta N,
\end{equation*}
where we used \eqref{eq:boot-local} and \(0<\vartheta<1/2\) in the last inequality. This proves \eqref{eq:improved-local}.
\end{proof}

\begin{lemma}\label{lem:bootstrap-second}
Under the assumptions of Proposition~\ref{prop:common-bootstrap}, the estimate \eqref{eq:improved-second} holds; that is,
\begin{equation*}
|W_{yy}(y,s)|\leq\frac{|y|}{\sqrt{1+y^2}}\left(\left\|\frac{\sqrt{1+z^2}}{z}W_{yy}(z,s_0)\right\|_{L_z^\infty}+\vartheta M^{1/8}\right)
\end{equation*}
for all \((y,s)\in\R\times[s_0,s_1]\).
\end{lemma}

\begin{proof}
Let \(\widetilde W:=W-U\), and define
\begin{equation*}
\widetilde Z(y,s):=\frac{\sqrt{1+y^2}}{y}W_{yy}(y,s),
\end{equation*}
with the quotient at the origin defined by continuity. Then, by \eqref{eq:second-derivative-common},
\begin{equation}\label{eq:weighted-second-equation}
\partial_s \widetilde Z + D^{\widetilde Z} \widetilde Z + \mathcal V \partial_y \widetilde Z = F^{\widetilde Z},
\end{equation}
where \(\mathcal V\) is given in \eqref{eq:remainder-transport-split}, and
\begin{equation*}
\begin{aligned}
D^{\widetilde Z}(y,s) &:= \frac{2\ell+3}{2} + \left(1+\frac2\ell\right) W_y(y,s) + \frac{1}{1+y^2} \left( \frac{2\ell+1}{2} + \frac{W(y,s)}y \right),\\
F^{\widetilde Z}(y,s) &:= \frac{\sqrt{1+y^2}}{y}\mathcal R_2(y,s) - \frac{\mathcal A(y,s)}{y(1+y^2)}\widetilde Z(y,s).
\end{aligned}
\end{equation*}

First, we estimate \(\widetilde Z\) in the region \(|y|\leq M^{-1}\). By Taylor's theorem, \eqref{eq:modulation-constraints-intro}, \eqref{eq:boot-origin}, and \eqref{eq:boot-fourth},
\begin{equation*}
|W_{yy}(y,s)|\leq|y||\partial_y^3W(0,s)|+\frac{y^2}{2}\|\partial_y^4W(\cdot,s)\|_{L^\infty}\leq(U^{(3)}(0)+1)|y|+\frac M2y^2.
\end{equation*}
Thus, for \(|y|\leq M^{-1}\),
\begin{equation}\label{eq:core-second}
|\widetilde Z(y,s)|\leq\sqrt{1+y^2}\left(U^{(3)}(0)+1+\frac M2|y|\right)\leq2(U^{(3)}(0)+2)\leq\frac{\vartheta M^{1/8}}8,
\end{equation}
provided \(M\) is sufficiently large that \(\vartheta M^{1/8}\geq16(U^{(3)}(0)+2)\). By the definition of \(\widetilde Z\), this proves \eqref{eq:improved-second} for \(|y|\leq M^{-1}\).

We next turn to the region \(|y|\geq M^{-1}\). To apply Lemma~\ref{lem:bootstrap-comparison}, we obtain a lower bound for \(D^{\widetilde Z}\), an upper bound for \(F^{\widetilde Z}\), and the required bound at spatial infinity.

For the damping coefficient, using \(W=U+\widetilde W\) and \eqref{eq:D2}, we write
\begin{equation*}
D^{\widetilde Z}(y,s)=D_{2,\ell}(y)+\left(1+\frac2\ell\right)\widetilde W_y(y,s)+\frac{\widetilde W(y,s)}{y(1+y^2)}.
\end{equation*}
Using \eqref{eq:boot-local}, \eqref{eq:bootstrap-primitive-bound}, and \eqref{eq:N-local}, we have
\begin{equation*}
\begin{split}
\left|\left(1+\frac2\ell\right)\widetilde W_y(y,s)+\frac{\widetilde W(y,s)}{y(1+y^2)}\right|&\leq\left(1+\frac2\ell\right)\frac{Ny^2}{1+y^2}+\frac{N}{|y|(1+y^2)}\int_0^{|y|}\frac{y'^2}{1+y'^2}\,dy'\\ &\leq\left(\frac43+\frac2\ell\right)\frac{Ny^2}{1+y^2}\leq\frac{c_{2,\ell}y^2}{4(1+y^2)}.
\end{split}
\end{equation*}
Together with \eqref{eq:bootstrap-profile-damping}, this gives
\begin{equation}\label{eq:second-damping-lower}
D^{\widetilde Z}(y,s)\geq\frac{3c_{2,\ell}y^2}{4(1+y^2)}.
\end{equation}

For the forcing term, using the definition of \(F^{\widetilde Z}\), the assumption \eqref{eq:bootstrap-second-forcing} gives
\begin{equation}\label{eq:second-forcing-bound}
|F^{\widetilde Z}(y,s)|\leq\left(\frac{c_{2,\ell}y^2}{4(1+y^2)}+a_2(s)\right)|\widetilde Z(y,s)|+e_2(s),\qquad |y|\geq M^{-1}.
\end{equation}
Combining \eqref{eq:weighted-second-equation}, \eqref{eq:second-damping-lower}, and \eqref{eq:second-forcing-bound}, we obtain
\begin{equation}\label{eq:second-comparison-inequality}
(\partial_s+\mathcal V\partial_y)|\widetilde Z|+\left(\frac{c_{2,\ell}y^2}{2(1+y^2)}-a_2(s)\right)|\widetilde Z|\leq e_2(s),\qquad |y|\geq M^{-1}.
\end{equation}
For \(|y|\geq M^{-1}\), we also have
\begin{equation}\label{eq:second-comparison-damping}
\frac{c_{2,\ell}y^2}{2(1+y^2)}\geq\frac{c_{2,\ell}}{2(1+M^2)}=\mu_2.
\end{equation}

Finally, we estimate \(\widetilde Z\) at spatial infinity. Applying Lemma~\ref{lem:transport-far-field} to \eqref{eq:second-comparison-inequality}, with \(R_0=M^{-1}\), \(\lambda_D=\mu_2\), \(a=a_2\), and \(e=e_2\), we obtain
\begin{equation}\label{eq:second-tail-control}
\limsup_{|y|\to\infty}|\widetilde Z(y,s)|\leq e^{\int_{s_0}^sa_2(s')\,ds'}\left(\|\widetilde Z(\cdot,s_0)\|_{L^\infty}+\int_{s_0}^{s}e^{-\mu_2(s-s')}e_2(s')\,ds'\right).
\end{equation}

We now apply Lemma~\ref{lem:bootstrap-comparison} to \eqref{eq:second-comparison-inequality} with
\begin{equation*}
f:=\widetilde Z,\qquad \Omega:=[-M^{-1},M^{-1}],\qquad K:=0,
\end{equation*}
and
\begin{equation*}
D(y,s):=\frac{c_{2,\ell}y^2}{2(1+y^2)}-a_2(s),\qquad H(y,s):=e_2(s).
\end{equation*}
Let \(\Phi\) be given by \eqref{eq:comparison-function} with
\begin{equation*}
L:=\|\widetilde Z(\cdot,s_0)\|_{L^\infty}+\frac{\vartheta M^{1/8}}8,\qquad \mu:=\mu_2,\qquad a:=a_2,\qquad e:=e_2.
\end{equation*}
Since \(\Phi(s_0)=L\), \eqref{eq:maximum-principle-initial} holds. Since \(\Phi(s)\geq L\geq\vartheta M^{1/8}/8\), \eqref{eq:core-second} gives \eqref{eq:maximum-principle-core}. The estimate \eqref{eq:second-tail-control} gives \eqref{eq:maximum-principle-tail}, and \eqref{eq:second-comparison-damping} together with \eqref{eq:comparison-function-derivative} gives \eqref{eq:comparison-barrier}. Hence, by Lemma~\ref{lem:bootstrap-comparison}, \eqref{eq:bootstrap-coefficient-integrals} and \eqref{eq:forcing-second-time}, we conclude that
\begin{equation*}
\|\widetilde Z(\cdot,s)\|_{L^\infty}\leq\Phi(s)\leq\left(1+\frac{\vartheta}{8}\right)\left(\|\widetilde Z(\cdot,s_0)\|_{L^\infty}+\frac{\vartheta M^{1/8}}4\right)\leq\|\widetilde Z(\cdot,s_0)\|_{L^\infty}+\vartheta M^{1/8},
\end{equation*}
where the last inequality follows from \eqref{eq:boot-second} at \(s=s_0\) and \(0<\vartheta<1/2\). This proves \eqref{eq:improved-second}.
\end{proof}

\begin{lemma}\label{lem:bootstrap-third}
Under the assumptions of Proposition~\ref{prop:common-bootstrap}, the estimate \eqref{eq:improved-third} holds; that is,
\begin{equation*}
\|\partial_y^3W(\cdot,s)\|_{L^\infty}\leq\|\partial_y^3W(\cdot,s_0)\|_{L^\infty}+\vartheta M^{3/4}
\end{equation*}
for all \(s\in[s_0,s_1]\).
\end{lemma}

\begin{proof}
We rewrite \eqref{eq:third-derivative-common} as
\begin{equation}\label{eq:third-transport-proof}
\partial_s\partial_y^3W+D_3^W\partial_y^3W+\mathcal V\partial_y^4W=F_3^W,
\end{equation}
where
\begin{equation*}
\begin{aligned}
&D_3^W(y,s) := 2\ell+2+\left(2+\frac2\ell\right)W_y(y,s),\\
&F_3^W(y,s) := -\left(1+\frac2\ell\right)W_{yy}(y,s)^2+\mathcal R_3(y,s).
\end{aligned}
\end{equation*}
We first consider the region \(|y|\leq\vartheta/(16M^{1/4})\). By \eqref{eq:boot-origin} and \eqref{eq:boot-fourth},
\begin{equation}\label{eq:core-third}
\begin{split}
|\partial_y^3W(y,s)|&\leq|\partial_y^3W(0,s)|+|y|\|\partial_y^4W(\cdot,s)\|_{L^\infty}\\
&\leq U^{(3)}(0)+1+M|y|\leq\frac{\vartheta M^{3/4}}8,\qquad |y|\leq\frac{\vartheta}{16M^{1/4}}
\end{split}
\end{equation}
for sufficiently large \(M>0\).

We next consider the region \(|y|\geq\vartheta/(16M^{1/4})\). By \eqref{eq:bootstrap-slope-lower}, we have
\begin{equation}\label{eq:third-damping-lower}
D_3^W(y,s)=\left(2+\frac2\ell\right)(\ell+W_y(y,s))\geq\left(1+\frac1\ell\right)\frac{c_{1,\ell}y^2}{1+y^2}.
\end{equation}
By \eqref{eq:bootstrap-third-forcing} and \eqref{eq:boot-second}, the term \(F_3^W\) satisfies
\begin{equation}\label{eq:third-forcing-bound}
\begin{split}
& |F_3^W(y,s)| \\
& \quad \leq \left(\frac{(1+1/\ell)c_{1,\ell}y^2}{4(1+y^2)}+a_3(s)\right)|\partial_y^3W(y,s)| + \left(1+\frac2\ell+\Gamma_3\right)|W_{yy}(y,s)|^2+e_3(s)\\
& \quad \leq \left(\frac{(1+1/\ell)c_{1,\ell}y^2}{4(1+y^2)}+a_3(s)\right)|\partial_y^3W(y,s)|   +\left(1+\frac2\ell+\Gamma_3\right)\frac{M^{1/4}y^2}{1+y^2}+e_3(s).
\end{split}
\end{equation}
Combining \eqref{eq:third-damping-lower} and \eqref{eq:third-forcing-bound} with \eqref{eq:third-transport-proof}, we obtain
\begin{equation}\label{eq:third-comparison-inequality}
\begin{split}
&(\partial_s+\mathcal V\partial_y)|\partial_y^3W|+\left(\frac{3(1+1/\ell)c_{1,\ell}y^2}{4(1+y^2)}-a_3(s)\right)|\partial_y^3W|\\
&\quad\leq\left(1+\frac2\ell+\Gamma_3\right)\frac{M^{1/4}y^2}{1+y^2}+e_3(s),\qquad |y|\geq\frac{\vartheta}{16M^{1/4}}.
\end{split}
\end{equation}
We estimate \(\partial_y^3W\) at spatial infinity. Since $\vartheta/(16M^{1/4})<1$ for sufficiently large $M>0$, the inequality \eqref{eq:third-comparison-inequality} holds on \(|y|\geq1\). On this region,
\begin{equation*}
\frac{3(1+1/\ell)c_{1,\ell}y^2}{4(1+y^2)} \geq\frac{3(1+1/\ell)c_{1,\ell}}8 =:\lambda_D
\end{equation*}
and
\begin{equation*}
\left(1+\frac2\ell+\Gamma_3\right)\frac{M^{1/4}y^2}{1+y^2}+e_3(s) \leq\left(1+\frac2\ell+\Gamma_3\right)M^{1/4}+e_3(s).
\end{equation*}
Applying Lemma~\ref{lem:transport-far-field} to \eqref{eq:third-comparison-inequality} with \(R_0=1\), \(\lambda_D\) as above, \(a=a_3\), and
\begin{equation*}
e(s)=\left(1+\frac2\ell+\Gamma_3\right)M^{1/4}+e_3(s),
\end{equation*}
we obtain
\begin{equation}\label{eq:third-tail-integral}
\begin{split}
& \limsup_{|y|\to\infty}|\partial_y^3W(y,s)|\\
& \quad \leq e^{\int_{s_0}^sa_3(s') \, ds'}\bigg(\|\partial_y^3W(\cdot,s_0)\|_{L^\infty}e^{-\lambda_D(s-s_0)} \\
&\qquad \qquad \qquad \qquad +\int_{s_0}^{s}e^{-\lambda_D(s-s')}\left[\left(1+\frac2\ell+\Gamma_3\right)M^{1/4}+e_3(s')\right] \, ds'\bigg).
\end{split}
\end{equation}
We choose \(M\) sufficiently large so that
\begin{equation} \label{eq:third-M-choice}
32\left(1+\frac2\ell+\Gamma_3\right)M^{1/4} \leq\left(1+\frac1\ell\right)c_{1,\ell}\vartheta M^{3/4}.
\end{equation}
Then it holds that
\begin{equation*}
\frac{1+2/\ell+\Gamma_3}{\lambda_D}M^{1/4} \leq\frac{\vartheta M^{3/4}}{12}.
\end{equation*}
Moreover, by \eqref{eq:bootstrap-damping-rates} and the definition of \(\lambda_D\),
\begin{equation*}
\frac{\mu_3}{\lambda_D}=\frac{4\vartheta^2}{3(\vartheta^2+256M^{1/2})}<1,
\end{equation*}
and hence, \(e^{-\lambda_D(s-s')}\leq e^{-\mu_3(s-s')}\). Using these estimates in \eqref{eq:third-tail-integral}, we obtain
\begin{equation}\label{eq:third-tail-control}
\begin{split}
& \limsup_{|y|\to\infty}|\partial_y^3W(y,s)| \\
& \quad \leq e^{\int_{s_0}^sa_3(s')\,ds'}\left(\|\partial_y^3W(\cdot,s_0)\|_{L^\infty}+\frac{\vartheta M^{3/4}}{12}+\int_{s_0}^{s}e^{-\mu_3(s-s')}e_3(s')\,ds'\right).
\end{split}
\end{equation}

We now apply Lemma~\ref{lem:bootstrap-comparison} to \eqref{eq:third-comparison-inequality} with
\begin{equation*}
f:=\partial_y^3W,\qquad \Omega:=\left[-\frac{\vartheta}{16M^{1/4}},\frac{\vartheta}{16M^{1/4}}\right],\qquad K:=0,
\end{equation*}
and
\begin{equation*}
D(y,s):=\frac{3(1+1/\ell)c_{1,\ell}y^2}{4(1+y^2)}-a_3(s),\qquad H(y,s):=\left(1+\frac2\ell+\Gamma_3\right)\frac{M^{1/4}y^2}{1+y^2}+e_3(s).
\end{equation*}
Let \(\Phi\) be given by \eqref{eq:comparison-function} with
\begin{equation*}
L:=\|\partial_y^3W(\cdot,s_0)\|_{L^\infty}+\frac{\vartheta M^{3/4}}8,\qquad \mu:=\mu_3,\qquad a:=a_3,\qquad e:=e_3.
\end{equation*}
Since \(\Phi(s_0)=L\), the initial condition \eqref{eq:maximum-principle-initial} holds. Moreover, \(\Phi(s)\geq L\geq\vartheta M^{3/4}/8\), so \eqref{eq:core-third} gives \eqref{eq:maximum-principle-core}. Since \(\vartheta M^{3/4}/12\leq\vartheta M^{3/4}/8\), \eqref{eq:third-tail-control} gives \eqref{eq:maximum-principle-tail}. It therefore remains to verify \eqref{eq:comparison-barrier}. By \eqref{eq:third-M-choice} and \(\Phi(s)\geq\vartheta M^{3/4}/8\),
\begin{equation*}
\left(1+\frac2\ell+\Gamma_3\right)\frac{M^{1/4}y^2}{1+y^2}\leq\frac{(1+1/\ell)c_{1,\ell}y^2}{4(1+y^2)}\Phi(s).
\end{equation*}
Moreover, by \eqref{eq:bootstrap-damping-rates},
\begin{equation*}
\frac{(1+1/\ell)c_{1,\ell}y^2}{2(1+y^2)}\geq\mu_3,\qquad |y|\geq\frac{\vartheta}{16M^{1/4}}.
\end{equation*}
Together with \eqref{eq:comparison-function-derivative}, these estimates yield
\begin{equation*}
\Phi'(s)+D(y,s)\Phi(s)\geq H(y,s),\qquad (y,s)\in\Omega^c\times[s_0,s_1],
\end{equation*}
which is \eqref{eq:comparison-barrier}. Lemma~\ref{lem:bootstrap-comparison} therefore gives
\begin{equation*}
\|\partial_y^3W(\cdot,s)\|_{L^\infty}\leq\Phi(s).
\end{equation*}
Using \eqref{eq:bootstrap-coefficient-integrals}, \eqref{eq:forcing-third-time}, and \eqref{eq:boot-third} at \(s=s_0\), we conclude that
\begin{equation*}
\begin{split}
\|\partial_y^3W(\cdot,s)\|_{L^\infty} \leq \Phi(s) &\leq \left(1+\frac{\vartheta}{8}\right)\left(\|\partial_y^3W(\cdot,s_0)\|_{L^\infty}+\frac{\vartheta M^{3/4}}4\right) \\
& \leq\|\partial_y^3W(\cdot,s_0)\|_{L^\infty}+\vartheta M^{3/4}.
\end{split}
\end{equation*}
This completes the proof.
\end{proof}

The next lemma proves \eqref{eq:improved-far}.

\begin{lemma}\label{lem:bootstrap-far}
Under the assumptions of Proposition~\ref{prop:common-bootstrap}, the estimate \eqref{eq:improved-far} holds; that is,
\begin{equation*}
\begin{aligned}
&\bigl(1+|y|^{2/(2\ell+1)}\bigr)|W_y(y,s)-U'(y)|\\
&\quad \leq \max \left\{\left\|\bigl(1+|z|^{2/(2\ell+1)}\bigr)(W_y-U')(z,s_0)\right\|_{L_z^\infty},N\bigl(1+m^{2/(2\ell+1)}\bigr)\right\}+\vartheta B_\ell
\end{aligned}
\end{equation*}
for all \((y,s)\in\R\times[s_0,s_1]\).
\end{lemma}

\begin{proof}
Set
\begin{equation*}
r:=\frac{2}{2\ell+1},\qquad \widetilde W:=W-U,\qquad \nu(y,s):=(1+|y|^r)\widetilde W_y(y,s),
\end{equation*}
and let
\begin{equation*}
L:=\max\left\{\|\nu(\cdot,s_0)\|_{L^\infty},N(1+m^r)\right\}.
\end{equation*}
Recall from Section~\ref{sec:bootstrap-assumptions} that \(m\geq\max\{1,m_0(\ell)\}\), where \(m_0(\ell)\) is the threshold in Proposition~\ref{prop:far-field}.

We note that, for \(|y|\leq m\), the assumption \eqref{eq:boot-local} directly gives
\begin{equation}\label{eq:far-interior-control}
|\nu(y,s)|\leq(1+|y|^r)\frac{Ny^2}{1+y^2}\leq N(1+m^r)\leq L.
\end{equation}
Thus \eqref{eq:improved-far} holds for \(|y|\leq m\), and it remains to prove
\begin{equation*}
|\nu(y,s)|\leq L+\vartheta B_\ell,\qquad |y|\geq m.
\end{equation*}

From \eqref{eq:f-full-transport}, we have
\begin{equation}\label{eq:far-weighted-equation}
\partial_s\nu+D^\nu\nu+\mathcal V\partial_y\nu=F^\nu-(1+|y|^r)U''(y)\int_0^y\frac{\nu(y',s)}{1+|y'|^r}\,dy',
\end{equation}
where \(\mathcal V\) is defined in \eqref{eq:remainder-transport-split}, and
\begin{subequations}
\begin{align}
D^\nu(y,s)&:=1+\frac{\widetilde W_y(y,s)+2U'(y)}{\ell}-\frac{r|y|^r}{1+|y|^r}\left(\frac{2\ell+1}{2}+\frac{U(y)+\widetilde W(y,s)}{y}\right),\label{eq:far-actual-damping}\\
F^\nu(y,s)&:=(1+|y|^r)\bigl(\mathcal R(y,s)-\mathcal A(y,s)U''(y)\bigr)+\frac{r|y|^r\mathcal A(y,s)}{y(1+|y|^r)}\nu(y,s).\label{eq:far-forcing-general}
\end{align}
\end{subequations}
We first estimate \(D^\nu\), \(F^\nu\), and the integral term in \eqref{eq:far-weighted-equation}. By \eqref{eq:boot-far} and \(\widetilde W(0,s)=0\),
\begin{equation}\label{eq:far-perturbation-bounds}
\widetilde W_y(y,s)\geq-\frac{B_\ell}{1+|y|^r},\qquad \left|\frac{\widetilde W(y,s)}{y}\right|\leq\frac{B_\ell}{|y|}\int_0^{|y|}\frac{dy'}{1+|y'|^r}.
\end{equation}
Substituting \eqref{eq:far-perturbation-bounds} into \eqref{eq:far-actual-damping}, we obtain
\begin{equation}\label{eq:far-damping-common}
\begin{aligned}
D^\nu(y,s)&\geq1-\frac{B_\ell}{\ell(1+|y|^r)}+\frac{2}{\ell}U'(y)-\frac{r|y|^r}{1+|y|^r}\left(\frac{2\ell+1}{2}+\frac{U(y)}{y}+\frac{B_\ell}{y}\int_0^y\frac{dy'}{1+|y'|^r}\right)\\
&=D_{\ell,B_\ell}^{\mathrm{far}}(y)\geq0,\qquad |y|\geq m,
\end{aligned}
\end{equation}
where \(D_{\ell,B}^{\mathrm{far}}\) is defined in \eqref{eq:far-damping}, and the last inequality follows from Proposition~\ref{prop:far-field}. Using the definition of \(F^\nu\), the assumption \eqref{eq:common-far-comparison} gives
\begin{equation}\label{eq:far-forcing-bound}
|F^\nu(y,s)|\leq(\Lambda_\ell-1)K_\ell^{\mathrm{far}}(y)\|\nu(\cdot,s)\|_{L^\infty}+a_\infty(s)|\nu(y,s)|+e_\infty(s),\qquad |y|\geq m,
\end{equation}
where \(K_\ell^{\mathrm{far}}\) is defined in \eqref{eq:far-kernel}. Moreover, \eqref{eq:far-kernel} gives
\begin{equation}\label{eq:far-integral-bound}
\left|(1+|y|^r)U''(y)\int_0^y\frac{\nu(y',s)}{1+|y'|^r}\,dy'\right|\leq K_\ell^{\mathrm{far}}(y)\|\nu(\cdot,s)\|_{L^\infty}.
\end{equation}
The estimates \eqref{eq:far-forcing-bound} and \eqref{eq:far-integral-bound}, together with \eqref{eq:far-weighted-equation}, yield
\begin{equation}\label{eq:far-comparison-inequality}
(\partial_s+\mathcal V\partial_y)|\nu|+(D^\nu-a_\infty(s))|\nu|\leq\Lambda_\ell K_\ell^{\mathrm{far}}(y)\|\nu(\cdot,s)\|_{L^\infty}+e_\infty(s),\qquad |y|\geq m.
\end{equation}

Next, we estimate \(\nu\) at spatial infinity. Since \eqref{eq:far-Usecond} and \(0<r<1\) imply \(|U''(y)|\leq C|y|^{-1-r}\) and \(\int_0^{|y|}(1+z^r)^{-1}\,dz\leq C|y|^{1-r}\) for \(|y|\geq1\), we obtain
\begin{equation}\label{eq:far-kernel-decay}
K_\ell^{\mathrm{far}}(y)=(1+|y|^r)|U''(y)|\int_0^{|y|}\frac{dz}{1+z^r}\leq C|y|^{-r}\longrightarrow0\quad\text{as }|y|\to\infty.
\end{equation}
For any \(R\geq m\), \eqref{eq:far-damping-common}, \eqref{eq:far-comparison-inequality}, and \eqref{eq:boot-far} give
\begin{equation*}
(\partial_s+\mathcal V\partial_y)|\nu|-a_\infty(s)|\nu|\leq e_R(s),\qquad |y|\geq R,
\end{equation*}
where
\begin{equation*}
e_R(s):=e_\infty(s)+\Lambda_\ell B_\ell\sup_{|y'|\geq R}K_\ell^{\mathrm{far}}(y').
\end{equation*}
Applying Lemma~\ref{lem:transport-far-field} with \(f=\nu\), \(R_0=R\), \(\lambda_D=0\), \(a=a_\infty\), and \(e=e_R\), and using \(a_\infty\geq0\), we obtain
\begin{equation}\label{eq:far-tail-before-limit}
\begin{aligned}
\limsup_{|y|\to\infty}|\nu(y,s)|\leq e^{\int_{s_0}^{s}a_\infty(s')\,ds'}\bigg(&\limsup_{|y|\to\infty}|\nu(y,s_0)|+\int_{s_0}^{s}e_\infty(s')\,ds'\\
&+\Lambda_\ell B_\ell(s-s_0)\sup_{|y'|\geq R}K_\ell^{\mathrm{far}}(y')\bigg).
\end{aligned}
\end{equation}
Since \(s-s_0\leq s_1-s_0\), \eqref{eq:far-kernel-decay} shows that the last term in \eqref{eq:far-tail-before-limit} tends to zero uniformly for \(s\in[s_0,s_1]\) as \(R\to\infty\). Letting \(R\to\infty\), we obtain
\begin{equation}\label{eq:far-tail-control}
\limsup_{|y|\to\infty}|\nu(y,s)|\leq e^{\int_{s_0}^{s}a_\infty(s')\,ds'}\left(\limsup_{|y|\to\infty}|\nu(y,s_0)|+\int_{s_0}^{s}e_\infty(s')\,ds'\right).
\end{equation}

We now use \eqref{eq:far-tail-control} to prove \eqref{eq:improved-far} on the whole line. By \eqref{eq:far-field-inequality} and \eqref{eq:far-damping-common},
\begin{equation}\label{eq:far-damping-kernel-comparison}
\Lambda_\ell K_\ell^{\mathrm{far}}(y)\leq D_{\ell,B_\ell}^{\mathrm{far}}(y)\leq D^\nu(y,s),\qquad |y|\geq m.
\end{equation}
Define
\begin{equation}\label{eq:far-comparison-function}
\Phi(s):=e^{\int_{s_0}^{s}a_\infty(s')\,ds'}\left(L+\int_{s_0}^{s}e_\infty(s')\,ds'\right).
\end{equation}
We apply Lemma~\ref{lem:bootstrap-comparison} to \eqref{eq:far-comparison-inequality} with
\begin{equation*}
f:=\nu,\qquad \Omega:=[-m,m],\qquad K(y,s):=\Lambda_\ell K_\ell^{\mathrm{far}}(y),
\end{equation*}
and
\begin{equation*}
D(y,s):=D^\nu(y,s)-a_\infty(s),\qquad H(y,s):=e_\infty(s).
\end{equation*}
Since \(\Phi(s_0)=L\), the initial condition \eqref{eq:maximum-principle-initial} holds. Since \(\Phi(s)\geq L\), \eqref{eq:far-interior-control} gives \eqref{eq:maximum-principle-core}. Moreover, \eqref{eq:far-tail-control} and \(\limsup_{|y|\to\infty}|\nu(y,s_0)|\leq\|\nu(\cdot,s_0)\|_{L^\infty}\leq L\) give \eqref{eq:maximum-principle-tail}. Direct differentiation of \eqref{eq:far-comparison-function} gives
\begin{equation*}
\Phi'(s)=a_\infty(s)\Phi(s)+e^{\int_{s_0}^{s}a_\infty(s')\,ds'}e_\infty(s)\geq a_\infty(s)\Phi(s)+e_\infty(s).
\end{equation*}
Together with \eqref{eq:far-damping-kernel-comparison}, this yields
\begin{equation*}
\Phi'(s)+(D(y,s)-K(y,s))\Phi(s)\geq H(y,s),\qquad (y,s)\in\Omega^c\times[s_0,s_1],
\end{equation*}
which is \eqref{eq:comparison-barrier}. Lemma~\ref{lem:bootstrap-comparison} therefore gives
\begin{equation*}
\|\nu(\cdot,s)\|_{L^\infty}\leq\Phi(s),\qquad s\in[s_0,s_1].
\end{equation*}

By \eqref{eq:boot-far} and \eqref{eq:N-second-far}, \(L\leq B_\ell\). Using \eqref{eq:bootstrap-coefficient-integrals} and \eqref{eq:common-forcing-smallness} in \eqref{eq:far-comparison-function}, we obtain
\begin{equation*}
\|\nu(\cdot,s)\|_{L^\infty}\leq\left(1+\frac{\vartheta}{8}\right)\left(L+\frac{\vartheta B_\ell}{8}\right)\leq L+\vartheta B_\ell,
\end{equation*}
where the last inequality follows from \(L\leq B_\ell\) and \(0<\vartheta<1/2\). By the definitions of \(L\) and \(\nu\), this proves \eqref{eq:improved-far}.
\end{proof}

\subsection{Proof of Theorem~\ref{thm:common-singularity}}\label{sec:proof-common-singularity}

By the local existence and the continuity of the weighted quantities in \emph{(i)} and the initial bounds in \emph{(ii)}, the bootstrap assumptions \eqref{eq:common-bootstrap-assumptions} hold on some interval \([s_0,s_1]\) with \(s_1>s_0\). By Proposition~\ref{prop:common-bootstrap} and assumptions \emph{(ii)}--\emph{(iv)}, the improved bounds \eqref{eq:bootstrap-improved-bounds} and the modulation estimates \eqref{eq:common-singularity-modulation} hold on this interval.

Let \([s_0,S)\) be the maximal interval on which the self-similar formulation and the bootstrap assumptions hold. The same argument gives \eqref{eq:bootstrap-improved-bounds} and \eqref{eq:common-singularity-modulation} on every finite interval contained in \([s_0,S)\). If \(S<\infty\), assumption \emph{(v)} extends the solution and the modulation functions beyond \(S\). The strict margin in \eqref{eq:bootstrap-improved-bounds} and the continuity properties in \emph{(i)} then imply that the bootstrap assumptions also hold beyond \(S\), contradicting the maximality of \(S\). Hence \(S=\infty\). It follows that \(W\) is defined for all \(s\geq s_0\), the self-similar formulation holds, and \eqref{eq:bootstrap-improved-bounds} holds on \(\R\times[s_0,\infty)\). Moreover, the modulation estimates \eqref{eq:common-singularity-modulation} hold for \(t=t(s)\), \(s\geq s_0\), with constants independent of \(s\).

We now turn to the original variables. By \eqref{eq:common-singularity-modulation}, \(1-\dot\tau(t)\geq c>0\), so \(t=t(s)\) is continuously differentiable. Differentiating \(s(t)=-\log(\tau(t)-t)\), we have
\begin{equation*}
\frac{dt}{ds}=\frac{e^{-s}}{1-\dot\tau(t(s))}.
\end{equation*}
By \eqref{eq:common-singularity-modulation}, \(t(s)\) has a finite limit as \(s\to\infty\). We denote this limit by
\begin{equation*}
T_*:=t_0+\int_{s_0}^{\infty}\frac{e^{-s'}}{1-\dot\tau(t(s'))}\,ds'>t_0.
\end{equation*}
The modulation estimates \eqref{eq:common-singularity-modulation} also give
\begin{equation*}
\frac{e^{-s}}{C}\leq T_*-t(s)=\int_s^{\infty}\frac{e^{-s'}}{1-\dot\tau(t(s'))}\,ds'\leq\frac{e^{-s}}{c}.
\end{equation*}
Since \(\tau(t(s))-t(s)=e^{-s}\), we have
\begin{equation*}
\tau(t)-t\sim T_*-t\quad\text{as }t\uparrow T_*.
\end{equation*}
Note that the bound on \(\int_{t_0}^{t}|\dot\xi(t')|\,dt'\) in \eqref{eq:common-singularity-modulation} is uniform for \(t<T_*\). Hence \(\dot\xi\in L^1([t_0,T_*))\), and there exists \(x_*\in\R\) such that
\begin{equation*}
\xi(t)\longrightarrow x_*\quad\text{as }t\uparrow T_*.
\end{equation*}
This proves \eqref{eq:common-singularity-time-center}. The same modulation estimates give the boundedness of \(\kappa\). Since the self-similar formulation and the regularity properties in \emph{(i)} hold on every finite self-similar time interval, \(w\) extends as a solution to \(\R\times[t_0,T_*)\).

Finally, \eqref{eq:boot-far}, \eqref{eq:boot-second}, and \eqref{eq:far-Uprime} imply
\begin{equation*}
|W_y(y,s)|\leq\frac{C}{1+|y|^{2/(2\ell+1)}},\qquad \sup_{|y|\leq1}|W_{yy}(y,s)|\leq M^{1/8}
\end{equation*}
for all \(y\in\R\) and \(s\geq s_0\). Thus \(W\) satisfies \eqref{eq:regularity-profile-bounds}. Applying Proposition~\ref{prop:sharp-regularity} with \eqref{eq:common-singularity-time-center} and the boundedness of \(\kappa\), we obtain the remaining conclusions. This completes the proof.

\section{Proofs of the profile inequalities}\label{sec:coercivity}

In this section, we prove the inequalities stated in Section~\ref{sec:profile-inequalities}. Throughout this section, we fix a real parameter $\ell>1/2$ and take $U=U_\beta$, with $\beta$ given by \eqref{eq:normalized-profile-choice}.

\subsection{Derivative bounds and localized inequalities}\label{sec:localized-inequality-proofs}

We first prove Lemma~\ref{lem:derivative-bounds}.

\begin{proof}[Proof of Lemma~\ref{lem:derivative-bounds}]
By \eqref{eq:Uprime-rho}, we have \(-\ell\leq U'(y)<0\) for all \(y\in\R\), and \(U'(0)=-\ell\). This proves \eqref{eq:Uprime-norm}.

To prove \eqref{eq:Usecond-norm}, we use the expression for \(U''\) in \eqref{eq:Uprime-rho}. A direct differentiation gives
\begin{equation*}
\frac{d}{d\rho}\left(\frac{\rho}{(1+\rho^2)^{\ell+2}}\right)=\frac{1-(2\ell+3)\rho^2}{(1+\rho^2)^{\ell+3}}.
\end{equation*}
Thus \(U''\) attains its maximum on \(y\geq0\) at \(\rho=(2\ell+3)^{-1/2}\). By the odd symmetry of \(U''\), we obtain
\begin{equation*}
\|U''\|_{L^\infty}= \ell(2\ell)^{\ell}\sqrt{\beta}\,\frac{(2\ell+3)^{\ell+3/2}}{(2\ell+4)^{\ell+2}},
\end{equation*}
which is \eqref{eq:Usecond-norm}. For the third derivative, \eqref{eq:Uthird-rho} gives
\begin{equation*}
 |U^{(3)}(y)|=\frac{\ell}{2}(2\ell)^{2\ell}\beta\frac{|1-(2\ell+3)\rho^2|}{(1+\rho^2)^{2\ell+3}}\leq \frac{\ell}{2}(2\ell)^{2\ell}\beta\frac{1+(2\ell+3)\rho^2}{(1+\rho^2)^{2\ell+3}}\leq \frac{\ell}{2}(2\ell)^{2\ell}\beta,
\end{equation*}
where the last inequality follows from \((1+\rho^2)^{2\ell+3}\geq1+(2\ell+3)\rho^2\). Equality holds at \(\rho=0\), proving \eqref{eq:Uthird-norm}.

Differentiating \eqref{eq:y-rho} and \eqref{eq:Uthird-rho}, we also obtain
\begin{equation*}
U^{(4)}(y)=\frac{\ell}{2}(2\ell)^{3\ell}\beta^{3/2}(2\ell+3)\frac{\rho\bigl((2\ell+2)\rho^2-2\bigr)}{(1+\rho^2)^{3\ell+4}}.
\end{equation*}
The right-hand side is bounded uniformly for $\rho\geq0$. By the odd symmetry of $U^{(4)}$, this proves \eqref{eq:Ufourth-bound} and completes the proof.
\end{proof}

We next turn to Proposition~\ref{prop:damping-inequalities}.

\begin{proof}[Proof of Proposition~\ref{prop:damping-inequalities}]
For \(\ell=2\), the estimates in \eqref{eq:CH-damping} follow from \cite[Lemmas~6.1--6.2]{KKY-CH}. Let \(\ell>1/2\), with \(\ell\ne2\). By the odd symmetry of \(U\), all the quantities in the statement are even in \(y\). Thus it is enough to consider \(y>0\); the bounds at \(y=0\) follow by continuity. For $y>0$, write $\rho=\rho_\beta(y)$, where $\rho_\beta$ is defined in \eqref{eq:y-rho}.

We first consider $1/2<\ell<2$. For $\rho\geq\sqrt{3}$, \eqref{eq:Uprime-rho} gives
\begin{equation*}
\ell+U'(y)=\frac{\ell\rho^2}{1+\rho^2}\geq\frac{3\ell}{4}\geq\frac{3\ell}{4}\frac{y^2}{1+y^2}.
\end{equation*}
Thus it remains to prove \eqref{eq:small-ell-first-damping} for $0<\rho<\sqrt3$. By \eqref{eq:y-rho},
\begin{equation*}
y=\frac{2}{(2\ell)^\ell\sqrt{\beta}}\int_0^\rho(1+s^2)^\ell\,ds\leq\frac{2\rho(1+\rho^2)^\ell}{(2\ell)^\ell\sqrt{\beta}}\leq\frac{2\rho}{\sqrt{\beta}}\left(\frac2\ell\right)^\ell.
\end{equation*}
Since $(2/\ell)^\ell<3$ for $1/2<\ell<2$ and $\sqrt{\beta}\geq24$ by \eqref{eq:normalized-profile-choice}, we have $y\leq\rho/4<\rho$. Hence, by \eqref{eq:Uprime-rho},
\begin{equation*}
\ell+U'(y)=\frac{\ell\rho^2}{1+\rho^2}\geq\ell\frac{y^2}{1+y^2}\geq\frac{3\ell}{4}\frac{y^2}{1+y^2}.
\end{equation*}
This proves \eqref{eq:small-ell-first-damping}. Integrating \eqref{eq:small-ell-first-damping} from $0$ to $y$ and using $U(0)=0$, we obtain
\begin{equation}\label{Uovy}
\frac{U(y)}y\geq-\ell+\frac{3\ell}{4}\left(1-\frac{\arctan y}{y}\right).
\end{equation}
Substituting \eqref{eq:small-ell-first-damping} and \eqref{Uovy} into \eqref{eq:D1} and \eqref{eq:D2}, we have
\begin{equation*}
D_{1,\ell}(y)\geq\frac12\frac{y^2}{1+y^2}+\frac{3\ell}{2(1+y^2)}\left(1-\frac{\arctan y}{y}\right)\geq\frac12\frac{y^2}{1+y^2}
\end{equation*}
and
\begin{equation*}
D_{2,\ell}(y)\geq\left(\frac{3\ell}{4}+1\right)\frac{y^2}{1+y^2}+\frac{3\ell}{4(1+y^2)}\left(1-\frac{\arctan y}{y}\right)\geq\left(\frac{3\ell}{4}+1\right)\frac{y^2}{1+y^2}.
\end{equation*}
Here we used $\arctan y\leq y$ for $y>0$. This proves \eqref{eq:small-ell-D1-lower} and \eqref{eq:small-ell-D2-lower}.

We next consider $\ell>2$, for which $\beta=1$ by \eqref{eq:normalized-profile-choice}. Since $U'(0)=-\ell$ by \eqref{eq:origin-derivatives}, we have
\begin{equation}\label{eq:first-damping-equivalent}
\ell+U'(y)-(\ell-1)\frac{y^2}{1+y^2}=\int_0^y\left(U''(y')-\frac{2(\ell-1)y'}{(1+y'^2)^2}\right)\,dy'.
\end{equation}
We first consider the region where $-\ell<U'(y)\leq-1$, or equivalently $0<\rho^2\leq\ell-1$. In this region, \eqref{eq:y-rho} gives
\begin{equation}\label{eq:damping-y-rho-bound}
y=\frac{2}{(2\ell)^{\ell}}\int_0^\rho(1+s^2)^{\ell}\,ds\leq\frac{2\rho(1+\rho^2)^{\ell}}{(2\ell)^{\ell}}.
\end{equation}
Using \eqref{eq:Uprime-rho} and \eqref{eq:damping-y-rho-bound}, we obtain
\begin{equation} \label{weighted_U''}
\frac{(1+y^2)^2}{y}U''(y)\geq\frac{U''(y)}y\geq\frac{\ell(2\ell)^{2\ell}}{2(1+\rho^2)^{2\ell+2}}\geq\frac{\ell(2\ell)^{2\ell}}{2\ell^{2\ell+2}}=\frac{2^{2\ell-1}}{\ell}\geq2\ell>2(\ell-1),
\end{equation}
where we used $2^{\ell-1}\geq\ell$ for $\ell\geq2$. Since $\rho_\beta$ is increasing, \eqref{weighted_U''} also holds with $y$ replaced by any $y'\in(0,y]$. Thus the integrand in \eqref{eq:first-damping-equivalent} is nonnegative, proving \eqref{eq:first-damping} in this region. On the other hand, when \(-1\leq U'(y)<0\),
\begin{equation*}
\ell+U'(y)-(\ell-1)\frac{y^2}{1+y^2}\geq(\ell-1)\left(1-\frac{y^2}{1+y^2}\right)\geq0.
\end{equation*}
This proves \eqref{eq:first-damping} for all \(y>0\).

It remains to prove \eqref{eq:D1-lower} and \eqref{eq:D2-lower} for $\ell>2$. By \eqref{eq:first-damping},
\begin{equation}\label{eq:Uprime-lower}
U'(y)\geq-1-\frac{\ell-1}{1+y^2}.
\end{equation}
Using \(U(0)=0\), we also obtain
\begin{equation}\label{eq:U-lower}
U(y)=\int_0^yU'(y')\,dy'\geq-y-(\ell-1)\int_0^y\frac{dy'}{1+y'^2}=-y-(\ell-1)\arctan y.
\end{equation}
Substituting \eqref{eq:Uprime-lower} and \eqref{eq:U-lower} into \eqref{eq:D1}, we have
\begin{equation*}
D_{1,\ell}(y)\geq\frac{1}{1+y^2}\left(\frac{\ell-2}{\ell}y^2+2(\ell-1)\left(1-\frac{\arctan y}{y}\right)\right)\geq\frac{\ell-2}{\ell}\frac{y^2}{1+y^2},
\end{equation*}
where we used $\arctan y\leq y$ for $y>0$. Similarly, \eqref{eq:D2} gives
\begin{equation*}
D_{2,\ell}(y)\geq\frac{1}{1+y^2}\left(\left(\ell+\frac12-\frac{2}{\ell}\right)y^2+(\ell-1)\left(1-\frac{\arctan y}{y}\right)\right)\geq\left(\ell+\frac12-\frac{2}{\ell}\right)\frac{y^2}{1+y^2}.
\end{equation*}
This completes the proof.
\end{proof}

We conclude this subsection with the proof of Proposition~\ref{prop:local-integral}.

\begin{proof}[Proof of Proposition~\ref{prop:local-integral}]
The case $\ell=2$ is proved in \cite[Lemma~6.3]{KKY-CH}. We therefore consider $\ell>1/2$ with $\ell\ne2$. By the odd symmetry of $U$ and continuity at $y=0$, it suffices to prove \eqref{eq:local-integral-inequality} for $y>0$. We begin by expressing $D_{1,\ell}$ and $yU''$ in terms of the parameter $\rho$ in \eqref{eq:y-rho}. By \eqref{eq:y-rho}--\eqref{eq:Uprime-rho}, we have
\begin{equation}\label{eq:D1-rho-local}
D_{1,\ell}(y)=1-\frac{2}{1+\rho^2}+\frac{R(\rho)}{1+y^2}=\frac{\rho^2-1}{1+\rho^2}+\frac{R(\rho)}{1+y^2},
\end{equation}
and
\begin{equation}\label{eq:yUsecond-rho-local}
yU''(y)=\frac{2\ell\rho\int_0^\rho(1+s^2)^{\ell}\,ds}{(1+\rho^2)^{\ell+2}}=\frac{2\ell\rho^2}{R(\rho)(1+\rho^2)^2},
\end{equation}
where
\begin{equation}\label{eq:R-definition-local}
R(\rho):=\frac{\rho(1+\rho^2)^{\ell}}{\int_0^\rho(1+s^2)^{\ell}\,ds}.
\end{equation}
Since \(U''(y)>0\) for \(y>0\) by \eqref{eq:Uprime-rho}, \eqref{eq:local-kernel} can be written as
\begin{equation}\label{eq:local-kernel-factorization}
K_{\ell}^{\mathrm{loc}}(y)=yU''(y)\frac{1+y^2}{y^3}\int_0^y\frac{y'^2}{1+y'^2}\,dy'.
\end{equation}
Note that, for $0\leq y'\leq y$, we have
\begin{equation*}
\frac{(y')^2}{1+(y')^2}\leq\frac{y^2}{1+y^2},\qquad \frac{(y')^2}{1+(y')^2}\leq(y')^2.
\end{equation*}
Using these bounds in \eqref{eq:local-kernel-factorization}, we obtain
\begin{equation}\label{eq:g-bounds}
0<K_{\ell}^{\mathrm{loc}}(y)\leq yU''(y),\qquad K_{\ell}^{\mathrm{loc}}(y)\leq\frac{y(1+y^2)}{3}U''(y),\qquad y>0.
\end{equation}

We next derive a lower bound for $R(\rho)$. Integration by parts gives
\begin{equation*}
R(\rho)=1+2\ell\frac{\int_0^\rho\frac{s^2}{1+s^2}(1+s^2)^\ell\,ds}{\int_0^\rho(1+s^2)^\ell\,ds}.
\end{equation*}
Since both $s^2/(1+s^2)$ and $(1+s^2)^\ell$ are increasing on $[0,\rho]$, Chebyshev's integral inequality gives
\begin{equation*}
\frac{\int_0^\rho\frac{s^2}{1+s^2}(1+s^2)^\ell\,ds}{\int_0^\rho(1+s^2)^\ell\,ds}\geq\frac1\rho\int_0^\rho\frac{s^2}{1+s^2}\,ds=1-\frac{\arctan\rho}{\rho}.
\end{equation*}
Moreover, since $1+s^2\leq1+\rho^2$ for $0\leq s\leq\rho$, we have
\begin{equation*}
1-\frac{\arctan\rho}{\rho}=\frac1\rho\int_0^\rho\frac{s^2}{1+s^2}\,ds\geq\frac1{\rho(1+\rho^2)}\int_0^\rho s^2\,ds=\frac{\rho^2}{3(1+\rho^2)}.
\end{equation*}
Combining these inequalities, we obtain
\begin{equation}\label{eq:R-lower-local}
R(\rho)\geq1+2\ell\left(1-\frac{\arctan\rho}{\rho}\right)\geq1+\frac{2\ell\rho^2}{3(1+\rho^2)}.
\end{equation}

With these preliminary estimates, we now prove \eqref{eq:local-integral-inequality}. We consider separately the two cases $1/2<\ell<2$ and $\ell>2$.

We first consider $1/2<\ell<2$ and divide the range $\rho>0$ into $0<\rho\leq1$, $1\leq\rho\leq\sqrt3$, and $\rho\geq\sqrt3$. For $0<\rho\leq\sqrt3$, \eqref{eq:y-rho} gives
\begin{equation*}
y=\frac{2}{(2\ell)^\ell\sqrt{\beta}}\int_0^\rho(1+s^2)^\ell\,ds\leq\frac{2\rho(1+\rho^2)^\ell}{(2\ell)^\ell\sqrt{\beta}}\leq\frac{2\rho}{\sqrt{\beta}}\left(\frac2\ell\right)^\ell\leq\frac{\rho}{4},
\end{equation*}
where we used $(2/\ell)^\ell<3$ for $1/2<\ell<2$ and $\sqrt{\beta}\geq24$ by \eqref{eq:normalized-profile-choice}. For $0<\rho\leq1$, \eqref{eq:D1-rho-local}, \eqref{eq:R-lower-local}, and $y^2\leq\rho^2/16$ give
\begin{equation*}
D_{1,\ell}(y)\geq\frac{(2+2\ell/3)\rho^2-(1-\rho^2)y^2}{(1+\rho^2)(1+y^2)}\geq\frac{2(\ell+2)\rho^2}{3(1+\rho^2)(1+y^2)}.
\end{equation*}
By \eqref{eq:yUsecond-rho-local} and \eqref{eq:g-bounds}, together with $R(\rho)\geq1$, we obtain
\begin{equation*}
\frac{K_\ell^{\mathrm{loc}}(y)}{D_{1,\ell}(y)}\leq\frac{\ell}{\ell+2}\frac{(1+y^2)^2}{1+\rho^2}\leq\frac{\ell}{\ell+2},
\end{equation*}
where we used $(1+\rho^2/16)^2\leq1+\rho^2$ for $0<\rho\leq1$. Next, for $1\leq\rho\leq\sqrt3$, \eqref{eq:D1-rho-local} gives $D_{1,\ell}(y)\geq R(\rho)/(1+y^2)$, while \eqref{eq:R-lower-local} gives $R(\rho)\geq1+\ell/3$. Using $\rho^2/(1+\rho^2)^2\leq1/4$ and $y^2\leq3/16$, we obtain
\begin{equation*}
\frac{K_\ell^{\mathrm{loc}}(y)}{D_{1,\ell}(y)}\leq\frac{\ell(1+y^2)^2}{6R(\rho)^2}\leq\frac{\ell}{4(1+\ell/3)^2},
\end{equation*}
where we used $(1+3/16)^2<3/2$. Finally, for $\rho\geq\sqrt3$, \eqref{eq:D1-rho-local} gives $D_{1,\ell}(y)\geq1/2$. Since $(\arctan\rho)/\rho$ is decreasing for $\rho>0$, we have $(\arctan\rho)/\rho\leq\pi/(3\sqrt3)<5/8$. Hence \eqref{eq:R-lower-local} gives $R(\rho)\geq1+3\ell/4$. Using $\rho^2/(1+\rho^2)^2\leq3/16$ and \eqref{eq:yUsecond-rho-local}--\eqref{eq:g-bounds}, we obtain
\begin{equation*}
\frac{K_\ell^{\mathrm{loc}}(y)}{D_{1,\ell}(y)}\leq\frac{3\ell}{4R(\rho)}\leq\frac{3\ell}{3\ell+4}.
\end{equation*}
By \eqref{eq:small-ell-D1-lower} and
$\max\{\frac{\ell}{\ell+2},\frac{\ell}{4(1+\ell/3)^2},\frac{3\ell}{3\ell+4}\}<\frac35$, we conclude that
\begin{equation*}
K_\ell^{\mathrm{loc}}(y)\leq\frac35D_{1,\ell}(y)\leq\frac34\left(D_{1,\ell}(y)-\frac1{10}\frac{y^2}{1+y^2}\right).
\end{equation*}
This proves \eqref{eq:local-integral-inequality} for $1/2<\ell<2$, with $\lambda_\ell^{\mathrm{loc}}=4/3$ and $\eta_\ell=1/10$.

We next consider $\ell>2$, for which $\beta=1$ by \eqref{eq:normalized-profile-choice}. As in the previous case, we divide the range $\rho>0$ into $0<\rho\leq1$, $1\leq\rho\leq\sqrt3$, and $\rho\geq\sqrt3$. For any \(0\leq s\leq\sqrt3\), we have \(((1+s^2)/(2\ell))^{\ell}\leq((1+s^2)/4)^2\), and hence \eqref{eq:y-rho} gives
\begin{equation}\label{eq:y-bound-small-n}
y\leq\frac{\rho}{8}\left(1+\frac{2\rho^2}{3}+\frac{\rho^4}{5}\right),\qquad 0\leq\rho\leq\sqrt3.
\end{equation}
For $0<\rho\leq1$, \eqref{eq:y-bound-small-n} gives $y\leq7\rho/30<\rho/4$. By \eqref{eq:D1-rho-local}, \eqref{eq:yUsecond-rho-local}, \eqref{eq:g-bounds}, and \eqref{eq:R-lower-local}, we obtain
\begin{equation*}
D_{1,\ell}(y)\geq\frac{(2+2\ell/3)\rho^2-(1-\rho^2)y^2}{(1+\rho^2)(1+y^2)}\geq\frac{2(\ell+2)\rho^2}{3(1+\rho^2)(1+y^2)},\quad K_\ell^{\mathrm{loc}}(y)\leq\frac{2\ell \rho^2(1+y^2)}{3(1+\rho^2)^2}.
\end{equation*}
Since \((1+\rho^2/16)^2\leq1+\rho^2\) for \(0\leq \rho^2\leq1\), these bounds imply
\begin{equation*}
\frac {K_\ell^{\mathrm{loc}}(y)}{D_{1,\ell}(y)}\leq\frac{\ell}{\ell+2}\frac{(1+y^2)^2}{1+\rho^2}\leq\frac{\ell}{\ell+2}.
\end{equation*}
For $1\leq\rho\leq\sqrt3$, \eqref{eq:y-bound-small-n} gives $y\leq3\sqrt3/5$, and hence $y^2\leq27/25<5/4$. Since $(\arctan\rho)/\rho$ is decreasing for $\rho>0$, \eqref{eq:R-lower-local} and $\pi/4<11/14$ give $R(\rho)\geq1+3\ell/7$. Moreover, \eqref{eq:D1-rho-local}, \eqref{eq:yUsecond-rho-local}, and \eqref{eq:g-bounds} give
\begin{equation*}
D_{1,\ell}(y)\geq\frac{R(\rho)}{1+y^2},\qquad K_\ell^{\mathrm{loc}}(y)\leq\frac{2\ell\rho^2(1+y^2)}{3R(\rho)(1+\rho^2)^2}.
\end{equation*}
Thus, using $\rho^2/(1+\rho^2)^2\leq1/4$, we obtain
\begin{equation*}
\frac {K_\ell^{\mathrm{loc}}(y)}{D_{1,\ell}(y)}\leq\frac{\ell(1+y^2)^2}{6R(\rho)^2}\leq\frac{27\ell}{32(1+3\ell/7)^2}.
\end{equation*}
Finally, for $\rho\geq\sqrt3$, \eqref{eq:D1-rho-local} gives $D_{1,\ell}(y)\geq1/2$. Since $(\arctan\rho)/\rho$ is decreasing for $\rho>0$, we have $(\arctan\rho)/\rho\leq\pi/(3\sqrt3)<5/8$. Hence \eqref{eq:R-lower-local} gives $R(\rho)\geq1+3\ell/4$. Using $\rho^2/(1+\rho^2)^2\leq3/16$ and \eqref{eq:yUsecond-rho-local}--\eqref{eq:g-bounds}, we obtain
\begin{equation*}
\frac{K_\ell^{\mathrm{loc}}(y)}{D_{1,\ell}(y)}\leq\frac{3\ell}{4R(\rho)}\leq\frac{3\ell}{3\ell+4}.
\end{equation*}
Since $\max\{\frac{\ell}{\ell+2},\frac{27\ell}{32(1+3\ell/7)^2},\frac{3\ell}{3\ell+4}\}=\frac{3\ell}{3\ell+4}$ for $\ell>2$, we obtain
\begin{equation}\label{eq:local-ratio-bound'}
K_\ell^{\mathrm{loc}}(y)\leq\frac{3\ell}{3\ell+4}D_{1,\ell}(y),\qquad y>0.
\end{equation}
To obtain \eqref{eq:local-integral-inequality}, let $d_\ell:=3\ell/(3\ell+4)\in(0,1)$ and choose
\begin{equation}\label{eq:eta-local-choice}
\eta_\ell:=\frac{\ell-2}{2\ell}(1-d_\ell)>0.
\end{equation}
By \eqref{eq:D1-lower},
\begin{equation*}
D_{1,\ell}(y)-\eta_\ell\frac{y^2}{1+y^2}\geq\left(1-\frac{\ell\eta_\ell}{\ell-2}\right)D_{1,\ell}(y)=\frac{1+d_\ell}{2}D_{1,\ell}(y).
\end{equation*}
Combining this with \eqref{eq:local-ratio-bound'}, we obtain
\begin{equation}\label{eq:subtracted-D1}
K_\ell^{\mathrm{loc}}(y)\leq\frac{2d_\ell}{1+d_\ell}\left(D_{1,\ell}(y)-\eta_\ell\frac{y^2}{1+y^2}\right).
\end{equation}
Since $d_\ell\in(0,1)$, we may choose
\begin{equation}\label{eq:local-lambda-delta-choice}
\lambda_\ell^{\mathrm{loc}}:=\frac{1+d_\ell}{2d_\ell}>1.
\end{equation}
Multiplying \eqref{eq:subtracted-D1} by $\lambda_\ell^{\mathrm{loc}}$ proves \eqref{eq:local-integral-inequality} for $\ell>2$.

The cases $y<0$ and $y=0$ follow by evenness and continuity, respectively. This completes the proof.
\end{proof}

\subsection{The far-field inequality and the choice of \texorpdfstring{$m_0(\ell)$}{m0(ell)}}\label{sec:far-field-inequality-proof}

We now prove Proposition~\ref{prop:far-field}.

\begin{proof}[Proof of Proposition~\ref{prop:far-field}]
By the odd symmetry of $U$, both sides of \eqref{eq:far-field-inequality} are even in $y$. Thus it is enough to consider $y>0$. We first compute the far-field limits needed to compare the two sides of \eqref{eq:far-field-inequality}. Write
\begin{equation*}
\gamma_{\ell}:=\frac{2\ell-1}{2\ell+1}.
\end{equation*}
By l'H\^opital's rule and \eqref{eq:Jr}, we have
\begin{equation}\label{eq:Jr-far-limit}
\lim_{y\to\infty}y^{-\gamma_{\ell}}J_\ell(y)=\lim_{y\to\infty}\frac{y^{2/(2\ell+1)}}{\gamma_{\ell}(1+y^{2/(2\ell+1)})}=\frac1{\gamma_{\ell}}.
\end{equation}
For $y>0$, \eqref{eq:Uprime-rho} gives $U''(y)>0$. Therefore, using \eqref{eq:far-kernel}, \eqref{eq:far-Usecond}, and \eqref{eq:Jr-far-limit}, we obtain
\begin{equation}\label{eq:far-kernel-limit}
\begin{split}
\lim_{y\to\infty}y^{2/(2\ell+1)}K_{\ell}^{\mathrm{far}}(y)&=\lim_{y\to\infty}(1+y^{-2/(2\ell+1)})\bigl(y^{(2\ell+3)/(2\ell+1)}U''(y)\bigr)\bigl(y^{-\gamma_{\ell}}J_\ell(y)\bigr)\\
&=\frac{\frac{2}{2\ell+1}\Theta_\ell\beta^{-1/(2\ell+1)}}{\gamma_{\ell}}=\frac{2\Theta_\ell\beta^{-1/(2\ell+1)}}{2\ell-1}.
\end{split}
\end{equation}
Next, multiplying \eqref{eq:far-damping} by $y^{2/(2\ell+1)}$, we obtain
\begin{equation}\label{eq:far-damping-rescaled}
\begin{split}
y^{2/(2\ell+1)}D_{\ell,B}^{\mathrm{far}}(y)&=\frac{(1-B/\ell)y^{2/(2\ell+1)}}{1+y^{2/(2\ell+1)}}+\frac{2y^{2/(2\ell+1)}U'(y)}{\ell}\\
&\quad-\frac{\frac{2}{2\ell+1}y^{2/(2\ell+1)}}{1+y^{2/(2\ell+1)}}\left(y^{-\gamma_{\ell}}U(y)+B y^{-\gamma_{\ell}}J_\ell(y)\right).
\end{split}
\end{equation}
By \eqref{eq:far-U} and \eqref{eq:far-Uprime}, with $\beta$ as in \eqref{eq:normalized-profile-choice},
\begin{equation*}
\lim_{y\to\infty}y^{-\gamma_{\ell}}U(y)=-\frac{\Theta_\ell\beta^{-1/(2\ell+1)}}{\gamma_{\ell}},\qquad \lim_{y\to\infty}y^{2/(2\ell+1)}U'(y)=-\Theta_\ell\beta^{-1/(2\ell+1)}.
\end{equation*}
Substituting these limits and \eqref{eq:Jr-far-limit} into \eqref{eq:far-damping-rescaled}, we obtain
\begin{equation}\label{eq:far-damping-limit}
\lim_{y\to\infty}y^{2/(2\ell+1)}D_{\ell,B}^{\mathrm{far}}(y)=1-\frac{B}{\ell}-\frac{2\Theta_\ell\beta^{-1/(2\ell+1)}}{\ell}+\frac{2}{2\ell-1}(\Theta_\ell\beta^{-1/(2\ell+1)}-B).
\end{equation}
Combining \eqref{eq:far-kernel-limit} and \eqref{eq:far-damping-limit}, we find that
\begin{equation}\label{eq:far-gap}
\begin{split}
g_{\ell}&:=\lim_{y\to\infty}y^{2/(2\ell+1)}\left(D_{\ell,B_{\ell}}^{\mathrm{far}}(y)-\Lambda_{\ell}K_{\ell}^{\mathrm{far}}(y)\right)\\
&=1-\frac{2\Theta_\ell\beta^{-1/(2\ell+1)}}{\ell}-B_{\ell}\left(\frac1{\ell}+\frac2{2\ell-1}\right)-(\Lambda_{\ell}-1)\frac{2\Theta_\ell\beta^{-1/(2\ell+1)}}{2\ell-1}.
\end{split}
\end{equation}

We claim that $g_\ell>0$ for all $\ell>1/2$. We consider separately the four cases $1/2<\ell<2$, $\ell=2$, $2<\ell<4$, and $\ell\geq4$. First, by \eqref{eq:theta-ell},
\begin{equation}\label{eq:theta-less-half}
\Theta_\ell^{2\ell+1}=\frac{\ell}{2^{2\ell-2}(2\ell+1)^2}=2^{-(2\ell+1)}\frac{8\ell}{(2\ell+1)^2}<2^{-(2\ell+1)},
\end{equation}
since $(2\ell+1)^2-8\ell=(2\ell-1)^2>0$. Thus $\Theta_\ell<1/2$. For $1/2<\ell<2$, substituting \eqref{eq:B-choice} into \eqref{eq:far-gap} and using \eqref{eq:normalized-profile-choice}, we obtain
\begin{equation*}
g_\ell=1-\Theta_\ell\beta^{-1/(2\ell+1)}\left(\frac4\ell+\frac6{2\ell-1}\right)\geq\frac12.
\end{equation*}
For $\ell\geq2$, we have $\beta=1$. For $\ell=2$, we use the sharper bound $\Theta_2=50^{-1/5}<11/24$, which follows from $24^5<50\cdot11^5$. By \eqref{eq:B-choice},
\begin{equation*}
g_2=\frac6{13}-\Theta_2-\frac1{10000}\frac{2\Theta_2}{3}>\frac6{13}-\frac{15001}{15000}\frac{11}{24}>0.
\end{equation*}
For $2<\ell<4$, we first note that $\Theta_\ell/\ell$ is decreasing for $\ell>1/2$. Indeed, by \eqref{eq:theta-ell},
\begin{equation*}
\frac{d}{d\ell}\log\frac{\Theta_\ell}{\ell}=\frac{-4\ell-6+2\log((2\ell+1)^2/(8\ell))}{(2\ell+1)^2}<0,
\end{equation*}
where we used $(2\ell+1)^2/(8\ell)\leq2\ell$ and $\log(2\ell)\leq2\ell-1$. Hence, for $\ell>2$,
\begin{equation} \label{Theta_w}
\Theta_\ell\left(\frac3\ell+\frac2{2\ell-1}\right)=\frac{\Theta_\ell}{\ell}\left(4+\frac1{2\ell-1}\right)\leq\frac{13}{6}\Theta_2<\frac{143}{144}<1.
\end{equation}
Substituting \eqref{eq:B-choice} into \eqref{eq:far-gap} gives
\begin{equation*}
g_\ell=\frac13\left[1-\Theta_\ell\left(\frac3\ell+\frac2{2\ell-1}\right)\right]\left(3-\frac1\ell-\frac{2(1+\Theta_\ell)}{2\ell-1}\right)>\frac12\left[1-\Theta_\ell\left(\frac3\ell+\frac2{2\ell-1}\right)\right]>0,
\end{equation*}
where we used $\Theta_\ell<1/2$ and $\ell>2$ for the first inequality, and \eqref{Theta_w} for the second. Finally, for $\ell\geq4$, \eqref{eq:B-choice}, \eqref{eq:far-gap}, and \eqref{eq:theta-less-half} give
\begin{equation*}
g_\ell>1-\frac2\ell-\frac2{2\ell-1}-\frac1{62(2\ell-1)}\geq\frac3{14}-\frac1{434}=\frac{46}{217}>\frac15.
\end{equation*}
This proves the claim. It follows from \eqref{eq:far-gap} that \eqref{eq:far-field-inequality} holds for all sufficiently large $y>0$. Moreover, the choices in \eqref{eq:B-choice} satisfy $\Lambda_\ell>1$ for every $\ell>1/2$, since $\Theta_\ell(3/\ell+2/(2\ell-1))<1$ for $2<\ell<4$.

To specify \(m_0(\ell)\), set
\begin{equation*}
G_{\ell}(\rho):=y_{\ell}(\rho)^{2/(2\ell+1)}\left(D_{\ell,B_{\ell}}^{\mathrm{far}}(y_{\ell}(\rho))-\Lambda_{\ell}K_{\ell}^{\mathrm{far}}(y_{\ell}(\rho))\right),
\end{equation*}
where \(y_\ell:=\rho_\beta^{-1}\), with \(\beta\) given by \eqref{eq:normalized-profile-choice}; see \eqref{eq:y-rho}. By \eqref{eq:far-gap}, \(G_{\ell}(\rho)\to g_{\ell}>0\) as \(\rho\to\infty\). We may therefore define
\begin{equation}\label{eq:rho0-definition}
\rho_0(\ell):=\min\left\{j\in\N:j\geq1\text{ and }G_{\ell}(\rho)\geq\frac{g_{\ell}}{2}\text{ for all }\rho\geq j\right\}.
\end{equation}
The set in \eqref{eq:rho0-definition} is nonempty. Since \(y_{\ell}\) is strictly increasing and onto \([0,\infty)\), the choice
\begin{equation}\label{eq:m0-definition}
m_0(\ell):=y_{\ell}(\rho_0(\ell))
\end{equation}
gives \(y^{2/(2\ell+1)}(D_{\ell,B_{\ell}}^{\mathrm{far}}(y)-\Lambda_{\ell}K_{\ell}^{\mathrm{far}}(y))\geq g_{\ell}/2>0\) for all \(y\geq m_0(\ell)\). By evenness, \eqref{eq:far-field-inequality} holds for all \(|y|\geq m_0(\ell)\). This completes the proof.
\end{proof}

\begin{remark}\label{rmk_beta}
The choice of $\beta$ in \eqref{eq:normalized-profile-choice}, together with $B_\ell$ and $\Lambda_\ell$ in \eqref{eq:B-choice}, ensures the following for every $\ell>1/2$: (i) the estimates in Propositions~\ref{prop:damping-inequalities} and~\ref{prop:local-integral}; (ii) $g_\ell>0$ in the proof of Proposition~\ref{prop:far-field}; and (iii) $B_\ell>\Theta_\ell\beta^{-1/(2\ell+1)}$ and $\Lambda_\ell>1$. The first inequality in (iii) allows a positive $\vartheta$ in \eqref{eq:common-singularity-theta}, while the second is used in Lemma~\ref{lem:bootstrap-far} to absorb the $K_\ell^{\mathrm{far}}$ term into $D_{\ell,B_\ell}^{\mathrm{far}}$.

For any fixed $\beta>0$, the conditions in (ii) and (iii) cannot hold for all $\ell>1/2$. Indeed, by \eqref{eq:far-gap} and (iii),
\begin{equation*}
g_\ell<1-\Theta_\ell\beta^{-1/(2\ell+1)}\left(\frac3\ell+\frac2{2\ell-1}\right)\longrightarrow-\infty\qquad\text{as }\ell\downarrow\frac12.
\end{equation*}
This explains the dependence of $\beta$ on $\ell$ in \eqref{eq:normalized-profile-choice}.
\end{remark}

\section{Singularity formation for the nonlocal transport equations}\label{sec:b-family}

We now prove Theorem~\ref{thm:b-cusp-intro}. The equation and its nonlocal formulation are given in \eqref{eq:b-family}--\eqref{eq:b-slope}. We first specify a class of initial data and state the result for this class. We then verify the hypotheses of Theorem~\ref{thm:common-singularity} and apply it to complete the proof. Throughout this section, $1<b<5$ and $g\in C^\infty(\mathbb R)$ with $g(0)=0$ are fixed.

\subsection{Initial conditions and the main result}\label{subsec:b-theorem}

Let $U=U_\beta$ be the smooth odd solution of the profile equation \eqref{eq:profile} corresponding to $\ell=2/(b-1)$, with $\beta$ given by \eqref{eq:normalized-profile-choice}. Choose $m,\vartheta,N$ as in Section~\ref{sec:common-bootstrap} and Remark~\ref{rem:common-singularity-parameters}, and set $\Gamma_3=\Gamma_4=0$. With $M_0$ given by Proposition~\ref{prop:common-bootstrap}, choose $M\geq M_0$ sufficiently large that \eqref{eq:bootstrap-initial-second}, \eqref{eq:bootstrap-initial-third}, and \eqref{eq:bootstrap-initial-fourth} hold with strict margins when $W(\cdot,s_0)$ is replaced by $U$.

Let $\varepsilon>0$ be a parameter to be chosen sufficiently small. We consider smooth initial data $w_0\in H^r(\mathbb R)$ for every $r\geq0$, and assume that $(1+x^2)^{1/2}w_0\in H^5(\mathbb R)$. We impose the matching conditions
\begin{equation}\label{eq:b-initial-matching}
w_0'(0)=-\frac{\ell}{\varepsilon},\qquad w_0''(0)=0,
\end{equation}
and require
\begin{equation}\label{eq:b-initial-origin}
\left|\varepsilon^{2\ell+2}w_0^{(3)}(0)-U^{(3)}(0)\right|\leq1-2\vartheta.
\end{equation}
The derivatives satisfy
\begin{subequations}\label{eq:b-initial-derivatives}
\begin{align}
&|w_0''(x)|\leq(1-2\vartheta)M^{1/8}\varepsilon^{-(2\ell+3)/2}\frac{|x|}{\sqrt{\varepsilon^{2\ell+1}+x^2}},\qquad x\in\mathbb R,\label{eq:b-initial-second}\\
&\|w_0^{(3)}\|_{L^\infty}\leq(1-2\vartheta)M^{3/4}\varepsilon^{-(2\ell+2)},\label{eq:b-initial-third}\\
&\|w_0^{(4)}\|_{L^\infty}\leq(1-2\vartheta)M\varepsilon^{-(6\ell+5)/2}.\label{eq:b-initial-fourth}
\end{align}
\end{subequations}
We assume that the rescaled initial slope is close to $U'$:
\begin{equation}\label{eq:b-initial-slope}
\left|\varepsilon w_0'(x)-U'\left(\frac{x}{\varepsilon^{(2\ell+1)/2}}\right)\right|\leq(1-2\vartheta)\min\left\{\frac{Nx^2}{\varepsilon^{2\ell+1}+x^2},\frac{\varepsilon B_\ell}{\varepsilon+|x|^{2/(2\ell+1)}}\right\}
\end{equation}
for all $x\in\mathbb R$. Finally, we impose
\begin{equation}\label{eq:b-initial-energy}
\|(1+x^2)^{1/(2\ell+1)}w_0\|_{L^\infty}\leq\frac14.
\end{equation}

The compatibility of these assumptions is established in Appendix~\ref{subsec:b-data}, where we construct admissible initial data. For each sufficiently small $\varepsilon>0$, the admissible set has non-empty interior relative to the class of smooth data satisfying the above regularity assumptions and \eqref{eq:b-initial-matching}, equipped with the weighted $H^5$ topology.

We now state the main result for these initial data.

\begin{theorem}\label{thm:b-cusp}
Fix $1<b<5$ and $g\in C^\infty(\mathbb R)$ with $g(0)=0$. Let $m,\vartheta,N,M$ be as specified in this subsection. There exists $\varepsilon_0>0$ such that the following statement holds. 

For each $0<\varepsilon\leq\varepsilon_0$, suppose that the initial data $w_0$ satisfy
\begin{equation*}
w_0\in H^r(\mathbb R)\quad\text{for every }r\geq0,\qquad (1+x^2)^{1/2}w_0\in H^5(\mathbb R),
\end{equation*}
together with \eqref{eq:b-initial-matching}--\eqref{eq:b-initial-energy}. Then the initial value problem \eqref{eq:b-family}--\eqref{eq:b-initial-data-intro} admits a unique smooth solution
\begin{equation*}
w\in C([0,T_*);H^5(\mathbb R))\cap C^1([0,T_*);H^4(\mathbb R)),
\end{equation*}
which blows up in the $C^1$-norm at $t=T_*$ for some $T_*\in(0,2\varepsilon]$. Moreover, assertions \emph{(i)}--\emph{(iii)} of Theorem~\ref{thm:b-cusp-intro} hold for $w$, where $x_*\in\mathbb R$ is the blow-up location.
\end{theorem}

The remainder of this section is devoted to the proof of Theorem~\ref{thm:b-cusp}.

\subsection{Local existence and blow-up criterion}\label{subsec:b-energy}

We first establish local existence and the continuation property for the nonlocal equations.

\begin{lemma}\label{lem:b-local}
Let \(w_0\in H^r(\mathbb R)\) for every \(r\geq0\), and suppose that \((1+x^2)^{1/2}w_0\in H^5(\mathbb R)\). There exists a maximal existence time \(T_{\max}>0\) such that the initial value problem \eqref{eq:b-family}--\eqref{eq:b-initial-data-intro} admits a unique solution \(w\) such that
\begin{equation*}
w\in C([0,T_{\max});H^r(\mathbb R))\cap C^1([0,T_{\max});H^{r-1}(\mathbb R)) \quad \text{for every } r \geq 0
\end{equation*}
and
\begin{equation*}
(1+x^2)^{1/2}w\in C([0,T_{\max});H^5(\mathbb R)).
\end{equation*}
Moreover, the following blow-up criterion holds:
\begin{equation}\label{eq:b-local-continuation}
T_{\max}<\infty\quad\Longrightarrow\quad\limsup_{t\nearrow T_{\max}}\|w_x(\cdot,t)\|_{L^\infty}=\infty.
\end{equation}
\end{lemma}

The proof of Lemma~\ref{lem:b-local} is given in Appendix~\ref{app:b-local}.

\subsection{Preliminary estimates}\label{subsec:b-preliminary}

We establish bounds for $P$, $P_x$, and $F$, together with an estimate for $w$. These estimates will be used in the proof of Theorem~\ref{thm:b-cusp}.

For the estimates below, we define
\begin{equation}\label{eq:b-pressure-rate}
\mathcal H_\ell(r):=
\begin{cases}
1,&\ell>3/2\\
1-\log r,&\ell=3/2\\
r^{\ell-3/2},&1/2<\ell<3/2
\end{cases}
\end{equation}
for $0<r\leq1$. For each $\ell>1/2$, we have
\begin{equation}\label{eq:b-pressure-rate-limit}
r\mathcal H_\ell(r)\longrightarrow0\quad\text{as }r\downarrow0
\end{equation}
and, for $0<\varepsilon\leq1$,
\begin{equation}\label{eq:b-pressure-rate-smallness}
\int_0^\varepsilon\mathcal H_\ell(r)\,dr=
\begin{cases}
\varepsilon,&\ell>3/2,\\
\varepsilon(2-\log\varepsilon),&\ell=3/2,\\
\frac{\varepsilon^{\ell-1/2}}{\ell-1/2},&1/2<\ell<3/2.
\end{cases}
\end{equation}
In particular, the integral tends to zero as $\varepsilon\downarrow0$.

\begin{lemma}\label{lem:b-pressure}
Let $0<\varepsilon<1/2$, let $w$ be a smooth solution of \eqref{eq:b-family}, and let $W$ be defined by \eqref{eq:modulated-variables}--\eqref{eq:modulated-unknown}. Assume that $s=s(t)$ is strictly increasing, and denote its inverse by $t=t(s)$. Suppose that \eqref{eq:common-bootstrap-assumptions} holds on a finite interval $[s_0,s_1]$ and that
\begin{equation}\label{eq:b-auxiliary-bootstrap}
0<\tau(t)-t\leq2\varepsilon,\qquad |\xi(t)|\leq1,\qquad \lVert (1+x^2)^{1/(2\ell+1)}w(\cdot,t)\rVert_{L^\infty}\leq1
\end{equation}
for all $t\in[t(s_0),t(s_1)]$. 
Then there exists a constant $C>0$, independent of $\varepsilon$ and $s_1$, such that
\begin{equation}\label{eq:b-physical-slope}
|w_x(x,t)|\leq\frac{C}{\tau(t)-t+|x-\xi(t)|^{2/(2\ell+1)}}
\end{equation}
for all $x\in\mathbb R$ and $t\in[t(s_0),t(s_1)]$, and
\begin{equation}\label{eq:b-weighted-pressure}
\begin{split}
&\lVert (1+x^2)^{1/(2\ell+1)}F\rVert_{L^\infty}+\lVert (1+x^2)^{1/(2\ell+1)}P\rVert_{L^\infty}+\lVert (1+x^2)^{1/(2\ell+1)}P_x\rVert_{L^\infty}\\
&\qquad\leq C\mathcal H_\ell(\tau(t)-t)
\end{split}
\end{equation}
for all $t\in[t(s_0),t(s_1)]$.
\end{lemma}

\begin{proof}
By \eqref{eq:boot-far}, \eqref{eq:far-Uprime}, and Lemma~\ref{lem:derivative-bounds}, we have
\begin{equation*}
|W_y(y,s)|\leq|W_y(y,s)-U'(y)|+|U'(y)|\leq\frac{C}{1+|y|^{2/(2\ell+1)}}.
\end{equation*}
It follows from \eqref{eq:modulated-variables}--\eqref{eq:modulated-unknown} that
\begin{equation*}
\begin{split}
|w_x(x,t)|&=\frac{|W_y(y,s)|}{\tau(t)-t} \leq \frac{C}{(\tau(t)-t)(1+|y|^{2/(2\ell+1)})} =\frac{C}{\tau(t)-t+|x-\xi(t)|^{2/(2\ell+1)}},
\end{split}
\end{equation*}
which proves \eqref{eq:b-physical-slope}.

To estimate the terms on the left-hand side of \eqref{eq:b-weighted-pressure}, we recall from \eqref{P-nonlocal} and \eqref{eq:b-slope} that
\begin{equation*}
P=K*\left(g(w)+\frac{3-b}{2}w_x^2\right),\quad P_x=K'*\left(g(w)+\frac{3-b}{2}w_x^2\right),\quad F=g(w)-P.
\end{equation*}
We first estimate the terms involving $g(w)$. For $x,z\in\mathbb R$, we have
\begin{equation}\label{eq:b-weight-product}
(1+x^2)^{1/(2\ell+1)}\leq 2^{1/(2\ell+1)}(1+(x-z)^2)^{1/(2\ell+1)}(1+z^2)^{1/(2\ell+1)}.
\end{equation}
Since $g(0)=0$ and $\|w(\cdot,t)\|_{L^\infty}\leq1$ by \eqref{eq:b-auxiliary-bootstrap}, the mean value theorem gives
\begin{equation} \label{g-ptwsbd}
|g(w(x,t))|\leq\|g'\|_{L^\infty([-1,1])}|w(x,t)|\leq C|w(x,t)|.
\end{equation}
Using this bound, \eqref{eq:b-weight-product}, and $K(x)=|K'(x)|=\frac12e^{-|x|}$ for $x\neq0$, we obtain
\begin{equation} \label{Kg-bd}
\begin{split}
&(1+x^2)^{1/(2\ell+1)} \big(|(K*g(w))(x,t)|+|(K'*g(w))(x,t)|\big)\\
&\quad\leq C\int_{\mathbb R}(1+(x-z)^2)^{1/(2\ell+1)}e^{-|x-z|}(1+z^2)^{1/(2\ell+1)}|w(z,t)|\,dz\\
&\quad\leq C\|(1+x^2)^{1/(2\ell+1)}w(\cdot,t)\|_{L^\infty}\int_{\mathbb R}(1+z^2)^{1/(2\ell+1)}e^{-|z|}\,dz\\
&\quad\leq C\|(1+x^2)^{1/(2\ell+1)}w(\cdot,t)\|_{L^\infty},
\end{split}
\end{equation}
where we used
\begin{equation*}
\int_{\mathbb R}(1+z^2)^{1/(2\ell+1)}e^{-|z|}\,dz\leq2\int_0^\infty(1+z)e^{-z}\,dz=4.
\end{equation*}
Hence, by \eqref{g-ptwsbd} and \eqref{Kg-bd},
\begin{equation}\label{eq:b-amplitude-pressure}
\begin{split}
&\|(1+x^2)^{1/(2\ell+1)}g(w)\|_{L^\infty}+\|(1+x^2)^{1/(2\ell+1)}(K*g(w))\|_{L^\infty}\\
&\quad+\|(1+x^2)^{1/(2\ell+1)}(K'*g(w))\|_{L^\infty}\leq C\|(1+x^2)^{1/(2\ell+1)}w\|_{L^\infty}\leq C,
\end{split}
\end{equation}
where the last inequality follows from \eqref{eq:b-auxiliary-bootstrap}.

We next estimate $K*w_x^2$ and $K'*w_x^2$. Using \eqref{eq:b-weight-product} and the formula for $K$ in \eqref{P-nonlocal}, we have
\begin{equation}\label{eq:b-quadratic-split}
\begin{split}
&(1+x^2)^{1/(2\ell+1)}\left(|(K*w_x^2)(x,t)|+|(K'*w_x^2)(x,t)|\right)\\
&\quad\leq C\int_{|z-\xi(t)|\leq2}(1+(x-z)^2)^{1/(2\ell+1)}e^{-|x-z|}(1+z^2)^{1/(2\ell+1)}|w_x(z,t)|^2\,dz\\
&\qquad+C\int_{|z-\xi(t)|>2}(1+(x-z)^2)^{1/(2\ell+1)}e^{-|x-z|}(1+z^2)^{1/(2\ell+1)}|w_x(z,t)|^2\,dz\\
&\quad=:I+II.
\end{split}
\end{equation}
Note that if $|z-\xi(t)|\leq2$, then $|z|\leq|z-\xi(t)|+|\xi(t)|\leq3$ by \eqref{eq:b-auxiliary-bootstrap}. Thus, since $\ell>1/2$,
\begin{equation*}
(1+z^2)^{1/(2\ell+1)}\leq C,\qquad (1+(x-z)^2)^{1/(2\ell+1)}e^{-|x-z|}\leq1.
\end{equation*}
Together with \eqref{eq:b-physical-slope}, these bounds yield
\begin{equation*}
\begin{split}
I&\leq C\int_{|z-\xi(t)|\leq2}|w_x(z,t)|^2\,dz \leq C\int_0^2\frac{1}{\bigl(\tau(t)-t+r^{2/(2\ell+1)}\bigr)^2}\,dr.
\end{split}
\end{equation*}
Making the change of variables $r^{2/(2\ell+1)}=(\tau(t)-t)\zeta$ and using $0<\tau(t)-t<1$ from \eqref{eq:b-auxiliary-bootstrap}, we obtain
\begin{equation*}
\begin{split}
\int_0^2\frac{1}{\bigl(\tau(t)-t+r^{2/(2\ell+1)}\bigr)^2}\,dr&=\frac{2\ell+1}{2}(\tau(t)-t)^{\ell-3/2}\int_0^{2^{2/(2\ell+1)}/(\tau(t)-t)}\frac{\zeta^{\ell-1/2}}{(1+\zeta)^2}\,d\zeta\\
&\leq C(\tau(t)-t)^{\ell-3/2}\left(1+\int_1^{2^{2/(2\ell+1)}/(\tau(t)-t)}\zeta^{\ell-5/2}\,d\zeta\right)\\
&\leq C\begin{cases}
1,&\ell>3/2\\
1-\log(\tau(t)-t),&\ell=3/2\\
(\tau(t)-t)^{\ell-3/2},&1/2<\ell<3/2
\end{cases}\\
&=C\mathcal H_\ell(\tau(t)-t),
\end{split}
\end{equation*}
where the last equality follows from \eqref{eq:b-pressure-rate}. Hence $I\leq C\mathcal H_\ell(\tau(t)-t)$.

For the term $II$, we have $|z|\leq\frac32|z-\xi(t)|$ whenever $|z-\xi(t)|>2$ by \eqref{eq:b-auxiliary-bootstrap}. Hence \eqref{eq:b-physical-slope} gives
\begin{equation*}
(1+z^2)^{1/(2\ell+1)}|w_x(z,t)|^2\leq\frac{C(1+z^2)^{1/(2\ell+1)}}{|z-\xi(t)|^{4/(2\ell+1)}}\leq\frac{C}{|z-\xi(t)|^{2/(2\ell+1)}}\leq C.
\end{equation*}
Substituting this bound into $II$, we have
\begin{equation*}
\begin{split}
II&\leq C\int_{|z-\xi(t)|>2}(1+(x-z)^2)^{1/(2\ell+1)}e^{-|x-z|}\,dz\\
&\leq C\int_{\mathbb R}(1+(x-z)^2)^{1/(2\ell+1)}e^{-|x-z|}\,dz\\
&=C\int_{\mathbb R}(1+z^2)^{1/(2\ell+1)}e^{-|z|}\,dz\leq C.
\end{split}
\end{equation*}
Combining the bounds for $I$ and $II$ with \eqref{eq:b-quadratic-split}, and using $\mathcal H_\ell(\tau(t)-t)\geq1$, we obtain
\begin{equation}\label{eq:b-quadratic-pressure}
\begin{split}
&\|(1+x^2)^{1/(2\ell+1)}(K*w_x^2)\|_{L^\infty}+\|(1+x^2)^{1/(2\ell+1)}(K'*w_x^2)\|_{L^\infty} \leq C\mathcal H_\ell(\tau(t)-t).
\end{split}
\end{equation}

Finally, applying \eqref{eq:b-amplitude-pressure} and \eqref{eq:b-quadratic-pressure} to \eqref{P-nonlocal}, we obtain
\begin{equation*}
\| (1+x^2)^{1/(2\ell+1)}P\|_{L^\infty}+\|(1+x^2)^{1/(2\ell+1)}P_x\|_{L^\infty}\leq C\mathcal H_\ell(\tau(t)-t).
\end{equation*}
Since $F=g(w)-P$ by \eqref{eq:b-slope}, this bound and \eqref{eq:b-amplitude-pressure} give
\begin{equation*}
\begin{split}
\|(1+x^2)^{1/(2\ell+1)}F\|_{L^\infty}&\leq\|(1+x^2)^{1/(2\ell+1)}g(w)\|_{L^\infty}+\|(1+x^2)^{1/(2\ell+1)}P\|_{L^\infty}\\
&\leq C\mathcal H_\ell(\tau(t)-t).
\end{split}
\end{equation*}
This proves \eqref{eq:b-weighted-pressure}.
\end{proof}

\begin{lemma}\label{lem:b-amplitude-estimate}
Under the assumptions of Lemma~\ref{lem:b-pressure}, suppose further that $t(s_0)=0$ and $w(\cdot,0)=w_0$. Then there exists a constant $C>0$, independent of $\varepsilon$ and $s_1$, such that
\begin{equation}\label{eq:b-weighted-amplitude}
\|(1+x^2)^{1/(2\ell+1)}w(\cdot,t)\|_{L^\infty}\leq e^{Ct}\left(\|(1+x^2)^{1/(2\ell+1)}w_0\|_{L^\infty}+C\int_0^t\mathcal H_\ell(\tau(t')-t')\,dt'\right)
\end{equation}
for all $t\in[0,t(s_1)]$.
\end{lemma}

\begin{proof}
Let $\psi$ be the characteristic curve satisfying
\begin{equation*}
\partial_t\psi(x,t)=w(\psi(x,t),t),\qquad \psi(x,0)=x.
\end{equation*}
The bounds on $w$ and $w_x$ in Lemma~\ref{lem:b-pressure} imply that $\psi(\cdot,t)$ is a diffeomorphism of $\mathbb R$ for each $t\in[0,t(s_1)]$. Multiplying \eqref{u-nonlocal} by $(1+x^2)^{1/(2\ell+1)}$, we obtain
\begin{equation}\label{eq:b-amplitude-transport}
\begin{split}
& \partial_t\bigl((1+x^2)^{1/(2\ell+1)}w\bigr)+w\partial_x\bigl((1+x^2)^{1/(2\ell+1)}w\bigr) \\
& \qquad =\frac{2x}{2\ell+1}(1+x^2)^{-2\ell/(2\ell+1)}w^2 -(1+x^2)^{1/(2\ell+1)}P_x.
\end{split}
\end{equation}
Since $\|w(\cdot,t)\|_{L^\infty}\leq1$ by \eqref{eq:b-auxiliary-bootstrap}, we have
\begin{equation*}
\begin{split}
\left|\frac{2x}{2\ell+1}(1+x^2)^{-2\ell/(2\ell+1)}w^2\right|&=\frac{2|x||w|}{(2\ell+1)(1+x^2)}(1+x^2)^{1/(2\ell+1)}|w|\\
&\leq C(1+x^2)^{1/(2\ell+1)}|w|.
\end{split}
\end{equation*}
Integrating \eqref{eq:b-amplitude-transport} along $\psi$, taking the supremum over $x\in\mathbb R$, and using \eqref{eq:b-weighted-pressure}, we obtain
\begin{equation*}
\begin{split}
\|(1+x^2)^{1/(2\ell+1)}w(\cdot,t)\|_{L^\infty}&\leq\|(1+x^2)^{1/(2\ell+1)}w_0\|_{L^\infty} +C\int_0^t\|(1+x^2)^{1/(2\ell+1)}w(\cdot,t')\|_{L^\infty}\,dt'\\
&\quad+C\int_0^t\mathcal H_\ell(\tau(t')-t')\,dt'.
\end{split}
\end{equation*}
The estimate \eqref{eq:b-weighted-amplitude} then follows from Gronwall's inequality. This completes the proof.

\end{proof}

\subsection{Proof of Theorem~\ref{thm:b-cusp}}\label{subsec:b-remainders}

To apply Theorem~\ref{thm:common-singularity}, we verify its hypotheses below.

\subsubsection{Self-similar formulation}\label{subsubsec:b-modulated}

We begin by verifying assumption~\emph{(i)} of Theorem~\ref{thm:common-singularity}, concerning local existence and the self-similar formulation. For related self-similar formulations, see \cite[Section~3.2]{KKY-CH}, \cite[Section~2.1]{KKS-rSV}, and \cite[Section~3.1]{KKKS-Novikov}.

By Lemma~\ref{lem:b-local}, the initial value problem \eqref{eq:b-family}--\eqref{eq:b-initial-data-intro} admits a unique smooth solution
\begin{equation} \label{w_reg}
w\in C([0,T];H^5(\mathbb R))\cap C^1([0,T];H^4(\mathbb R))
\end{equation}
for some $T>0$. Moreover, the solution satisfies
\begin{equation} \label{weighted_reg}
(1+x^2)^{1/2}w\in C([0,T];H^5(\mathbb R)).
\end{equation}

We define the modulation functions $\tau,\xi,\kappa$ by the system of ODEs
\begin{equation}\label{eq:b-modulation-tau}
\begin{aligned}
\dot\tau(t)&=\frac{(\tau(t)-t)^2}{\ell}F(\xi(t),t),\\
\dot\xi(t)&=\kappa(t)-\frac{F_x(\xi(t),t)}{\partial_x^3w(\xi(t),t)},\\
\dot\kappa(t)&=-P_x(\xi(t),t)+\frac{\ell F_x(\xi(t),t)}{(\tau(t)-t)\partial_x^3w(\xi(t),t)}
\end{aligned}
\end{equation}
with the initial data
\begin{equation*}
\tau(0)=\varepsilon,\qquad \xi(0)=0,\qquad \kappa(0)=w_0(0).
\end{equation*}
These equations are obtained by requiring the normalization \eqref{eq:modulation-constraints-intro} to be preserved under the evolution, as verified below. By \eqref{eq:b-initial-origin}, \eqref{eq:origin-derivatives}, and \eqref{eq:normalized-profile-choice}, we have \(w_0^{(3)}(0)>0\). By Lemma~\ref{lem:b-local} and Sobolev embedding, the functions \(F\), \(F_x\), \(P_x\), and \(\partial_x^3w\) are continuous in \((x,t)\) and locally Lipschitz in \(x\) on \(\mathbb R\times[0,T]\). Hence the right-hand side of \eqref{eq:b-modulation-tau} is continuous in \(t\) and locally Lipschitz in \((\tau,\xi,\kappa)\) wherever \(\tau-t>0\) and \(\partial_x^3w(\xi,t)>0\). Since these inequalities hold at \(t=0\), the ODE system has a unique local solution. After decreasing \(T>0\) if necessary, this solution is defined on \([0,T]\) and satisfies both inequalities there. The smoothness of \(w\) and the modulation equations imply that \(\tau\), \(\xi\), and \(\kappa\) are smooth.

To show that $s(t)=-\log(\tau(t)-t)$ is strictly increasing, we first estimate $\dot\tau(0)$. Since $\dot\tau(0)=\varepsilon^2F(0,0)/\ell$ by \eqref{eq:b-modulation-tau}, we use the initial data to bound $F(0,0)$. The initial bound \eqref{eq:b-initial-slope}, together with \eqref{eq:far-Uprime} and Lemma~\ref{lem:derivative-bounds}, gives
\begin{equation} \label{w0'bd}
|w_0'(x)|\leq\frac{C}{\varepsilon+|x|^{2/(2\ell+1)}}.
\end{equation}
Since $F=g(w)-P$ and $P=K*(g(w)+\frac{3-b}{2}w_x^2)$ with $K(x)=\frac12e^{-|x|}$, we have
\begin{equation*}
\begin{aligned}
|F(0,0)|
&=\left|g(w_0(0))-\frac12\int_{\mathbb R}e^{-|z|}\left(g(w_0(z))+\frac{3-b}{2}(w_0'(z))^2\right)\,dz\right|\\
&\le C+C\int_0^\infty\frac{e^{-r}}{(\varepsilon+r^{2/(2\ell+1)})^2}\,dr\\
&\le C\left(1+\frac{1}{\varepsilon^2}\int_0^{\varepsilon^{(2\ell+1)/2}}dr+\int_{\varepsilon^{(2\ell+1)/2}}^1r^{-4/(2\ell+1)}\,dr\right)\\
&\le C\left(1+\varepsilon^{\ell-3/2}+\int_{\varepsilon^{(2\ell+1)/2}}^1r^{-4/(2\ell+1)}\,dr\right)
\le C\mathcal H_\ell(\varepsilon),
\end{aligned}
\end{equation*}
Here we used \eqref{eq:b-initial-energy}, $g(0)=0$, and \eqref{w0'bd}; the last inequality follows from \eqref{eq:b-pressure-rate}. Consequently,
\begin{equation*}
|\dot\tau(0)|=\frac{\varepsilon^2}{\ell}|F(0,0)|\le C\varepsilon^2\mathcal H_\ell(\varepsilon)<\frac14,
\end{equation*}
where the last inequality holds by \eqref{eq:b-pressure-rate-limit} for sufficiently small $\varepsilon>0$. By continuity, decreasing $T$ further if necessary, we have $|\dot\tau(t)|<1/2$ on $[0,T]$. Hence
\begin{equation*}
\frac{ds}{dt}=\frac{1-\dot\tau(t)}{\tau(t)-t}>0,\qquad t\in[0,T].
\end{equation*}

Thus the change of variables \eqref{eq:modulated-variables}--\eqref{eq:modulated-unknown} is well-defined on $[0,T]$. Set $s_0:=s(0)=-\log\varepsilon$, and denote the inverse of $s=s(t)$ by $t=t(s)$. The function $W$ is defined on $\mathbb R\times[s_0,s(T)]$ by
\begin{equation*}
w(x,t)-\kappa(t)=(\tau(t)-t)^{(2\ell-1)/2}W\left(\frac{x-\xi(t)}{(\tau(t)-t)^{(2\ell+1)/2}},-\log(\tau(t)-t)\right).
\end{equation*}
The change of variables \eqref{eq:modulated-variables} gives
\begin{equation*}
\partial_t =(1-\dot\tau)e^s\partial_s+\left(-\dot\xi e^{(2\ell+1)s/2}+\frac{2\ell+1}{2}y(1-\dot\tau)e^s\right)\partial_y, \qquad \partial_x=e^{(2\ell+1)s/2}\partial_y.
\end{equation*}
Using these relations in \eqref{eq:b-slope}, we obtain
\begin{equation*}
\partial_sW_y+W_y+\left(\frac{2\ell+1}{2}y+\frac{W+e^{(2\ell-1)s/2}(\kappa-\dot\xi)}{1-\dot\tau}\right)W_{yy}+\frac{W_y^2}{\ell(1-\dot\tau)}=\frac{G}{1-\dot\tau},
\end{equation*}
where
\begin{equation}\label{eq:b-G}
G(y,s):=e^{-2s}F\left(\xi(t)+e^{-(2\ell+1)s/2}y,t\right).
\end{equation}
Thus $W_y$ satisfies \eqref{eq:perturbed-slope-equation} with
\begin{equation}\label{eq:b-AR}
\mathcal V=\frac{2\ell+1}{2}y+W+\mathcal A, \quad \mathcal A=\frac{\dot\tau W+e^{(2\ell-1)s/2}(\kappa-\dot\xi)}{1-\dot\tau},\quad \mathcal R=\frac{G-\ell^{-1}\dot\tau W_y^2}{1-\dot\tau}.
\end{equation}

The modulation equations \eqref{eq:b-modulation-tau} preserve the normalization \eqref{eq:modulation-constraints-intro}. Indeed, evaluating \eqref{u-nonlocal}, \eqref{eq:b-slope}, and the $x$-derivative of \eqref{eq:b-slope} along $x=\xi(t)$, we obtain
\begin{equation*}
\begin{aligned}
\frac{d}{dt}w(\xi(t),t)&=(\dot\xi-w(\xi(t),t))w_x(\xi(t),t)-P_x(\xi(t),t),\\
\frac{d}{dt}w_x(\xi(t),t)&=(\dot\xi-w(\xi(t),t))w_{xx}(\xi(t),t)-\frac1\ell w_x(\xi(t),t)^2+F(\xi(t),t),\\
\frac{d}{dt}w_{xx}(\xi(t),t)&=(\dot\xi-w(\xi(t),t))w_{xxx}(\xi(t),t) \\
& \quad - \left(1+\frac2\ell\right)w_x(\xi(t),t)w_{xx}(\xi(t),t) + F_x(\xi(t),t).
\end{aligned}
\end{equation*}
By \eqref{eq:b-modulation-tau}, the triple $(\kappa(t),-\ell/(\tau(t)-t),0)$ solves this ODE system. Since the initial values agree by \eqref{eq:b-initial-matching} and the initial data for the modulation functions, uniqueness implies
\begin{equation*}
w(\xi(t),t)=\kappa(t),\qquad w_x(\xi(t),t)=-\frac{\ell}{\tau(t)-t},\qquad w_{xx}(\xi(t),t)=0.
\end{equation*}
Consequently, by the change of variables \eqref{eq:modulated-variables}--\eqref{eq:modulated-unknown}, $W$ satisfies the normalization \eqref{eq:modulation-constraints-intro} on $[s_0,s(T)]$.

It remains to establish the regularity and continuity properties in assumption \emph{(i)}. Let $[s_0,s_1]$ be any finite interval contained in the self-similar time interval of existence. Lemma~\ref{lem:b-local}, the smoothness of the modulation functions, and \eqref{eq:b-AR} imply that $W,\mathcal V,\mathcal R$ are smooth on $\mathbb R\times[s_0,s_1]$. By \eqref{w_reg}, Sobolev embedding, and the change of variables \eqref{eq:modulated-variables}--\eqref{eq:modulated-unknown}, we have
\begin{equation} \label{W_reg}
\partial_y^jW\in C([s_0,s_1];L^\infty(\mathbb R)),\qquad j=0,1,2,3,4.
\end{equation}
In particular, $W$ is bounded on this interval, and \eqref{eq:b-AR} gives
\begin{equation*}
|\mathcal V(y,s)|\leq C(1+|y|),\qquad (y,s)\in\mathbb R\times[s_0,s_1],
\end{equation*}
which is \eqref{eq:common-singularity-velocity}.

Together with \eqref{eq:modulated-variables}--\eqref{eq:modulated-unknown} and Sobolev embedding, \eqref{weighted_reg} implies
\begin{equation*}
|W_y(y,s)|\leq\frac{C}{(1+y^2)^{1/2}},\qquad (y,s)\in\mathbb R\times[s_0,s_1].
\end{equation*}
Since $2/(2\ell+1)<1$, the function $(1+|y|^{2/(2\ell+1)})W_y(y,s)$ tends to zero as $|y|\to\infty$, uniformly for $s\in[s_0,s_1]$.  Combining this uniform decay with \eqref{W_reg}, we obtain
\begin{equation*}
(1+|y|^{2/(2\ell+1)})W_y\in C([s_0,s_1];L^\infty(\mathbb R)).
\end{equation*}
Recalling from \eqref{eq:far-Uprime} that $(1+|y|^{2/(2\ell+1)})U'\in L^\infty(\mathbb R)$, we have
\begin{equation}\label{eq:b-weighted-slope-continuity}
(1+|y|^{2/(2\ell+1)})(W_y-U')\in C([s_0,s_1];L^\infty(\mathbb R)).
\end{equation}

By \eqref{eq:modulation-constraints-intro}, \eqref{eq:origin-derivatives}, and Taylor's formula, we have
\begin{equation*}
\begin{aligned}
& \frac{1+y^2}{y^2}(W_y(y,s)-U'(y)) =(1+y^2)\int_0^1(1-\theta)\bigl((\partial_y^3W)(\theta y,s)-U^{(3)}(\theta y)\bigr)\,d\theta,\\
& \frac{\sqrt{1+y^2}}{y}W_{yy}(y,s) =\sqrt{1+y^2}\int_0^1(\partial_y^3W)(\theta y,s)\,d\theta
\end{aligned}
\end{equation*}
for $y\neq0$. The right-hand sides extend continuously to $y=0$ and are continuous in $s$ with values in $L^\infty(\{|y|\leq1\})$. On $\{|y|\geq1\}$, the factors $(1+y^2)/y^2$ and $\sqrt{1+y^2}/y$ are bounded, so the required $L^\infty$ continuity follows from \eqref{W_reg}. Combining these observations with \eqref{W_reg} and \eqref{eq:b-weighted-slope-continuity}, we conclude that all the weighted quantities in assumption~\emph{(i)} of Theorem~\ref{thm:common-singularity} belong to $C([s_0,s_1];L^\infty(\mathbb R))$. This verifies assumption~\emph{(i)}.

\subsubsection{Initial bounds}\label{subsubsec:b-initial-bounds}

We verify assumption~\emph{(ii)} of Theorem~\ref{thm:common-singularity}.

By the change of variables \eqref{eq:modulated-variables}--\eqref{eq:modulated-unknown}, we have
\begin{equation}\label{eq:b-initial-profile}
W(y,s_0)=\varepsilon^{-(2\ell-1)/2}\left(w_0\left(\varepsilon^{(2\ell+1)/2}y\right)-w_0(0)\right),\qquad s_0=-\log\varepsilon.
\end{equation}
In particular,
\begin{equation}\label{eq:b-initial-rescaled-derivatives}
\partial_y^jW(y,s_0) = \varepsilon^{\left(j(2\ell+1)-(2\ell-1)\right)/2}w_0^{(j)}\left(\varepsilon^{(2\ell+1)/2}y\right),\qquad j=1,2,3,4.
\end{equation}
Using \eqref{eq:b-initial-slope} and \eqref{eq:b-initial-rescaled-derivatives}, we obtain
\begin{equation*}
\begin{split}
|W_y(y,s_0)-U'(y)|&=\left|\varepsilon w_0'\left(\varepsilon^{(2\ell+1)/2}y\right)-U'(y)\right|\\
&\leq(1-2\vartheta)\min\left\{\frac{Ny^2}{1+y^2},\frac{B_\ell}{1+|y|^{2/(2\ell+1)}}\right\}.
\end{split}
\end{equation*}
Similarly, \eqref{eq:b-initial-origin} and \eqref{eq:b-initial-derivatives}, together with \eqref{eq:b-initial-rescaled-derivatives}, give \eqref{eq:bootstrap-initial-second}--\eqref{eq:bootstrap-initial-fourth}. Hence the initial data for $W$ satisfy \eqref{eq:bootstrap-initial-bounds}. This proves assumption~\emph{(ii)} of Theorem~\ref{thm:common-singularity}.

Since $\tau(0)=\varepsilon$ and $\xi(0)=0$, \eqref{eq:b-initial-energy} also gives \eqref{eq:b-auxiliary-bootstrap} at $t=0$ with strict inequalities. Combining this with \eqref{eq:bootstrap-initial-bounds}, \eqref{weighted_reg}, and the continuity established in Section~\ref{subsubsec:b-modulated}, we obtain some $s_1\in(s_0,s(T)]$ such that \eqref{eq:common-bootstrap-assumptions} holds for all $s\in[s_0,s_1]$ and \eqref{eq:b-auxiliary-bootstrap} holds for all $t\in[0,t(s_1)]$.

\subsubsection{Modulation estimates}\label{subsubsec:b-modulation}

We next verify assumption~\emph{(iv)} of Theorem~\ref{thm:common-singularity} before turning to assumption~\emph{(iii)}, since the remainder bounds require estimates for the modulation terms.

Let $[s_0,s_1]$ be a finite interval on which \eqref{eq:common-bootstrap-assumptions} holds. By Section~\ref{subsubsec:b-initial-bounds}, \eqref{eq:b-auxiliary-bootstrap} holds on $[0,t(\sigma_1)]$ for some $\sigma_1\in(s_0,s_1]$. We first prove \eqref{eq:common-singularity-modulation} on any such interval. We then improve \eqref{eq:b-auxiliary-bootstrap} and extend the estimates to $[0,t(s_1)]$ by a continuity argument.

\smallskip
\noindent$\bullet$\ \textit{Estimate of $1-\dot\tau$.}
Applying \eqref{eq:b-weighted-pressure} to \eqref{eq:b-modulation-tau} and using \eqref{eq:b-pressure-rate-limit} and $e^{-s}=\tau(t)-t\leq2\varepsilon$, we obtain
\begin{equation}\label{eq:b-tau-derivative-bound}
|\dot\tau(t)|=\frac{e^{-2s}}{\ell}|F(\xi(t),t)|\leq Ce^{-2s}\mathcal H_\ell(e^{-s})\leq\frac14
\end{equation}
for sufficiently small $\varepsilon>0$. Hence
\begin{equation}\label{eq:b-tau-factor-bound}
\frac12\leq1-\dot\tau(t)\leq\frac32.
\end{equation}
Using \eqref{eq:modulated-variables} and $t(s_0)=0$, we also obtain
\begin{equation}\label{eq:b-time-window}
\frac{dt}{ds}=\frac{e^{-s}}{1-\dot\tau}\leq2e^{-s},\qquad t(s)\leq2\int_{s_0}^s e^{-\sigma}\,d\sigma<2\varepsilon.
\end{equation}

\smallskip
\noindent$\bullet$\ \textit{Estimate of $\kappa$.}
By \eqref{eq:modulated-unknown}, \eqref{eq:modulation-constraints-intro}, and \eqref{eq:b-auxiliary-bootstrap}, we have
\begin{equation}\label{eq:b-kappa-bound}
|\kappa(t)|=|w(\xi(t),t)|\leq\|w(\cdot,t)\|_{L^\infty}\leq\|(1+x^2)^{1/(2\ell+1)}w(\cdot,t)\|_{L^\infty}\leq1.
\end{equation}

\smallskip
\noindent$\bullet$\ \textit{Estimate of $\xi$.}
From \eqref{eq:b-modulation-tau} and the identity $\partial_x^3w(\xi(t),t)=e^{(2\ell+2)s}\partial_y^3W(0,s)$, we have
\begin{equation*}
|\dot\xi(t)|\leq|\kappa(t)|+\frac{e^{-(2\ell+2)s}|F_x(\xi(t),t)|}{|\partial_y^3W(0,s)|}.
\end{equation*}
By \eqref{eq:b-kappa-bound}, the first term is bounded by $1$. To estimate the second term, we use \eqref{eq:boot-origin}, which gives
\begin{equation}\label{eq:b-center-third}
\partial_y^3W(0,s)\geq U^{(3)}(0)-1>0.
\end{equation}
Since $F_x=g'(w)w_x-P_x$ and $w_x(\xi(t),t)=-\ell e^s$, \eqref{eq:b-auxiliary-bootstrap} and \eqref{eq:b-weighted-pressure} give
\begin{equation}\label{eq:b-Fx-center-bound}
\begin{split}
|F_x(\xi(t),t)|&\leq\ell e^s\|g'\|_{L^\infty([-1,1])}+\|P_x(\cdot,t)\|_{L^\infty} \leq Ce^s+C\mathcal H_\ell(e^{-s})\leq Ce^s,
\end{split}
\end{equation}
where we used $1\leq\mathcal H_\ell(e^{-s})\leq e^s$ from \eqref{eq:b-pressure-rate}. Applying \eqref{eq:b-center-third} and \eqref{eq:b-Fx-center-bound} to \eqref{eq:b-modulation-tau}, we obtain
\begin{equation}\label{eq:b-xi-difference-bound}
|\dot\xi(t)-\kappa(t)|=\frac{e^{-(2\ell+2)s}|F_x(\xi(t),t)|}{\partial_y^3W(0,s)}\leq Ce^{-(2\ell+1)s}.
\end{equation}
Combining \eqref{eq:b-kappa-bound} and \eqref{eq:b-xi-difference-bound}, and using $s\geq s_0=-\log\varepsilon$, we obtain
\begin{equation}\label{eq:b-xi-speed-bound}
|\dot\xi(t)|\leq1+C\varepsilon^{2\ell+1}\leq2
\end{equation}
for sufficiently small $\varepsilon>0$. Integrating \eqref{eq:b-xi-speed-bound} and using \eqref{eq:b-time-window}, we obtain
\begin{equation}\label{eq:b-xi-integral-bound}
\int_0^t|\dot\xi(t')|\,dt'\leq2t<4\varepsilon.
\end{equation}
This proves \eqref{eq:common-singularity-modulation} for $0\leq t\leq t(\sigma_1)$.

\smallskip
\noindent$\bullet$\ \textit{Improvement of \eqref{eq:b-auxiliary-bootstrap} and continuation.}
Since $\tau(0)=\varepsilon$ and $\frac{d}{dt}(\tau(t)-t)=-(1-\dot\tau(t))$, \eqref{eq:b-tau-factor-bound} gives
\begin{equation}\label{eq:b-time-gap-improved}
0<\tau(t)-t\leq\varepsilon.
\end{equation}
Since $\xi(0)=0$, \eqref{eq:b-xi-integral-bound} gives
\begin{equation}\label{eq:b-xi-improved}
|\xi(t)|\leq\int_0^t|\dot\xi(t')|\,dt'<4\varepsilon<\frac12
\end{equation}
for sufficiently small $\varepsilon>0$.

To improve the bound on $w$ in \eqref{eq:b-auxiliary-bootstrap}, we apply \eqref{eq:b-weighted-amplitude}. By \eqref{eq:b-time-window}, the integral on the right-hand side of \eqref{eq:b-weighted-amplitude} satisfies
\begin{equation}\label{eq:b-physical-forcing-integral}
\begin{split}
\int_0^{t(s)}\mathcal H_\ell(\tau(t')-t')\,dt'&=\int_{s_0}^s\frac{e^{-\sigma}\mathcal H_\ell(e^{-\sigma})}{1-\dot\tau(t(\sigma))}\,d\sigma \leq2\int_{e^{-s}}^\varepsilon\mathcal H_\ell(r)\,dr\leq2\int_0^\varepsilon\mathcal H_\ell(r)\,dr.
\end{split}
\end{equation}
Applying \eqref{eq:b-initial-energy}, \eqref{eq:b-time-window}, and \eqref{eq:b-physical-forcing-integral} to \eqref{eq:b-weighted-amplitude}, we obtain
\begin{equation}\label{eq:b-low-bounds}
\|(1+x^2)^{1/(2\ell+1)}w(\cdot,t)\|_{L^\infty}\leq e^{2C\varepsilon}\left(\frac14+2C\int_0^\varepsilon\mathcal H_\ell(r)\,dr\right)\leq\frac12
\end{equation}
for sufficiently small $\varepsilon>0$. Here we used \eqref{eq:b-pressure-rate-smallness}.

The bounds \eqref{eq:b-time-gap-improved}, \eqref{eq:b-xi-improved}, and \eqref{eq:b-low-bounds} strictly improve \eqref{eq:b-auxiliary-bootstrap} on $[0,t(\sigma_1)]$. The constants and the choice of $\varepsilon$ are independent of $\sigma_1$ and $s_1$. The quantities in \eqref{eq:b-auxiliary-bootstrap} are continuous in time by \eqref{weighted_reg}, Sobolev embedding, and the smoothness of the modulation functions. A continuity argument therefore yields \eqref{eq:b-auxiliary-bootstrap} on $[0,t(s_1)]$. By \eqref{eq:b-tau-factor-bound}, \eqref{eq:b-kappa-bound}, and \eqref{eq:b-xi-integral-bound}, we have
\begin{equation*}
\frac12\leq1-\dot\tau(t)\leq\frac32,\qquad |\kappa(t)|\leq1,\qquad \int_0^t|\dot\xi(t')|\,dt'\leq4\varepsilon
\end{equation*}
for all $t\in[0,t(s_1)]$. This verifies assumption~\emph{(iv)} of Theorem~\ref{thm:common-singularity}.

\smallskip
\noindent$\bullet$\ \textit{Further estimates.}
We also record estimates for the modulation functions that will be used in Section~\ref{subsubsec:b-remainder-bounds}. Applying \eqref{eq:b-weighted-pressure} and \eqref{eq:b-xi-difference-bound} to \eqref{eq:b-modulation-tau}, we obtain
\begin{equation}\label{eq:b-kappa-derivative-bound}
\begin{split}
|\dot\kappa(t)|&\leq|P_x(\xi(t),t)|+\ell e^s|\dot\xi(t)-\kappa(t)| \leq C\mathcal H_\ell(e^{-s})+Ce^{-2\ell s}\leq C\mathcal H_\ell(e^{-s}).
\end{split}
\end{equation}
Combining \eqref{eq:b-tau-derivative-bound}, \eqref{eq:b-xi-difference-bound}, and \eqref{eq:b-kappa-derivative-bound}, we obtain
\begin{equation}\label{eq:b-modulation-bounds}
\begin{aligned}
|\dot\tau|&\leq Ce^{-2s}\mathcal H_\ell(e^{-s}),&\qquad |\dot\kappa|&\leq C\mathcal H_\ell(e^{-s}),\\
|\dot\xi-\kappa|&\leq Ce^{-(2\ell+1)s},& |e^{(2\ell-1)s/2}(\kappa-\dot\xi)|&\leq Ce^{-(2\ell+3)s/2}
\end{aligned}
\end{equation}
for all $s\in[s_0,s_1]$, with constants independent of $\varepsilon$ and $s_1$.

\subsubsection{Bounds on the model-specific remainders}\label{subsubsec:b-remainder-bounds}

Let $[s_0,s_1]$ be a finite interval contained in the self-similar time interval of existence on which \eqref{eq:common-bootstrap-assumptions} holds. By Remark~\ref{rem:sufficient-remainder-bounds}, it suffices to verify \eqref{eq:sufficient-remainder-bounds} and \eqref{eq:sufficient-remainder-smallness} to establish assumption~\emph{(iii)} of Theorem~\ref{thm:common-singularity}. The functions $\mathcal A$ and $\mathcal R$ are given by \eqref{eq:b-AR}. By Section~\ref{subsubsec:b-modulation}, \eqref{eq:b-auxiliary-bootstrap} holds on $[0,t(s_1)]$, so we may apply \eqref{eq:b-weighted-pressure}. We first estimate the nonlocal term $G$ in \eqref{eq:b-G} and its first three $y$-derivatives, and then combine these estimates with \eqref{eq:b-modulation-bounds} to verify the sufficient condition with the choice $e(s)=Ce^{-s}\mathcal H_\ell(e^{-s})$.

\smallskip
\noindent$\bullet$\ \textit{Estimates for $G$.}
By \eqref{eq:b-G} and \eqref{eq:b-weighted-pressure}, we have
\begin{equation}\label{eq:b-G-zero-bound}
\|G(\cdot,s)\|_{L^\infty}=e^{-2s}\|F(\cdot,t)\|_{L^\infty}\leq Ce^{-2s}\mathcal H_\ell(e^{-s}).
\end{equation}
To estimate the derivatives of $G$, we use \eqref{P-nonlocal} and \eqref{eq:b-slope}, which give
\begin{equation*}
P_{xx}=P-g(w)-\frac{3-b}{2}w_x^2=-F-\frac{3-b}{2}w_x^2.
\end{equation*}
Differentiating $F=g(w)-P$, we obtain
\begin{equation}\label{eq:b-F-derivatives-one}
F_x=g'(w)w_x-P_x,\qquad F_{xx}=F+\left(g''(w)+\frac{3-b}{2}\right)w_x^2+g'(w)w_{xx},
\end{equation}
and
\begin{equation}\label{eq:b-F-derivatives-two}
\partial_x^3F=F_x+g^{(3)}(w)w_x^3+\bigl(3g''(w)+3-b\bigr)w_xw_{xx}+g'(w)\partial_x^3w.
\end{equation}
Using \eqref{eq:b-F-derivatives-one}--\eqref{eq:b-F-derivatives-two} in \eqref{eq:b-G}, together with the change of variables \eqref{eq:modulated-variables}--\eqref{eq:modulated-unknown}, we obtain
\begin{align}
G_y &=e^{-(2\ell+3)s/2}g'(w)W_y-e^{-(2\ell+5)s/2}P_x, \label{eq:b-Gy}\\
G_{yy} &=e^{-(2\ell+1)s}G+e^{-(2\ell+1)s}\left(g''(w)+\frac{3-b}{2}\right)W_y^2+e^{-(2\ell+3)s/2}g'(w)W_{yy}, \label{eq:b-Gyy}
\end{align}
and
\begin{equation}\label{eq:b-Gyyy}
\begin{aligned}
\partial_y^3G&=e^{-(2\ell+1)s}G_y+e^{-(6\ell+1)s/2}g^{(3)}(w)W_y^3\\
&\quad+e^{-(2\ell+1)s}\bigl(3g''(w)+3-b\bigr)W_yW_{yy}+e^{-(2\ell+3)s/2}g'(w)\partial_y^3W.
\end{aligned}
\end{equation}
Here $w$ and $P_x$ in \eqref{eq:b-Gy}--\eqref{eq:b-Gyyy} are evaluated at $(\xi(t)+e^{-(2\ell+1)s/2}y,t)$, where $t=t(s)$.

Since $\|w(\cdot,t)\|_{L^\infty}\leq1/2$ by \eqref{eq:b-low-bounds}, the functions $g'(w)$, $g''(w)$, and $g^{(3)}(w)$ are uniformly bounded. Applying \eqref{eq:common-bootstrap-assumptions}, \eqref{eq:Uprime-norm}, \eqref{eq:b-weighted-pressure}, and \eqref{eq:b-G-zero-bound} to \eqref{eq:b-Gy}--\eqref{eq:b-Gyyy}, we obtain
\begin{equation}\label{eq:b-G-derivative-estimates}
\begin{aligned}
\|G_y\|_{L^\infty}&\leq Ce^{-(2\ell+3)s/2}+Ce^{-(2\ell+5)s/2}\mathcal H_\ell(e^{-s})\leq Ce^{-(2\ell+3)s/2},\\
\|G_{yy}\|_{L^\infty}&\leq e^{-(2\ell+1)s}\|G\|_{L^\infty}+Ce^{-(2\ell+1)s}+Ce^{-(2\ell+3)s/2}\leq Ce^{-(2\ell+3)s/2},\\
\|\partial_y^3G\|_{L^\infty}&\leq e^{-(2\ell+1)s}\|G_y\|_{L^\infty}+Ce^{-(6\ell+1)s/2}+Ce^{-(2\ell+1)s}+Ce^{-(2\ell+3)s/2}\leq Ce^{-(2\ell+3)s/2}.
\end{aligned}
\end{equation}
Here we used $1\leq\mathcal H_\ell(e^{-s})\leq e^s$ from \eqref{eq:b-pressure-rate} and
\begin{equation*}
2\ell+1>\frac{2\ell+3}{2},\qquad \frac{6\ell+1}{2}>\frac{2\ell+3}{2},\qquad \ell>\frac12.
\end{equation*}
Combining \eqref{eq:b-G-zero-bound} and \eqref{eq:b-G-derivative-estimates}, we have
\begin{equation}\label{eq:b-G-bounds}
\|G\|_{L^\infty}\leq Ce^{-2s}\mathcal H_\ell(e^{-s}),\qquad \sum_{j=1}^3\|\partial_y^jG\|_{L^\infty}\leq Ce^{-(2\ell+3)s/2}.
\end{equation}

For the weighted estimate of $G$, the change of variables \eqref{eq:modulated-variables} and \eqref{eq:b-xi-improved} give
\begin{equation}\label{eq:b-weight-conversion}
1+|y|^{2/(2\ell+1)}=1+e^s|x-\xi(t)|^{2/(2\ell+1)}\leq Ce^s(1+x^2)^{1/(2\ell+1)}.
\end{equation}
Applying \eqref{eq:b-weight-conversion} and \eqref{eq:b-weighted-pressure} to \eqref{eq:b-G}, we obtain
\begin{equation}\label{eq:b-G-weighted}
\|(1+|y|^{2/(2\ell+1)})G\|_{L^\infty}\leq Ce^{-s}\|(1+x^2)^{1/(2\ell+1)}F\|_{L^\infty}\leq Ce^{-s}\mathcal H_\ell(e^{-s}).
\end{equation}

\smallskip
\noindent$\bullet$\ \textit{Estimates for $\mathcal A$ and $\mathcal R$.}
Using \eqref{eq:boot-far}, \eqref{eq:far-Uprime}, and Lemma~\ref{lem:derivative-bounds}, together with $W(0,s)=0$ from \eqref{eq:modulation-constraints-intro}, we obtain
\begin{equation}\label{eq:b-W-growth}
|W_y(y,s)|\leq\frac{C}{1+|y|^{2/(2\ell+1)}},\qquad |W(y,s)|\leq C|y|.
\end{equation}
By \eqref{eq:b-AR}, \eqref{eq:b-tau-factor-bound}, and \eqref{eq:b-modulation-bounds},
\begin{equation}\label{eq:b-A-pointwise}
|\mathcal A(y,s)|\leq Ce^{-2s}\mathcal H_\ell(e^{-s})|W(y,s)|+Ce^{-(2\ell+3)s/2}.
\end{equation}
It follows from \eqref{eq:b-W-growth} and \eqref{eq:b-A-pointwise} that
\begin{equation}\label{eq:b-A-sufficient}
\left\|\frac{\mathcal A}{1+|y|}\right\|_{L^\infty}\leq C\left(e^{-2s}\mathcal H_\ell(e^{-s})+e^{-(2\ell+3)s/2}\right)\leq Ce^{-2s}\mathcal H_\ell(e^{-s}),
\end{equation}
where we used $(2\ell+3)/2>2$ and $\mathcal H_\ell(e^{-s})\geq1$.

For $\mathcal R$, \eqref{eq:b-W-growth} gives
\begin{equation*}
\|(1+|y|^{2/(2\ell+1)})W_y^2\|_{L^\infty}\leq C.
\end{equation*}
Using this bound together with \eqref{eq:b-G-weighted}, \eqref{eq:b-tau-factor-bound}, and \eqref{eq:b-modulation-bounds} in \eqref{eq:b-AR}, we obtain
\begin{equation}\label{eq:b-R-sufficient}
\begin{split}
\|(1+|y|^{2/(2\ell+1)})\mathcal R\|_{L^\infty}
&\leq2\|(1+|y|^{2/(2\ell+1)})G\|_{L^\infty}+\frac{2}{\ell}|\dot\tau|\|(1+|y|^{2/(2\ell+1)})W_y^2\|_{L^\infty}\\
&\leq Ce^{-s}\mathcal H_\ell(e^{-s}).
\end{split}
\end{equation}

For the differentiated remainders, we use \eqref{eq:differentiated-remainders}. Since $\dot\tau$, $\kappa$, and $\dot\xi$ depend only on $s$, \eqref{eq:b-AR} gives
\begin{equation}\label{eq:b-A-derivative-identity}
\partial_y^j\mathcal A=\frac{\dot\tau}{1-\dot\tau}\partial_y^jW,\qquad j=1,2,3.
\end{equation}
Substituting \eqref{eq:b-A-derivative-identity} and \eqref{eq:b-AR} into \eqref{eq:differentiated-remainders}, we obtain
\begin{align}
\mathcal R_2 &=\frac{G_y-\dot\tau(1+2/\ell)W_yW_{yy}}{1-\dot\tau}, \label{eq:b-R2}\\
\mathcal R_3 &=\frac{G_{yy}-\dot\tau\bigl((1+2/\ell)W_{yy}^2+(2+2/\ell)W_y\partial_y^3W\bigr)}{1-\dot\tau}, \label{eq:b-R3}
\end{align}
and
\begin{equation}\label{eq:b-R4}
\mathcal R_4=\frac{\partial_y^3G-\dot\tau\bigl((4+6/\ell)W_{yy}\partial_y^3W+(3+2/\ell)W_y\partial_y^4W\bigr)}{1-\dot\tau}.
\end{equation}
Taking the $L^\infty$-norm in \eqref{eq:b-R2}--\eqref{eq:b-R4} and using \eqref{eq:common-bootstrap-assumptions}, \eqref{eq:b-tau-factor-bound}, \eqref{eq:b-G-bounds}, and \eqref{eq:b-modulation-bounds}, we obtain
\begin{equation}\label{eq:b-R-derivative-bounds}
\sum_{j=2}^4\|\mathcal R_j\|_{L^\infty}\leq C\left(\sum_{j=1}^3\|\partial_y^jG\|_{L^\infty}+|\dot\tau|\right)\leq Ce^{-2s}\mathcal H_\ell(e^{-s}).
\end{equation}

\smallskip
\noindent$\bullet$\ \textit{Verification of the sufficient condition.}
Combining \eqref{eq:b-A-sufficient}, \eqref{eq:b-R-sufficient}, and \eqref{eq:b-R-derivative-bounds}, we obtain \eqref{eq:sufficient-remainder-bounds} with
\begin{equation*}
e(s):=Ce^{-s}\mathcal H_\ell(e^{-s}),
\end{equation*}
where $C>0$ is sufficiently large and independent of $\varepsilon$ and $s_1$. Since $s_0=-\log\varepsilon$, \eqref{eq:b-pressure-rate-smallness} gives
\begin{equation*}
\int_{s_0}^{s_1}e(s)\,ds\leq C\int_{s_0}^{\infty}e^{-s}\mathcal H_\ell(e^{-s})\,ds=C\int_0^\varepsilon\mathcal H_\ell(r)\,dr\longrightarrow0\qquad\text{as }\varepsilon\downarrow0.
\end{equation*}
Thus, for sufficiently small $\varepsilon>0$, independently of $s_1$, the bound \eqref{eq:sufficient-remainder-smallness} holds. By Remark~\ref{rem:sufficient-remainder-bounds}, this verifies assumption~\emph{(iii)} of Theorem~\ref{thm:common-singularity}.

\subsubsection{Continuation and completion of the proof}\label{subsubsec:b-continuation}

It remains to verify assumption~\emph{(v)} of Theorem~\ref{thm:common-singularity}. Suppose that $s_0<S<\infty$, that the solution and the modulation functions are defined for all $s \in [s_0,S)$ with the properties in \emph{(i)}, and that \eqref{eq:common-bootstrap-assumptions} holds on this interval. Since $t(s)$ is strictly increasing and satisfies \eqref{eq:b-time-window}, the limit $T_S:=\lim_{s\uparrow S}t(s)$ exists and is finite. For any $s_1\in(s_0,S)$, Section~\ref{subsubsec:b-modulation} ensures that \eqref{eq:b-auxiliary-bootstrap} holds on $[0,t(s_1)]$. Thus \eqref{eq:b-physical-slope} gives
\begin{equation}\label{eq:b-continuation-slope}
\|w_x(\cdot,t(s))\|_{L^\infty}\leq\frac{C}{\tau(t(s))-t(s)}=Ce^s\leq Ce^S
\end{equation}
for all $s\in[s_0,S)$. The bound \eqref{eq:b-continuation-slope} and the blow-up criterion \eqref{eq:b-local-continuation} imply that $T_S<T_{\max}$. Hence $w$ extends smoothly beyond $T_S$, with the higher Sobolev regularity and the weighted regularity preserved by Lemma~\ref{lem:b-local}.

We next continue the modulation functions using \eqref{eq:b-modulation-tau}. Since $e^{-s}\geq e^{-S}>0$, the bounds \eqref{eq:b-modulation-bounds} and \eqref{eq:b-xi-speed-bound} imply that $\dot\tau$, $\dot\xi$, and $\dot\kappa$ are bounded on $[0,T_S)$. Thus $\tau$, $\xi$, and $\kappa$ have finite limits as $t\uparrow T_S$. Extending these functions continuously to $T_S$, we obtain from \eqref{eq:modulated-variables}, \eqref{eq:modulated-unknown}, and \eqref{eq:b-center-third} that
\begin{equation}\label{eq:b-continuation-nondegeneracy}
\tau(T_S)-T_S=e^{-S}>0,\qquad \partial_x^3w(\xi(T_S),T_S)\geq e^{(2\ell+2)S}\bigl(U^{(3)}(0)-1\bigr)>0.
\end{equation}
By \eqref{eq:b-continuation-nondegeneracy} and the local well-posedness theory for ODEs, the solution $(\tau,\xi,\kappa)$ of \eqref{eq:b-modulation-tau} extends beyond $T_S$. Moreover, \eqref{eq:b-tau-factor-bound} and continuity give $1-\dot\tau(t)>0$ for $t$ sufficiently close to $T_S$. Hence the change of variables \eqref{eq:modulated-variables}--\eqref{eq:modulated-unknown} extends beyond $s=S$. The argument in Section~\ref{subsubsec:b-modulated} applies on the extended interval and preserves the self-similar formulation, the normalization \eqref{eq:modulation-constraints-intro}, and the regularity and continuity properties in \emph{(i)}. This verifies assumption~\emph{(v)}.

All hypotheses of Theorem~\ref{thm:common-singularity} are now satisfied with $t_0=0$ and $\ell=2/(b-1)$. Applying that theorem, we conclude that $W$ is defined for all $s\geq s_0=-\log\varepsilon$ and satisfies \eqref{eq:bootstrap-improved-bounds}. In the original variables, there exist $T_*>0$ and $x_*\in\mathbb R$ such that $w$ remains smooth on $[0,T_*)$ and
\begin{equation*}
\tau(t)-t\sim T_*-t,\qquad \xi(t)\longrightarrow x_*\quad\text{as }t\uparrow T_*.
\end{equation*}
Letting $s\to\infty$ in \eqref{eq:b-time-window}, we obtain $T_*\leq2\varepsilon$. By \eqref{eq:b-low-bounds}, we also have
\begin{equation*}
\sup_{0\leq t<T_*}\|w(\cdot,t)\|_{L^\infty(\mathbb R)}\leq\frac12.
\end{equation*}
Since $\gamma_\ell=(2\ell-1)/(2\ell+1)=(5-b)/(b+3)=\gamma_b$, the conclusions of Theorem~\ref{thm:common-singularity}, together with \eqref{eq:b-low-bounds}, give assertions \emph{(i)}--\emph{(iii)} of Theorem~\ref{thm:b-cusp-intro}. This completes the proof of Theorem~\ref{thm:b-cusp}.
\qed

\appendix

\section{Admissible initial data}\label{subsec:b-data}

We construct a non-empty relatively open class of smooth initial data satisfying \eqref{eq:b-initial-matching}--\eqref{eq:b-initial-energy} for sufficiently small \(\varepsilon>0\). The relative topology is induced by the norm \(\|(1+x^2)^{1/2}w_0\|_{H^5}\) for smooth data satisfying the regularity assumptions in Section~\ref{subsec:b-theorem} and the matching conditions \eqref{eq:b-initial-matching}. The construction follows \cite[Appendix~A]{KKS-rSV} and \cite[Section~5.2]{KKKS-Novikov}.

\subsection{Construction of an initial datum} \label{app:A1}

To construct an admissible initial datum, it suffices to construct an even function \(v\in C_c^\infty(\mathbb R)\) such that
\begin{equation}\label{eq:b-g-admissible}
v(y)=U'(y)\ \text{near }y=0,\qquad \int_{\mathbb R}v(y)\,dy=0,\qquad |v(y)-U'(y)|\leq\rho b_*(y),
\end{equation}
for some \(\rho<1\), with derivative bounds corresponding to \eqref{eq:b-initial-derivatives}, where
\begin{equation} \label{bstar}
b_*(y):=(1-2\vartheta)\min\left\{\frac{Ny^2}{1+y^2},\frac{B_\ell}{1+|y|^{2/(2\ell+1)}}\right\}.
\end{equation}
Once such a function \(v\) has been constructed, we set
\begin{equation}\label{eq:b-Phi}
\Phi(y):=\int_{-\infty}^y v(y')\,dy',\qquad w_0^\varepsilon(x):=\varepsilon^{(2\ell-1)/2}\Phi\left(\frac{x}{\varepsilon^{(2\ell+1)/2}}\right).
\end{equation}
Then, since \(v\) is compactly supported and has zero integral, \eqref{eq:b-Phi} gives \(\Phi\in C_c^\infty(\mathbb R)\), and hence \(w_0^\varepsilon\in C_c^\infty(\mathbb R)\). Moreover,
\begin{equation}\label{eq:b-data-rescaled-slope}
\varepsilon(w_0^\varepsilon)'\left(\varepsilon^{(2\ell+1)/2}y\right)=v(y).
\end{equation}
Thus the third condition in \eqref{eq:b-g-admissible} and the fact that \(\rho<1\) give \eqref{eq:b-initial-slope}. We now construct such a function \(v\) and verify \eqref{eq:b-initial-matching}--\eqref{eq:b-initial-derivatives} and \eqref{eq:b-initial-energy} for \(w_0^\varepsilon\).

Let \(\chi\in C_c^\infty(\mathbb R)\) be an even cutoff function such that \(0\leq\chi\leq1\), \(\chi=1\) for \(|y|\leq1\), and \(\chi=0\) for \(|y|\geq2\). By \eqref{eq:common-singularity-theta}, we may choose \(\rho_0\) such that
\begin{equation*}
\frac{\Theta_\ell\beta^{-1/(2\ell+1)}}{(1-2\vartheta)B_\ell}<\rho_0<1.
\end{equation*}
By \eqref{eq:far-Uprime}, there exists \(R>1\) sufficiently large such that
\begin{equation*}
|U'(y)|\leq\rho_0b_*(y)\qquad\text{for }|y|\geq R.
\end{equation*}
Set \(\chi_R(y):=\chi(y/R)\) and \(v_1:=\chi_RU'\). Then \(v_1=U'\) for \(|y|\leq R\), \(v_1\in C_c^\infty(\mathbb R)\) is even, and
\begin{equation*}
|v_1(y)-U'(y)|\leq\rho_0b_*(y).
\end{equation*}
Moreover, since \(U'<0\) on \(\mathbb R\), we have
\begin{equation*}
\int_{\mathbb R}v_1(y)\,dy<0.
\end{equation*}

We now add a nonnegative correction to \(v_1\) in the far field so that the resulting function satisfies the zero-integral condition in \eqref{eq:b-g-admissible}. We choose an even nonnegative function \(\psi\in C_c^\infty(\mathbb R)\) supported in \(\{1<|y|<2\}\) such that \(\int_{\mathbb R}\psi(y)\,dy=1\) and set
\begin{equation*}
v_2(y):=-\frac{1}{L}\left(\int_{\mathbb R}v_1(y')\,dy'\right)\psi\left(\frac{y}{L}\right),\qquad L>2R.
\end{equation*}
Then \(v_2\geq0\) and
\begin{equation*}
\int_{\mathbb R}v_2(y)\,dy=-\int_{\mathbb R}v_1(y)\,dy.
\end{equation*}
By \eqref{bstar}, \(b_*(y)\sim(1-2\vartheta)B_\ell|y|^{-2/(2\ell+1)}\) as \(|y|\to\infty\). Since \(2/(2\ell+1)<1\), we have
\begin{equation*}
\sup_{L\leq|y|\leq2L}\frac{v_2(y)}{b_*(y)}\leq C_RL^{2/(2\ell+1)-1}\longrightarrow0\qquad\text{as }L\to\infty.
\end{equation*}
If we choose \(\rho\) to satisfy \(\rho_0<\rho<1\), then
\begin{equation*}
0\leq v_2(y)\leq(\rho-\rho_0)b_*(y)
\end{equation*}
for \(L\) sufficiently large. Setting \(v:=v_1+v_2\), we obtain an even function \(v\in C_c^\infty(\mathbb R)\) satisfying \eqref{eq:b-g-admissible}.

We next show that \(R\) and \(L\) can be chosen so that the data defined by \eqref{eq:b-Phi} also satisfy the derivative bounds \eqref{eq:b-initial-derivatives}. By Proposition~\ref{prop:profile-family} and the proof of Lemma~\ref{lem:derivative-bounds}, we have
\begin{equation*}
\lim_{|y|\to\infty}U^{(j)}(y)=0,\qquad j=1,2,3,4.
\end{equation*}
Hence, by the definition of \(\chi_R\),
\begin{equation} \label{eq:b-data-cutoff-limit}
\|v_1-U'\|_{W^{3,\infty}}=\|(\chi_R-1)U'\|_{W^{3,\infty}}\leq C_\chi\sum_{j=1}^4\sup_{|y|\geq R}|U^{(j)}(y)|\longrightarrow0\qquad\text{as }R\to\infty,
\end{equation}
where \(C_\chi>0\) is independent of \(R\). Moreover, the definition of \(v_2\) gives, for any fixed \(R>1\),
\begin{equation} \label{eq:b-data-correction-limit}
\|\partial_y^jv_2\|_{L^\infty}=\left|\int_{\mathbb R}v_1(y)\,dy\right|L^{-1-j}\|\psi^{(j)}\|_{L^\infty}\longrightarrow0\qquad\text{as }L\to\infty,\quad j=0,1,2,3.
\end{equation}
Using \(v=U'\) for \(|y|\leq R\), \(R>1\), and \(v=v_1+v_2\), we obtain
\begin{equation}\label{eq:b-data-weighted-approximation}
\begin{split}
\left\|\frac{\sqrt{1+y^2}}{y}v'\right\|_{L^\infty} & \leq \left\|\frac{\sqrt{1+y^2}}{y}U''\right\|_{L^\infty}  + \left\|\frac{\sqrt{1+y^2}}{y}(v'-U'')\right\|_{L^\infty} \\
& \leq\left\|\frac{\sqrt{1+y^2}}{y}U''\right\|_{L^\infty}+\sqrt2\bigl(\|v_1-U'\|_{W^{3,\infty}}+\|v_2\|_{W^{3,\infty}}\bigr)
\end{split}
\end{equation}
and
\begin{equation}\label{eq:b-data-higher-approximation}
\|v^{(j)}\|_{L^\infty}\leq\|U^{(j+1)}\|_{L^\infty}+\|v_1-U'\|_{W^{3,\infty}}+\|v_2\|_{W^{3,\infty}},\qquad j=2,3.
\end{equation}
Combining \eqref{eq:b-data-cutoff-limit} and \eqref{eq:b-data-correction-limit} with \eqref{eq:b-data-weighted-approximation}--\eqref{eq:b-data-higher-approximation} and the choice of \(M\) in Section~\ref{subsec:b-theorem}, we may first choose \(R\) sufficiently large and then \(L>2R\) sufficiently large so that
\begin{equation}\label{eq:b-data-second}
\left\|\frac{\sqrt{1+y^2}}{y}v'(y)\right\|_{L^\infty}<(1-2\vartheta)M^{1/8}
\end{equation}
and
\begin{equation}\label{eq:b-data-higher}
\|v''\|_{L^\infty}<(1-2\vartheta)M^{3/4},\qquad \|v^{(3)}\|_{L^\infty}<(1-2\vartheta)M.
\end{equation}

We verify \eqref{eq:b-initial-matching}--\eqref{eq:b-initial-derivatives} and \eqref{eq:b-initial-energy} for the data defined by \eqref{eq:b-Phi}. Since \(v\) is even and has zero integral, \(\Phi\) is odd. Since \(v=U'\) near the origin and \(\Phi(0)=U(0)=0\), we have \(\Phi=U\) near the origin. Thus \eqref{eq:origin-derivatives} gives
\begin{equation*}
(w_0^\varepsilon)'(0)=-\frac{\ell}{\varepsilon},\qquad (w_0^\varepsilon)''(0)=0,\qquad \varepsilon^{2\ell+2}(w_0^\varepsilon)^{(3)}(0)=U^{(3)}(0).
\end{equation*}
This proves \eqref{eq:b-initial-matching} and \eqref{eq:b-initial-origin}. By \eqref{eq:b-Phi}, for \(j=2,3,4\), we have
\begin{equation*}
(w_0^\varepsilon)^{(j)}(x)=\varepsilon^{[(2\ell-1)-j(2\ell+1)]/2}v^{(j-1)}\left(\frac{x}{\varepsilon^{(2\ell+1)/2}}\right).
\end{equation*}
Together with \eqref{eq:b-data-second}--\eqref{eq:b-data-higher}, these identities give \eqref{eq:b-initial-derivatives}.

To verify \eqref{eq:b-initial-energy}, we use \eqref{eq:b-Phi} and the fact that \(\Phi\) is compactly supported and independent of \(\varepsilon\) to obtain
\begin{equation*}
\|(1+x^2)^{1/(2\ell+1)}w_0^\varepsilon\|_{L^\infty}=\varepsilon^{(2\ell-1)/2}\|(1+\varepsilon^{2\ell+1}y^2)^{1/(2\ell+1)}\Phi(y)\|_{L_y^\infty}=O\left(\varepsilon^{(2\ell-1)/2}\right).
\end{equation*}
As \(\ell>1/2\), the right-hand side tends to zero as \(\varepsilon\downarrow0\). Hence \eqref{eq:b-initial-energy} holds with a strict inequality for all sufficiently small \(\varepsilon>0\). This completes the construction of an admissible initial datum.

\subsection{Relatively open class}

We show that there exists a non-empty relatively open set of admissible initial data containing \(w_0^\varepsilon\). Fix a sufficiently small \(\varepsilon>0\). We consider a perturbation \(f\) of \(w_0^\varepsilon\) such that \(w_0^\varepsilon+f\) satisfies the regularity assumptions in Section~\ref{subsec:b-theorem} and the matching conditions \eqref{eq:b-initial-matching}. Note that \(w_0^\varepsilon\) also satisfies \eqref{eq:b-initial-matching} and hence \(f'(0)=f''(0)=0\).

We first estimate \(f\) in terms of \(\|(1+x^2)^{1/2}f\|_{H^5}\). Since \((1+x^2)^{-1/2}\) and its derivatives are bounded, we have
\begin{equation}\label{eq:b11}
\|f\|_{H^5}\leq C\|(1+x^2)^{1/2}f\|_{H^5}.
\end{equation}
Using \eqref{eq:b11} and Sobolev embedding, we obtain
\begin{equation}\label{eq:b-data-weighted-first}
\begin{split}
\|(1+x^2)^{1/2}f'\|_{L^\infty} &\leq\left\|\partial_x\left((1+x^2)^{1/2}f\right)\right\|_{L^\infty}+\|f\|_{L^\infty}\\ 
&\leq C\left(\|(1+x^2)^{1/2}f\|_{H^2}+\|f\|_{H^1}\right)\leq C\|(1+x^2)^{1/2}f\|_{H^5}
\end{split}
\end{equation}
and
\begin{equation}\label{eq:b-data-perturbation-higher}
\|f''\|_{W^{2,\infty}} \leq C\|f\|_{H^5}\leq C\|(1+x^2)^{1/2}f\|_{H^5}.
\end{equation}
Moreover, since \(f'(0)=f''(0)=0\), Taylor's theorem yields
\begin{equation}\label{eq:b-data-taylor}
|f'(x)|\leq\frac{x^2}{2}\|f^{(3)}\|_{L^\infty},\qquad |f''(x)|\leq|x|\|f^{(3)}\|_{L^\infty}.
\end{equation}
Combining \eqref{eq:b-data-weighted-first}, \eqref{eq:b-data-perturbation-higher}, and \eqref{eq:b-data-taylor}, we obtain
\begin{equation}\label{eq:b-data-perturbation-bounds}
\begin{split}
|f'(x)| &\leq\min\left\{\frac{x^2}{2}\|f^{(3)}\|_{L^\infty},(1+x^2)^{-1/2}\|(1+x^2)^{1/2}f'\|_{L^\infty}\right\}\\ 
&\leq C\|(1+x^2)^{1/2}f\|_{H^5}\min\{x^2,(1+x^2)^{-1/2}\},\\ 
|f''(x)| &\leq\min\left\{|x|\|f^{(3)}\|_{L^\infty},\|f''\|_{L^\infty}\right\}\leq C\|(1+x^2)^{1/2}f\|_{H^5}\min\{|x|,1\}. 
\end{split}
\end{equation}

We compare \eqref{eq:b-data-perturbation-bounds} with the weights in \eqref{eq:b-initial-slope} and \eqref{eq:b-initial-second}. Since \(2/(2\ell+1)<1\), for \(x\neq0\) we have
\begin{equation*}
\begin{split} &\frac{\varepsilon^{2\ell+1}+x^2}{x^2}\min\{x^2,(1+x^2)^{-1/2}\}\leq\varepsilon^{2\ell+1}+1,\\ &\left(\varepsilon+|x|^{2/(2\ell+1)}\right)\min\{x^2,(1+x^2)^{-1/2}\}\leq\varepsilon+1,\\ &\frac{\sqrt{\varepsilon^{2\ell+1}+x^2}}{|x|}\min\{|x|,1\}\leq\sqrt{\varepsilon^{2\ell+1}+1}. \end{split}
\end{equation*}
Together with \eqref{eq:b-data-perturbation-bounds}, these inequalities yield
\begin{equation}\label{eq:b-data-perturbation-initial-weights}
\begin{split} &\varepsilon|f'(x)|\leq C_\varepsilon\|(1+x^2)^{1/2}f\|_{H^5}\min\left\{\frac{Nx^2}{\varepsilon^{2\ell+1}+x^2},\frac{\varepsilon B_\ell}{\varepsilon+|x|^{2/(2\ell+1)}}\right\},\\ &|f''(x)|\leq C_\varepsilon\|(1+x^2)^{1/2}f\|_{H^5}\frac{|x|}{\sqrt{\varepsilon^{2\ell+1}+x^2}}. \end{split}
\end{equation}

To verify \eqref{eq:b-initial-slope} for \(w_0^\varepsilon+f\), we use \eqref{eq:b-g-admissible}, \eqref{bstar}, \eqref{eq:b-data-rescaled-slope}, and \eqref{eq:b-data-perturbation-initial-weights} to obtain
\begin{equation*}
\begin{split} &\left|\varepsilon(w_0^\varepsilon+f)'(x)-U'\left(\frac{x}{\varepsilon^{(2\ell+1)/2}}\right)\right|\\ &\qquad\leq\left(\rho(1-2\vartheta)+C_\varepsilon\|(1+x^2)^{1/2}f\|_{H^5}\right)\min\left\{\frac{Nx^2}{\varepsilon^{2\ell+1}+x^2},\frac{\varepsilon B_\ell}{\varepsilon+|x|^{2/(2\ell+1)}}\right\}. \end{split}
\end{equation*}
Since \(\rho<1\), \eqref{eq:b-initial-slope} holds for \(w_0^\varepsilon+f\) when \(\|(1+x^2)^{1/2}f\|_{H^5}\) is sufficiently small. For \eqref{eq:b-initial-second}, using \eqref{eq:b-Phi} and \eqref{eq:b-data-perturbation-initial-weights}, we obtain
\begin{equation*}
|(w_0^\varepsilon+f)''(x)|\leq\left(\varepsilon^{-(2\ell+3)/2}\left\|\frac{\sqrt{1+y^2}}{y}v'(y)\right\|_{L^\infty}+C_\varepsilon\|(1+x^2)^{1/2}f\|_{H^5}\right)\frac{|x|}{\sqrt{\varepsilon^{2\ell+1}+x^2}}.
\end{equation*}
By the strict bound \eqref{eq:b-data-second}, \eqref{eq:b-initial-second} also holds for \(w_0^\varepsilon+f\) when \(\|(1+x^2)^{1/2}f\|_{H^5}\) is sufficiently small. At \(x=0\), equality in \eqref{eq:b-initial-slope} and \eqref{eq:b-initial-second} is preserved by the matching conditions \eqref{eq:b-initial-matching}.

It remains to verify \eqref{eq:b-initial-origin}, \eqref{eq:b-initial-third}--\eqref{eq:b-initial-fourth}, and \eqref{eq:b-initial-energy} for \(w_0^\varepsilon+f\). By \eqref{eq:b-data-perturbation-higher}, the norm \(\|(1+x^2)^{1/2}f\|_{H^5}\) controls \(|f^{(3)}(0)|\), \(\|f^{(3)}\|_{L^\infty}\), and \(\|f^{(4)}\|_{L^\infty}\). Since \(1/(2\ell+1)<1/2\), Sobolev embedding also yields
\begin{equation*}
\|(1+x^2)^{1/(2\ell+1)}f\|_{L^\infty} \leq \|(1+x^2)^{1/2}f\|_{L^\infty}\leq C\|(1+x^2)^{1/2}f\|_{H^5}.
\end{equation*}
By \eqref{eq:b-Phi} and the strict bounds \eqref{eq:b-data-higher}, the bounds \eqref{eq:b-initial-third}--\eqref{eq:b-initial-fourth} hold with strict inequalities for \(w_0^\varepsilon\). The same is true of \eqref{eq:b-initial-origin} and \eqref{eq:b-initial-energy}. Hence all these conditions hold for \(w_0^\varepsilon+f\) when \(\|(1+x^2)^{1/2}f\|_{H^5}\) is sufficiently small. Thus a relative neighborhood of \(w_0^\varepsilon\) is contained in the admissible set. The admissible set therefore has non-empty relative interior. Together with Theorem~\ref{thm:b-cusp}, this completes the proof of Theorem~\ref{thm:b-cusp-intro}.

\section{Proof of Lemma~\ref{lem:b-local}} \label{app:b-local}

We prove Lemma~\ref{lem:b-local} following the local well-posedness and continuation arguments of \cite[Sections~2--3]{EscherYin2008}. Relative to the standard $b$-equation treated in \cite{EscherYin2008}, the only additional term in the nonlocal formulation \eqref{eq:b-nonlocal} is
\begin{equation*}
\mathcal R(w):=-\partial_x(1-\partial_x^2)^{-1}h(w),\qquad h(w):=g(w)-\frac b2w^2.
\end{equation*}

Fix an integer $r\geq5$. We first verify the conditions on the additional term $\mathcal R$ required in the proof of \cite[Theorem~2.1]{EscherYin2008}. Let $j\geq4$ be an integer and $w\in H^j(\mathbb R)$. By Sobolev embedding, $w\in L^\infty(\mathbb R)$. Since $h(0)=0$, the mean value theorem gives
\begin{equation*}
\|h(w)\|_{L^2}\leq\sup_{|z|\leq\|w\|_{L^\infty}}|h'(z)|\|w\|_{L^2}<\infty.
\end{equation*}
By the Moser composition estimate \cite[Remark~2.2]{TexierZumbrun2011}, we obtain
\begin{equation}\label{moser}
\|h(w)\|_{H^j}\leq C_{h,j}(\|w\|_{L^\infty})\|w\|_{H^j}.
\end{equation}
Thus $h$ maps $H^j(\mathbb R)$ into itself. Moreover, the Fourier transform gives
\begin{equation*}
\widehat{\mathcal R(w)}(\xi)=-\frac{i\xi}{1+\xi^2}\widehat{h(w)}(\xi),
\end{equation*}
and hence
\begin{equation*}
\|\mathcal R(w)\|_{H^{j+1}}^2=\int_{\mathbb R}(1+\xi^2)^{j+1}\frac{\xi^2}{(1+\xi^2)^2}|\widehat{h(w)}(\xi)|^2\,d\xi\leq\|h(w)\|_{H^j}^2.
\end{equation*}
Therefore, for any $R\geq\|w\|_{L^\infty}$,
\begin{equation*}
\|\mathcal R(w)\|_{H^{j+1}}\leq C_{h,j,R}\|w\|_{H^j}.
\end{equation*}

We next establish the local Lipschitz continuity of $h$ and $\mathcal R$. Let $v\in H^j(\mathbb R)$ and suppose that $\|w\|_{H^j},\|v\|_{H^j}\leq M$. By the fundamental theorem of calculus,
\begin{equation*}
\begin{aligned}
h(w)-h(v)&=(w-v)\int_0^1h'(v+\theta(w-v))\,d\theta\\
&=h'(0)(w-v)+(w-v)\int_0^1\bigl(h'(v+\theta(w-v))-h'(0)\bigr)\,d\theta.
\end{aligned}
\end{equation*}
Since $H^j(\mathbb R)$ is an algebra, applying the Moser estimate to $h'-h'(0)$ gives, for every $\theta\in[0,1]$,
\begin{equation*}
\begin{aligned}
&\|(w-v)(h'(v+\theta(w-v))-h'(0))\|_{H^j}\\&\qquad\leq C\|w-v\|_{H^j}\|h'(v+\theta(w-v))-h'(0)\|_{H^j}\leq C_M\|w-v\|_{H^j},
\end{aligned}
\end{equation*}
where $C_M$ is independent of $\theta$. Consequently,
\begin{equation*}
\|h(w)-h(v)\|_{H^j}\leq C_M\|w-v\|_{H^j}.
\end{equation*}
Using the Fourier multiplier estimate, we also obtain
\begin{equation*}
\|\mathcal R(w)-\mathcal R(v)\|_{H^{j+1}}\leq\|h(w)-h(v)\|_{H^j}\leq C_M\|w-v\|_{H^j}.
\end{equation*}
Applying these bounds first with $j=r-1$ and then with $j=r$, we find that $\mathcal R$ is bounded on bounded subsets of $H^r$ and satisfies the source conditions (a) and (b) in the proof of \cite[Theorem~2.1]{EscherYin2008}. The transport operator $w\partial_x$ remains unchanged from the standard $b$-equation. We may therefore apply the Kato theory used in that proof to the present equation. This yields a unique maximal solution
\begin{equation*}
w\in C([0,T_r);H^r(\mathbb R))\cap C^1([0,T_r);H^{r-1}(\mathbb R)),
\end{equation*}
where $T_r$ denotes the maximal existence time in $H^r(\mathbb R)$.

We next show that $T_r=T_5$ for every integer $r\geq5$. By uniqueness in $H^5$, the $H^r$ solution agrees with the $H^5$ solution on their common interval of existence. Since $H^r(\mathbb R)\subset H^5(\mathbb R)$, the maximality of $T_5$ gives $T_r\leq T_5$. 

To prove that $T_r\geq T_5$, we establish an $H^r$ energy estimate. Let $\Gamma^r:=(1-\partial_x^2)^{r/2}$. We first derive the estimate for smooth solutions. Applying $\Gamma^r$ to \eqref{u-nonlocal}, we obtain
\begin{equation*}
\partial_t\Gamma^r w+w\partial_x\Gamma^r w+[\Gamma^r,w]w_x+\Gamma^r P_x=0,
\end{equation*}
where $[\Gamma^r,w]w_x:=\Gamma^r(ww_x)-w\Gamma^r w_x$. Multiplying by $\Gamma^r w$ and integrating over $\mathbb R$, we have
\begin{equation}\label{eq:b-sobolev-energy}
\begin{split}
\frac12\frac{d}{dt}\|w\|_{H^r}^2&=-\int_{\mathbb R}w\Gamma^r w\,\partial_x\Gamma^r w\,dx-\int_{\mathbb R}\Gamma^r w\,[\Gamma^r,w]w_x\,dx-\int_{\mathbb R}\Gamma^r w\,\Gamma^r P_x\,dx\\
&\leq\frac12\|w_x\|_{L^\infty}\|w\|_{H^r}^2+\|w\|_{H^r}\|[\Gamma^r,w]w_x\|_{L^2}+\|w\|_{H^r}\|P_x\|_{H^r}.
\end{split}
\end{equation}
By the Kato--Ponce commutator estimate \cite[Lemma~3.2]{EscherYin2008},
\begin{equation}\label{eq:b-commutator-bound}
\begin{split}
\|[\Gamma^r,w]w_x\|_{L^2}&\leq C_r\bigl(\|w_x\|_{L^\infty}\|w_x\|_{H^{r-1}}+\|w\|_{H^r}\|w_x\|_{L^\infty}\bigr) \leq C_r\|w_x\|_{L^\infty}\|w\|_{H^r}.
\end{split}
\end{equation}
For the nonlocal term $P_x$, the Fourier multiplier estimate for $\partial_x(1-\partial_x^2)^{-1}$ yields
\begin{equation*}
\|P_x\|_{H^r}\leq C\left\|g(w)+\frac{3-b}{2}w_x^2\right\|_{H^{r-1}}\leq C_b\left(\|g(w)\|_{H^{r-1}}+\|w_x^2\|_{H^{r-1}}\right).
\end{equation*}
Since $g\in C^\infty$, $g(0)=0$, and $\|w\|_{L^\infty}\leq R$, the Moser composition estimate yields
\begin{equation*}
\|g(w)\|_{H^{r-1}}\leq C_{g,r,R}\|w\|_{H^{r-1}}\leq C_{g,r,R}\|w\|_{H^r}.
\end{equation*}
Moreover, the Moser product estimate \cite[Lemma~3.1]{EscherYin2008} yields
\begin{equation*}
\|w_x^2\|_{H^{r-1}}\leq C_r\|w_x\|_{L^\infty}\|w_x\|_{H^{r-1}}\leq C_r\|w_x\|_{L^\infty}\|w\|_{H^r}.
\end{equation*}
Consequently,
\begin{equation}\label{eq:b-pressure-sobolev-bound}
\|P_x\|_{H^r}\leq C_{b,g,r,R}(1+\|w_x\|_{L^\infty})\|w\|_{H^r}.
\end{equation}
Combining \eqref{eq:b-sobolev-energy}, \eqref{eq:b-commutator-bound}, and \eqref{eq:b-pressure-sobolev-bound}, we obtain
\begin{equation}\label{eq:b-sobolev-tame}
\frac{d}{dt}\|w\|_{H^r}^2\leq C_{b,g,r,R}(1+\|w_x\|_{L^\infty})\|w\|_{H^r}^2.
\end{equation}
A standard approximation argument extends the integrated form of \eqref{eq:b-sobolev-tame}, derived above for smooth solutions, to the $H^r$ solution.

Suppose that $T_r<T_5$. Since the $H^5$ solution is continuous on $[0,T_r]$, the quantities $\|w\|_{L^\infty}$ and $\|w_x\|_{L^\infty}$ are bounded on $[0,T_r)$. Applying the integrated form of \eqref{eq:b-sobolev-tame} and Gronwall's inequality, we obtain a uniform $H^r$ bound on $[0,T_r)$. The local well-posedness theory gives an existence time bounded below on bounded subsets of $H^r$. We may therefore extend the $H^r$ solution beyond $T_r$, contradicting its maximality. Thus $T_r\geq T_5$, and hence $T_r=T_5$ for every integer $r\geq5$. We set $T_{\max}:=T_5$.

Moreover, \eqref{u-nonlocal} and \eqref{P-nonlocal}, together with the product and composition estimates, imply that
\begin{equation*}
w_t=-ww_x-\partial_x(1-\partial_x^2)^{-1}\left(g(w)+\frac{3-b}{2}w_x^2\right)\in C([0,T_{\max});H^{r-1}(\mathbb R))
\end{equation*}
for every integer $r\geq5$. By the continuous embeddings between Sobolev spaces, we conclude that
\begin{equation*}
w\in C([0,T_{\max});H^r(\mathbb R))\cap C^1([0,T_{\max});H^{r-1}(\mathbb R))
\end{equation*}
for every real $r\geq0$.

We next prove the regularity of $(1+x^2)^{1/2}w$. Multiplying \eqref{u-nonlocal} by $(1+x^2)^{1/2}$, we obtain
\begin{equation}\label{weight_eq}
\begin{aligned}
&\partial_t\bigl((1+x^2)^{1/2}w\bigr)+w\partial_x\bigl((1+x^2)^{1/2}w\bigr)\\&\qquad=\frac{x}{(1+x^2)^{1/2}}w^2-(1+x^2)^{1/2}\left[K'*\left(g(w)+\frac{3-b}{2}w_x^2\right)\right].
\end{aligned}
\end{equation}
We first derive weighted a priori estimates under the assumption that $(1+x^2)^{1/2}w\in H^5(\mathbb R)$. We then justify these estimates using bounded weights.

Since $(1+x^2)^{1/2}\leq\sqrt2(1+(x-z)^2)^{1/2}(1+z^2)^{1/2}$ holds for any $x,z\in\mathbb R$, we obtain
\begin{equation*}
\begin{aligned}
(1+x^2)^{1/2}|(K'*f)(x)|&\leq\int_{\mathbb R}(1+x^2)^{1/2}|K'(x-z)||f(z)|\,dz\\&\leq\sqrt2\int_{\mathbb R}(1+(x-z)^2)^{1/2}|K'(x-z)|(1+z^2)^{1/2}|f(z)|\,dz,
\end{aligned}
\end{equation*}
where $f:=g(w)+\frac{3-b}{2}w_x^2$. Hence, by Young's inequality, we obtain
\begin{equation*}
\begin{aligned}
\|(1+x^2)^{1/2}(K'*f)\|_{L^2}&\leq\sqrt2\|(1+x^2)^{1/2}K'\|_{L^1}\|(1+x^2)^{1/2}f\|_{L^2} \leq C\|(1+x^2)^{1/2}f\|_{L^2}.
\end{aligned}
\end{equation*}
Here the last inequality follows from $K'(x)=-\frac12\operatorname{sgn}(x)e^{-|x|}$. Similarly, for each integer $0\leq j\leq4$, we obtain
\begin{equation*}
\|\partial_x^j[(1+x^2)^{1/2}(K'*f)]\|_{L^2}\leq C\sum_{m=0}^j\|(1+x^2)^{1/2}\partial_x^mf\|_{L^2}\leq C\|(1+x^2)^{1/2}f\|_{H^j}.
\end{equation*}
In the last inequality, we used the boundedness of $(1+x^2)^{1/2}\partial_x^m(1+x^2)^{-1/2}$ for $0\leq m\leq j$. For the fifth derivative of $K'*f$, using $K''=K-\delta_0$, we obtain
\begin{equation*}
\partial_x^5(K'*f)=K''*\partial_x^4f=K*\partial_x^4f-\partial_x^4f,
\end{equation*}
and hence
\begin{equation*}
\begin{aligned}
\|(1+x^2)^{1/2}\partial_x^5(K'*f)\|_{L^2}&\leq\|(1+x^2)^{1/2}(K*\partial_x^4f)\|_{L^2}+\|(1+x^2)^{1/2}\partial_x^4f\|_{L^2}\\
&\leq C\|(1+x^2)^{1/2}f\|_{H^4}.
\end{aligned}
\end{equation*}
Combining the estimates for $K'*f$ and its derivatives, we obtain
\begin{equation}\label{nonlocal}
\begin{aligned}
\|(1+x^2)^{1/2}(K'*f)\|_{H^5}&\leq C\|(1+x^2)^{1/2}f\|_{H^4}\\&\leq C_b\left(\|(1+x^2)^{1/2}g(w)\|_{H^4}+\|(1+x^2)^{1/2}w_x^2\|_{H^4}\right).
\end{aligned}
\end{equation}

For the first term in \eqref{nonlocal}, we define $\tilde h(w):=\int_0^1g'(\theta w)\,d\theta$. Since $g(0)=0$, we have $g(w)=w\tilde h(w)$. For $\|w\|_{L^\infty}\leq R$, the Moser composition estimate gives
\begin{equation*}
\|\tilde h(w)-\tilde h(0)\|_{H^4}\leq C_{g,R}\|w\|_{H^4}.
\end{equation*}
Using the Moser product estimate and the bound on $\tilde h(w)-\tilde h(0)$, we obtain
\begin{equation}\label{first}
\begin{aligned}
\|(1+x^2)^{1/2}g(w)\|_{H^4}&=\|(1+x^2)^{1/2}w\tilde h(w)\|_{H^4}\\
&\leq|\tilde h(0)|\|(1+x^2)^{1/2}w\|_{H^4}+\|(1+x^2)^{1/2}w(\tilde h(w)-\tilde h(0))\|_{H^4}\\
&\leq C_{g,R}(1+\|w\|_{H^4})\|(1+x^2)^{1/2}w\|_{H^4}\\
&\leq C_{g,R}(1+\|w\|_{H^5})\|(1+x^2)^{1/2}w\|_{H^5}.
\end{aligned}
\end{equation}
To estimate the second term in \eqref{nonlocal}, we first use
\begin{equation*}
(1+x^2)^{1/2}w_x=\partial_x\bigl((1+x^2)^{1/2}w\bigr)-\frac{x}{1+x^2}(1+x^2)^{1/2}w
\end{equation*}
to obtain
\begin{equation*}
\begin{aligned}
\|(1+x^2)^{1/2}w_x\|_{H^4}&\leq\|\partial_x((1+x^2)^{1/2}w)\|_{H^4}+\left\|\frac{x}{1+x^2}(1+x^2)^{1/2}w\right\|_{H^4}\\&\leq C\|(1+x^2)^{1/2}w\|_{H^5}.
\end{aligned}
\end{equation*}
This implies that
\begin{equation}\label{second}
\|(1+x^2)^{1/2}w_x^2\|_{H^4}\leq C\|(1+x^2)^{1/2}w_x\|_{H^4}\|w_x\|_{H^4}\leq C\|w\|_{H^5}\|(1+x^2)^{1/2}w\|_{H^5}.
\end{equation}
Substituting \eqref{first} and \eqref{second} into \eqref{nonlocal}, we have
\begin{equation}\label{eq:b-weighted-Px-H5}
\|(1+x^2)^{1/2}P_x\|_{H^5}\leq C_{b,g,R}(1+\|w\|_{H^5})\|(1+x^2)^{1/2}w\|_{H^5},
\end{equation}
whenever $\|w\|_{L^\infty}\leq R$.

To prove the weighted regularity, we introduce the smooth bounded weights
\begin{equation*}
\omega_N(x):=\frac{N(1+x^2)^{1/2}}{N+(1+x^2)^{1/2}},\qquad N\geq1.
\end{equation*}
Then $\omega_N(x)\to(1+x^2)^{1/2}$ as $N\to\infty$. In the estimates below, all constants are independent of $N$. Applying the calculation leading to \eqref{weight_eq} with $\omega_N$ in place of $(1+x^2)^{1/2}$, we obtain
\begin{equation}\label{eq:b-cutoff-weight-equation}
\partial_t(\omega_Nw)+w\partial_x(\omega_Nw)=\omega_N'w^2-\omega_NP_x.
\end{equation}
Moreover,
\begin{equation}\label{eq:b-cutoff-weight-bounds}
\begin{aligned}
&\sup_{N\geq1}\left\|\frac{\omega_N'}{\omega_N}\right\|_{W^{5,\infty}}<\infty,\qquad \sup_{N\geq1}\sum_{j=1}^5\left\|\frac{\partial_x^j\omega_N}{\omega_N}\right\|_{L^\infty}<\infty,\\&\omega_N(x)\leq C(1+(x-y)^2)^{1/2}\omega_N(y),\qquad x,y\in\mathbb R.
\end{aligned}
\end{equation}
The derivative bounds in \eqref{eq:b-cutoff-weight-bounds} imply that $\omega_N\partial_x^k(\omega_N^{-1})$ is uniformly bounded for $0\leq k\leq4$. Thus, by the product rule,
\begin{equation*}
\omega_N\partial_x^j f=\sum_{k=0}^j\binom{j}{k}\bigl(\omega_N\partial_x^k(\omega_N^{-1})\bigr)\partial_x^{j-k}(\omega_Nf),\qquad 0\leq j\leq4,
\end{equation*}
and consequently
\begin{equation*}
\sum_{j=0}^4\|\omega_N\partial_x^jf\|_{L^2}\leq C\|\omega_Nf\|_{H^4}.
\end{equation*}
Using this bound, the last two bounds in \eqref{eq:b-cutoff-weight-bounds}, and the calculation leading to \eqref{nonlocal}, we obtain
\begin{equation*}
\|\omega_N(K'*f)\|_{H^5}\leq C\|\omega_Nf\|_{H^4}.
\end{equation*}
To estimate $\|\omega_Nf\|_{H^4}$, we use \eqref{first}--\eqref{second} with $\omega_N$ in place of $(1+x^2)^{1/2}$. In particular, the first bound in \eqref{eq:b-cutoff-weight-bounds} and the identity
\begin{equation*}
\omega_Nw_x=\partial_x(\omega_Nw)-\frac{\omega_N'}{\omega_N}\omega_Nw
\end{equation*}
give $\|\omega_Nw_x\|_{H^4}\leq C\|\omega_Nw\|_{H^5}$. Consequently, for $\|w\|_{L^\infty}\leq R$,
\begin{equation}\label{eq:b-cutoff-Px}
\|\omega_NP_x\|_{H^5}\leq C_{b,g,R}(1+\|w\|_{H^5})\|\omega_Nw\|_{H^5}.
\end{equation}
The first bound in \eqref{eq:b-cutoff-weight-bounds} and the product estimate also give
\begin{equation}\label{eq:b-cutoff-weight-source}
\|\omega_N'w^2\|_{H^5}=\left\|\frac{\omega_N'}{\omega_N}w(\omega_Nw)\right\|_{H^5}\leq C\|w\|_{H^5}\|\omega_Nw\|_{H^5}.
\end{equation}
The energy calculation for $\omega_Nw$ is justified by the boundedness of $\omega_N$ and its derivatives and the previously established regularity $w\in C([0,T_{\max});H^6)\cap C^1([0,T_{\max});H^5)$.

Fix $T<T_{\max}$. Applying $\Gamma^5$ to \eqref{eq:b-cutoff-weight-equation}, multiplying by $\Gamma^5(\omega_Nw)$, and integrating over $\mathbb R$, we use the commutator estimate \cite[Lemma~3.2]{EscherYin2008}, \eqref{eq:b-cutoff-Px}, and \eqref{eq:b-cutoff-weight-source} to obtain
\begin{equation}\label{eq:b-cutoff-weight-energy}
\left|\frac{d}{dt}\|\omega_Nw(t)\|_{H^5}^2\right|\leq C_T\|\omega_Nw(t)\|_{H^5}^2,\qquad 0\leq t\leq T.
\end{equation}
By Gronwall's inequality,
\begin{equation} \label{eq:b-cutoff-H5-bound}
\sup_{0\leq t\leq T}\|\omega_Nw(t)\|_{H^5}\leq C_T\|\omega_Nw_0\|_{H^5}\leq C_T\|(1+x^2)^{1/2}w_0\|_{H^5}.
\end{equation}
Here we used the uniform boundedness of $\omega_N/(1+x^2)^{1/2}$ in $W^{5,\infty}$.  For each $t\in[0,T]$, $\omega_Nw(t)\to(1+x^2)^{1/2}w(t)$ in the sense of distributions. By \eqref{eq:b-cutoff-H5-bound} and weak compactness in $H^5$, we obtain
\begin{equation}\label{eq:b-weighted-H5-uniform}
\sup_{0\leq t\leq T}\|(1+x^2)^{1/2}w(t)\|_{H^5}<\infty.
\end{equation}
For each $t\in[0,T]$, we also have $\omega_Nw(t)\to(1+x^2)^{1/2}w(t)$ in $H^5$, since $\omega_N/(1+x^2)^{1/2}$ and its derivatives up to order five are uniformly bounded and converge pointwise to those of $1$. Applying Gronwall's inequality to \eqref{eq:b-cutoff-weight-energy} and letting $N\to\infty$, we obtain
\begin{equation*}
\|(1+x^2)^{1/2}w(t)\|_{H^5}\leq e^{C_T|t-s|}\|(1+x^2)^{1/2}w(s)\|_{H^5},\qquad s,t\in[0,T].
\end{equation*}
Thus the $H^5$ norm of $(1+x^2)^{1/2}w(t)$ is continuous in $t$. The continuity of $w$ in $H^5$ gives continuity of $(1+x^2)^{1/2}w$ in the sense of distributions. This continuity and \eqref{eq:b-weighted-H5-uniform} imply weak continuity of $(1+x^2)^{1/2}w$ in $H^5$. The continuity of its norm then yields
\begin{equation*}
(1+x^2)^{1/2}w\in C([0,T];H^5(\mathbb R)).
\end{equation*}
Since $T<T_{\max}$ is arbitrary, the weighted regularity holds on $[0,T_{\max})$.

Finally, we prove the blow-up criterion \eqref{eq:b-local-continuation}. Define the energy $E$ by
\begin{equation}\label{eq:b-energy-def}
E(t):=\int_{\mathbb R}(w^2+w_x^2)\,dx.
\end{equation}
Let $H_g(v):=\int_0^v g(z)\,dz$, so that $\partial_xH_g(w)=g(w)w_x$. By \eqref{u-nonlocal} and integration by parts,\footnote{For any $R\geq\|w\|_{L^\infty}$, let $C_R:=\sup_{|z|\leq R}|g'(z)|$. Since $g(0)=0$, we have $|g(w)|\leq C_R|w|$ and $|H_g(w)|\leq C_R|w|^2/2$. Hence $H_g(w)\in L^1(\mathbb R)$ and $g(w)w_x\in L^1(\mathbb R)$, which justifies the integration by parts.} we obtain
\begin{equation}\label{eq:b-energy-identity}
\begin{split}
E'(t)&=2\int_{\mathbb R}(ww_t+w_xw_{xt})\,dx=2\int_{\mathbb R}(w-w_{xx})w_t\,dx\\&=-2\int_{\mathbb R}(w-w_{xx})(ww_x+P_x)\,dx=-\int_{\mathbb R}w_x^3\,dx+2\int_{\mathbb R}(w_x-w_{xxx})P\,dx\\&=-\int_{\mathbb R}w_x^3\,dx+2\int_{\mathbb R}w_x\left(g(w)+\frac{3-b}{2}w_x^2\right)\,dx=(2-b)\int_{\mathbb R}w_x^3\,dx,
\end{split}
\end{equation}
which implies that
\begin{equation*}
|E'(t)|\leq|2-b|\int_{\mathbb R}|w_x|^3\,dx\leq|2-b|\|w_x\|_{L^\infty}E(t).
\end{equation*}
Suppose that $T_{\max}<\infty$ and $\limsup_{t\nearrow T_{\max}}\|w_x(\cdot,t)\|_{L^\infty}<\infty$. Since $w\in C([0,T_{\max});H^5(\mathbb R))$, we have $\sup_{0\leq t<T_{\max}}\|w_x(\cdot,t)\|_{L^\infty}<\infty$. By \eqref{eq:b-energy-identity} and Gronwall's inequality,
\begin{equation*}
E(t)\leq E(0)\exp\left(|2-b|\int_0^t\|w_x(\cdot,\tau)\|_{L^\infty}\,d\tau\right).
\end{equation*}
Hence $\|w\|_{L^\infty}$ remains bounded on $[0,T_{\max})$. Applying the integrated form of \eqref{eq:b-sobolev-tame} with $r=5$ and using Gronwall's inequality again, we obtain a uniform $H^5$ bound on $[0,T_{\max})$. The local well-posedness argument gives a positive existence time depending only on an upper bound for the $H^5$ norm of the initial data. We may therefore restart the solution at times approaching $T_{\max}$ with a uniform positive existence time. This extends the solution beyond $T_{\max}$, contradicting its maximality. Thus \eqref{eq:b-local-continuation} follows.

\section{Estimates for transport-type equations} \label{app:transport-estimates}

We consider transport-type equations of the form
\begin{equation}\label{eq:general-transport}
\partial_s f(y,s)+\mathcal V(y,s)\partial_y f(y,s)+D(y,s)f(y,s)=\mathcal F(y,s),\qquad (y,s)\in\R\times[s_0,s_1],
\end{equation}
where \(\mathcal F\) may contain nonlocal terms. We record the far-field estimate and maximum principle for \eqref{eq:general-transport} used in Section~\ref{sec:perturbative}. The proofs follow \cite[Lemma~3.5]{BKK-delta} and Lemma~SM2.1 in its supplementary material, respectively. In the following lemmas, we work with differential inequalities for \(|f|\) obtained from such equations.

\subsection{Far-field behavior}

The following lemma gives a bound at spatial infinity for solutions of a transport-type equation.

\begin{lemma}\label{lem:transport-far-field}
Let \(I=[s_0,s_1]\), and let \(f\in C(\R\times I)\cap C(I;L^\infty(\R))\) be \(C^1\) on \(\{|y|>R_0\}\times I\) for some \(R_0\geq0\). Suppose that $f$ satisfies
\begin{equation}\label{eq:transport-tail-equation}
\partial_s|f(y,s)|+\mathcal V(y,s)\partial_y|f(y,s)|+D(y,s)|f(y,s)|\leq F(y,s)
\end{equation}
for \(|y|>R_0\) whenever \(f(y,s)\neq0\), where $\mathcal{V},\partial_y\mathcal{V} \in C (\mathbb{R} \times I)$ and $D,F \in C (\{|y| > R_0\} \times I)$. Assume also that
\begin{subequations}\label{eq:transport-tail-assumptions}
\begin{align}
&|\mathcal V(y,s)|\leq C(1+|y|),\qquad (y,s)\in\R\times I,\label{eq:transport-tail-growth}\\
&D(y,s)\geq\lambda_D-a(s),\qquad 0\leq F(y,s)\leq e(s),\qquad |y|>R_0,\quad s\in I,\label{eq:transport-tail-coefficients}
\end{align}
\end{subequations}
where \(C<\infty\), \(\lambda_D\geq0\), and \(a,e\) are nonnegative continuous functions on \(I\). Then
\begin{equation}\label{eq:transport-tail-bound}
\begin{split}
\limsup_{|y|\to\infty}|f(y,s)| &\leq e^{-\lambda_D(s-s_0)+\int_{s_0}^sa(\sigma)\,d\sigma}\limsup_{|y|\to\infty}|f(y,s_0)|\\
&\quad +\int_{s_0}^se^{-\lambda_D(s-\sigma)+\int_\sigma^sa(\rho)\,d\rho}e(\sigma)\,d\sigma,\qquad s\in I.
\end{split}
\end{equation}
\end{lemma}

\begin{proof}
Let \(\psi(y,s)\) be the solution to
\begin{equation}\label{eq:transport-tail-characteristics}
\partial_s\psi(y,s)=\mathcal V(\psi(y,s),s),\qquad \psi(y,s_0)=y.
\end{equation}
We first claim that \(\psi\) is defined on \(\R\times I\) and satisfies
\begin{equation}\label{eq:transport-tail-flow-bounds}
(1+|y|)e^{-C(s-s_0)}\leq1+|\psi(y,s)|\leq(1+|y|)e^{C(s-s_0)},\qquad (y,s)\in\R\times I.
\end{equation}
Indeed, on the interval of existence of \(\psi\), \eqref{eq:transport-tail-characteristics} and \eqref{eq:transport-tail-growth} give
\begin{equation*}
|\partial_s\psi(y,s)|\leq C(1+|\psi(y,s)|),
\end{equation*}
and hence
\begin{equation*}
\left|\frac{d}{ds}\log(1+|\psi(y,s)|)\right|\leq C
\end{equation*}
for almost every \(s\). Integrating from \(s_0\) to \(s\), we obtain \eqref{eq:transport-tail-flow-bounds} on the interval of existence. The upper bound in \eqref{eq:transport-tail-flow-bounds}, together with the continuity of \(\mathcal V\) and \(\partial_y\mathcal V\), implies that \(\psi\) is defined throughout \(I\). This proves the claim.

By \eqref{eq:transport-tail-flow-bounds}, if \(|y|>(R_0+1)e^{C(s_1-s_0)}-1\), then \(|\psi(y,s)|>R_0\) for all \(s\in I\). For such \(y\), set
\begin{equation*}
g(s):=|f(\psi(y,s),s)|.
\end{equation*}
Then, by \eqref{eq:transport-tail-equation},
\begin{equation}\label{eq:transport-tail-along-flow}
g'(s)+D(\psi(y,s),s)g(s)\leq F(\psi(y,s),s)
\end{equation}
whenever \(g(s)>0\). Moreover, \(g\) is absolutely continuous and \(g'=0\) almost everywhere on \(\{g=0\}\). Since \(F\geq0\), \eqref{eq:transport-tail-along-flow} also holds almost everywhere on \(\{g=0\}\). Thus \eqref{eq:transport-tail-along-flow} holds for almost every \(s\in I\).

Integrating \eqref{eq:transport-tail-along-flow} and using \eqref{eq:transport-tail-coefficients}, we obtain
\begin{equation}\label{eq:transport-tail-characteristic-bound}
\begin{split}
|f(\psi(y,s),s)| &\leq |f(y,s_0)|e^{-\int_{s_0}^{s}D(\psi(y,\sigma),\sigma)\,d\sigma}+\int_{s_0}^{s}e^{-\int_{\sigma}^{s}D(\psi(y,\rho),\rho)\,d\rho}F(\psi(y,\sigma),\sigma)\,d\sigma\\ 
&\leq |f(y,s_0)|e^{-\lambda_D(s-s_0)+\int_{s_0}^sa(\sigma)\,d\sigma}+\int_{s_0}^se^{-\lambda_D(s-\sigma)+\int_\sigma^sa(\rho)\,d\rho}e(\sigma)\,d\sigma
\end{split}
\end{equation}
for all \(y\) satisfying \(|y|>(R_0+1)e^{C(s_1-s_0)}-1\) and all \(s\in I\).

To finish the proof, we show that, for each \(s\in I\),
\begin{equation}\label{eq:transport-tail-change-variables}
\limsup_{|y|\to\infty}|f(\psi(y,s),s)|=\limsup_{|y|\to\infty}|f(y,s)|.
\end{equation}
We let \(\psi_s(y):=\psi(y,s)\). For each \(s\in I\), the map \(y\mapsto\psi_s(y)\) is strictly increasing and onto. Using these properties of \(\psi_s\), we have, for any \(R>0\),
\begin{equation}\label{eq:transport-tail-suprema}
\sup_{|y|\geq R}|f(\psi_s(y),s)|=\sup_{z\in(-\infty,\psi_s(-R)]\cup[\psi_s(R),\infty)}|f(z,s)|.
\end{equation}
Indeed, since \(\psi_s\) is increasing, the left-hand side is bounded above by the right-hand side. On the other hand, for each \(z\geq\psi_s(R)\), the surjectivity of \(\psi_s\) gives \(y\in\R\) such that \(\psi_s(y)=z\). The monotonicity of \(\psi_s\) then gives \(y\geq R\). Similarly, for each \(z\leq\psi_s(-R)\), there exists \(y\leq-R\) such that \(\psi_s(y)=z\). Hence
\begin{equation*}
|f(z,s)|\leq\sup_{|y|\geq R}|f(\psi_s(y),s)|
\end{equation*}
for all \(z\in(-\infty,\psi_s(-R)]\cup[\psi_s(R),\infty)\). This proves \eqref{eq:transport-tail-suprema}.

Set
\begin{equation*}
A:=\limsup_{|y|\to\infty}|f(\psi_s(y),s)|,\qquad B:=\limsup_{|y|\to\infty}|f(y,s)|.
\end{equation*}
Since \(\psi_s\) is strictly increasing and onto, \(\psi_s(R)\to\infty\) and \(\psi_s(-R)\to-\infty\) as \(R\to\infty\). Thus \(\min\{\psi_s(R),-\psi_s(-R)\}\to\infty\) as \(R\to\infty\). For sufficiently large \(R\), \eqref{eq:transport-tail-suprema} gives
\begin{equation*}
\sup_{|y|\geq R}|f(\psi_s(y),s)|\leq\sup_{|z|\geq \min\{\psi_s(R),-\psi_s(-R)\}}|f(z,s)|.
\end{equation*}
Taking \(R\to\infty\), we obtain \(A\leq B\).

To show \(A=B\), we suppose that \(A<B\). We choose \(\varepsilon>0\) sufficiently small so that \(A+\varepsilon<B-\varepsilon\). For sufficiently large \(R_1>0\), we have
\begin{equation*}
\sup_{|y|\geq R_1}|f(\psi_s(y),s)|<A+\varepsilon.
\end{equation*}
We then choose \(R_2>\max\{\psi_s(R_1),-\psi_s(-R_1)\}\). By \eqref{eq:transport-tail-suprema},
\begin{equation*}
\sup_{|y|\geq R_1}|f(\psi_s(y),s)|\geq\sup_{|z|\geq R_2}|f(z,s)|\geq B>B-\varepsilon.
\end{equation*}
This contradicts \(A+\varepsilon<B-\varepsilon\). Thus \eqref{eq:transport-tail-change-variables} holds. Taking the upper limit as \(|y|\to\infty\) in \eqref{eq:transport-tail-characteristic-bound}, we obtain \eqref{eq:transport-tail-bound} and finish the proof of Lemma~\ref{lem:transport-far-field}.
\end{proof}

\subsection{Maximum principle}

We next give a maximum principle that bounds \(\|f(\cdot,s)\|_{L^\infty}\) by a time-dependent comparison function.

\begin{lemma}\label{lem:bootstrap-comparison}
Let \(I=[s_0,s_1]\), and let \(\Omega\subset\R\) be a compact set. Suppose that \(f\in C(\R\times I)\cap C(I;L^\infty(\R))\) is \(C^1\) on \(\Omega^c\times I\) and satisfies
\begin{equation}\label{eq:comparison-inequality}
\partial_s|f(y,s)|+\mathcal V(y,s)\partial_y|f(y,s)|+D(y,s)|f(y,s)|\leq K(y,s)\|f(\cdot,s)\|_{L^\infty}+H(y,s)
\end{equation}
for \((y,s)\in\Omega^c\times I\) whenever \(f(y,s)\neq0\), where \(\mathcal V,D,K,H\in C(\Omega^c\times I)\), \(K\geq0\) is bounded, and \(D\) is bounded below. Assume that there exists a positive function \(\Phi\in C^1(I)\) satisfying
\begin{subequations}\label{eq:maximum-principle-assumptions}
\begin{align}
&\|f(\cdot,s_0)\|_{L^\infty}\leq\Phi(s_0),\label{eq:maximum-principle-initial}\\
&|f(y,s)|\leq\Phi(s),\qquad (y,s)\in\Omega\times I,\label{eq:maximum-principle-core}\\
&\Phi'(s)+(D(y,s)-K(y,s))\Phi(s)\geq H(y,s),\qquad (y,s)\in\Omega^c\times I,\label{eq:comparison-barrier}\\
&\limsup_{|y|\to\infty}|f(y,s)|\leq\Phi(s),\qquad s\in I.\label{eq:maximum-principle-tail}
\end{align}
\end{subequations}
Then
\begin{equation}\label{eq:maximum-principle-conclusion}
\|f(\cdot,s)\|_{L^\infty}\leq\Phi(s),\qquad s\in I.
\end{equation}
\end{lemma}

\begin{proof}
Since \(D\) is bounded below and \(K\) is bounded, we can choose \(c>0\) such that \(c+D-K\geq1\) on \(\Omega^c\times I\). For such $c>0$ and any \(\varepsilon>0\), we set
\begin{equation*}
\Phi_\varepsilon(s):=\Phi(s)+\varepsilon e^{c(s-s_0)}.
\end{equation*}
By \eqref{eq:comparison-barrier}, we have
\begin{equation}\label{eq:maximum-principle-strict-barrier}
\Phi_\varepsilon'(s)+(D(y,s)-K(y,s))\Phi_\varepsilon(s) \geq H(y,s)+\varepsilon e^{c(s-s_0)}>H(y,s)
\end{equation}
for \((y,s)\in\Omega^c\times I\).

Suppose to the contrary that \(\|f(\cdot,s)\|_{L^\infty}>\Phi_\varepsilon(s)\) for some \(s\in I\). Then, by the continuity of \(f\) in \(L^\infty\) and \eqref{eq:maximum-principle-initial}, there exists \(s_*\in(s_0,s_1)\) such that
\begin{equation}\label{eq:maximum-principle-first-time}
\|f(\cdot,s_*)\|_{L^\infty}=\Phi_\varepsilon(s_*), \qquad \|f(\cdot,s)\|_{L^\infty}<\Phi_\varepsilon(s) \quad\text{for }s_0\leq s<s_*.
\end{equation}
Moreover, \eqref{eq:maximum-principle-core} and \eqref{eq:maximum-principle-tail} imply that there exists \(y_*\in\Omega^c\) such that
\begin{equation*}
|f(y_*,s_*)|=\|f(\cdot,s_*)\|_{L^\infty} =\Phi_\varepsilon(s_*).
\end{equation*}
Since \(f(y_*,s_*)\neq0\), this implies \(\partial_y|f|(y_*,s_*)=0\).

Set \(Q(s):=\|f(\cdot,s)\|_{L^\infty}-\Phi_\varepsilon(s)\) and \(\mathcal V_*:=\mathcal V(y_*,s_*)\). By \eqref{eq:maximum-principle-first-time}, for sufficiently small \(h>0\), we have
\begin{equation*}
0\leq\frac{Q(s_*-h)-Q(s_*)}{-h} \leq\frac{|f(y_*-h\mathcal V_*,s_*-h)|-\Phi_\varepsilon(s_*-h)-|f(y_*,s_*)|+\Phi_\varepsilon(s_*)}{-h}.
\end{equation*}
Taking \(h\to0^+\) and using \(\partial_y|f|(y_*,s_*)=0\), we obtain
\begin{equation}\label{eq:maximum-principle-time-sign}
\partial_s|f|(y_*,s_*)\geq\Phi_\varepsilon'(s_*).
\end{equation}

On the other hand, applying \eqref{eq:comparison-inequality} at \((y_*,s_*)\) and using \eqref{eq:maximum-principle-strict-barrier},
we obtain
\begin{equation*}
\partial_s|f|(y_*,s_*)
\leq H(y_*,s_*)-(D(y_*,s_*)-K(y_*,s_*))\Phi_\varepsilon(s_*)
<\Phi_\varepsilon'(s_*),
\end{equation*}
which contradicts \eqref{eq:maximum-principle-time-sign}. Hence \(\|f(\cdot,s)\|_{L^\infty}\leq\Phi_\varepsilon(s)\) for all \(s\in I\). Letting \(\varepsilon\downarrow0\), we obtain \eqref{eq:maximum-principle-conclusion}.
\end{proof}

In the applications of Lemma~\ref{lem:bootstrap-comparison} in Section~\ref{sec:proof-common-bootstrap}, we use
\begin{equation}\label{eq:comparison-function}
\Phi(s):=\exp\left(\int_{s_0}^sa(\sigma)\,d\sigma\right)\left(L+\int_{s_0}^se^{-\mu(s-\sigma)}e(\sigma)\,d\sigma\right),\qquad L>0,\quad\mu\geq0,
\end{equation}
where \(a,e\) are nonnegative continuous functions. Direct differentiation gives
\begin{equation}\label{eq:comparison-function-derivative}
\Phi(s)\geq L,\qquad \Phi'(s)+(\mu-a(s))\Phi(s)=\exp\left(\int_{s_0}^sa(\sigma)\,d\sigma\right)(\mu L+e(s))\geq e(s).
\end{equation}
Thus \eqref{eq:comparison-barrier} holds whenever \(D-K\geq\mu-a\) and \(H\leq e\). The bound at spatial infinity is verified separately by Lemma~\ref{lem:transport-far-field} or directly from the bootstrap estimates.

\section{Formal leading-order reductions to the gHS equation}\label{app:formal-reductions}

We present formal leading-order reductions of two model families to the gHS equation \eqref{eq:ghs}. Under the assumptions stated below, these reductions suggest cusp formation with sharp \(C^{(4n-1)/(4n+1)}\) and \(C^{3/5}\) regularity at blow-up for the generalized Camassa--Holm--Novikov family and the asymptotic variational wave equations, respectively. These predictions remain formal: their rigorous justification within our framework requires verifying the hypotheses of Theorem~\ref{thm:common-singularity}, including the bounds on the model-specific terms, for each model.

\subsection{The generalized Camassa--Holm--Novikov family}\label{app:gchn-reduction}

We formally derive \eqref{eq:ghs} with \(\ell=2n\) from the generalized Camassa--Holm--Novikov family \cite{AncoDaSilvaFreire2015, ArmstrongChenHuangShen2026}:
\begin{equation}\label{eq:gchn}
M_t+w^nM_x+(n+1)w^{n-1}w_xM=0,\qquad M=w-w_{xx},\qquad n\in\N.
\end{equation}
We introduce the new unknown
\begin{equation*}
v(x,t)=w(x+c^nt,t)-c,
\end{equation*}
where \(c\neq0\) is a constant background value. Then \eqref{eq:gchn} becomes
\begin{equation}\label{eq:gchn-moving}
v_{xxt}+\bigl((c+v)^n-c^n\bigr)v_{xxx}+(n+1)(c+v)^{n-1}v_xv_{xx}=v_t+\bigl((n+2)(c+v)^n-c^n\bigr)v_x.
\end{equation}
Expanding the coefficients around \(v=0\), we have
\begin{equation*}
(c+v)^n-c^n=nc^{n-1}v+O(v^2),\qquad (c+v)^{n-1}=c^{n-1}+O(v).
\end{equation*}
Using these expansions, we can rewrite \eqref{eq:gchn-moving} as
\begin{equation*}
\partial_x\left(v_{xt}+c^{n-1}\left(nvv_{xx}+\frac12v_x^2\right)\right)=\mathcal{P}_x,
\end{equation*}
where \(\mathcal{P}_x\) collects the terms of higher nonlinear order or lower differential order. Integrating with respect to \(x\) and absorbing the resulting function of time into \(\mathcal{P}\), we obtain
\begin{equation*}
v_{xt}+c^{n-1}\left(nvv_{xx}+\frac12v_x^2\right)=\mathcal{P}.
\end{equation*}
Formally neglecting \(\mathcal{P}\) and applying the normalization \(u=nc^{n-1}v\), we recover \eqref{eq:ghs} with
\begin{equation*}
\ell=2n.
\end{equation*}

\subsection{Asymptotic variational wave equations}\label{app:bzz-reduction}

We formally derive \eqref{eq:ghs} with \(\ell=2\) from the asymptotic variational wave equations studied by Bressan, Zhang, and Zheng \cite{BressanZhangZheng2007}:
\begin{equation}\label{eq:bzz}
(q_t+f(q)_x)_x=\frac12f''(q)q_x^2,
\end{equation}
where \(f\) is a smooth flux. For \(f(q)=q^2/2\) and \(f(q)=q^3/3\), these equations arise as unidirectional weakly nonlinear approximations of the variational wave equation for nematic liquid crystals \cite{HunterSaxton,BressanZheng2006,BressanZhangZheng2007}.

Fix a reference value \(q_*\) with \(f''(q_*)\neq0\), and consider smooth solutions taking values near \(q_*\). We introduce
\begin{equation} \label{zH_def}
z=f'(q),\qquad H=q_t+f'(q)q_x.
\end{equation}
Since \(f''(q_*)\neq0\), the transformation \(q\mapsto z=f'(q)\) is locally invertible near \(q_*\) and preserves local H\"older regularity. By \eqref{eq:bzz} and the definitions \eqref{zH_def},
\begin{equation*}
H_x=\frac12f''(q)q_x^2,\qquad z_t+zz_x=f''(q)H.
\end{equation*}
Differentiating the second identity with respect to \(x\) and using the first, we obtain
\begin{equation}\label{eq:bzz-speed}
z_{xt}+zz_{xx}+\frac12z_x^2=\frac{f'''(q)}{f''(q)}Hz_x.
\end{equation}

Assume that \(H\) remains bounded in the region under consideration. Then the right-hand side of \eqref{eq:bzz-speed} is \(O(|z_x|)\), and is therefore of lower order than \(z_x^2\) as \(|z_x|\to\infty\). Formally neglecting this term, we recover \eqref{eq:ghs} with
\begin{equation*}
\ell=2.
\end{equation*}
For the quadratic flux, \(z=q\) and the right-hand side of \eqref{eq:bzz-speed} vanishes identically. For the cubic flux, \(z=q^2\) and the nondegeneracy condition is \(q_*\neq0\).

\section*{Declaration of AI assistance}

ChatGPT (OpenAI) provided assistance throughout the research and preparation of this manuscript. The authors led the research and directed all aspects of its development and presentation. All mathematical statements, calculations, and proofs were independently verified by the authors, who take full responsibility for the content of the manuscript.

\end{document}